\documentclass[fleqn,12pt]{article}
\usepackage{amssymb,amsmath,amsfonts}
\usepackage{color}
\usepackage{graphicx}
\usepackage{latexsym}
\usepackage{amsmath}
\usepackage{amssymb}
\usepackage{graphics}
\usepackage[dvips]{epsfig}
\usepackage[margin=1in]{geometry}

\numberwithin{equation}{section}
\newtheorem{thm}{Theorem}[section]

\newtheorem{lem}[thm]{Lemma}
\newtheorem{prop}[thm]{Proposition}

\newcommand{\ea}{\epsilon}
\newcommand{\ta}{\theta}

\newcommand{\da}{\delta}

\renewcommand{\aa}{\alpha}

\newcommand{\pl}{\partial}
\newcommand{\sa}{\sigma}

\newcommand{\ga}{\gamma}

\newcommand{\iy}{\infty}

\newcommand{\lt}{\left}
\newcommand{\rt}{\right}
\newcommand{\be}{\begin{equation}}
\newcommand{\bs}{\begin{split}}
\newcommand{\es}{\end{split}}
\newcommand{\ee}{\end{equation}}
\newcommand{\bee}{\begin{equation*}}
\newcommand{\eee}{\end{equation*}}

\newcommand{\ef}{\eqref}
\newcommand{\f}{\frac}
\newcommand{\te}{\tilde}

\begin{document}

\begin{center}
\large{ \bf  Global Solutions to the Gas-Vacuum Interface Problem of Isentropic Compressible Inviscid  Flows with Damping in Spherically Symmetric Motions and Physical Vacuum
}
\end{center}
\centerline{Huihui Zeng}
\begin{abstract}
For the physical vacuum free boundary problem with the sound speed being $C^{{1}/{2}}$-H$\ddot{\rm o}$lder continuous  near vacuum boundaries of the three-dimensional compressible Euler equations with damping,
the  global existence of   spherically symmetric smooth solutions is proved, which are shown to converge to  Barenblatt self-similar solutions of the
porous media equation with the same total masses when   initial data are  small perturbations of  Barenblatt  solutions. The pointwise  convergence with a rate of density, the convergence rate of velocity in supreme norm and the precise expanding rate of  physical vacuum boundaries are also given by constructing nonlinear functionals with space-time  weights featuring the behavior of solutions in large time and near the vacuum boundary and the center of symmetry, the nonlinear energy estimates and elliptic estimates. \end{abstract}

\section{Introduction}
Due to its great physical importance and mathematical challenges, the motion of physical vacuum in compressible fluids has received much attention recently (cf. \cite{7}-\cite{zhenlei1}, \cite{16}-\cite{17'}, \cite{23}-\cite{39}),  and significant progress has been made on the local well-posedness theory (cf. \cite{7}-\cite{10'}, \cite{ 16, 16', LXZ}). However, much less is known on the global existence and long time dynamics of solutions,  which are of fundamental importance in both physics and nonlinear partial differential equations.  This is the main issue we address in this work for the spherically symmetric motions of three-dimensional isentropic compressible inviscid flows with frictional damping.   Physical vacuum problems arise in many physical situations naturally, for example, in the study of the evolution and structure of gaseous stars (cf. \cite{6',cox}) for which vacuum boundaries are natural boundaries.  Another   situation in which the physical vacuum plays an important role is the gas-vacuum interface problem of compressible isentropic Euler equations with damping (cf. \cite{23}-\cite{25}, \cite{38, 39}). In three dimensions, this problem is given as follows:
\be\label{2.1} \begin{split}
&  \rho_t  + {\rm div}(\rho {\bf u}) = 0 &  {\rm in}& \ \ \Omega(t), \\
 &    (\rho {\bf u})_t  + {\rm div}(\rho {\bf u}\otimes {\bf u})+\nabla_{\bf x} p(\rho) = -\rho {\bf u}  & {\rm in}& \ \ \Omega(t),\\
 &\rho>0 &{\rm in }  & \ \ \Omega(t),\\
 & \rho=0    &    {\rm on}& \  \ \Gamma(t):=\pl \Omega(t),\\
 &    \mathcal{V}(\Gamma(t))={\bf u}\cdot {\bf n}, & &\\
&(\rho,{\bf u})=(\rho_0, {\bf u}_0) & {\rm on} & \ \  \Omega:= \Omega(0).
 \end{split} \ee
Here $({\bf x},t)\in \mathbb{R}^3\times [0,\iy)$,  $\rho $, ${\bf u} $, and $p$ denote, respectively, the space and time variable, density, velocity and  pressure; $\Omega(t)\subset \mathbb{R}^3$, $\Gamma(t)$, $\mathcal{V}(\Gamma(t))$ and ${\bf n}$ represent, respectively, the changing volume occupied by a gas at time $t$, moving interface of fluids and vacuum states, normal velocity of $\Gamma(t)$ and exterior unit normal vector to $\Gamma(t)$.  We consider a polytropic gas: the equation of state is  given by
\be\label{polytropic} p(\rho)=\rho^{\gamma} \  \ {\rm for} \ \  \gamma>1 \ee
with  the adiabatic constant   set to be unity for convenience.  Equations $\ef{2.1}_{1,2}$ describe the balance laws of mass and momentum, respectively; conditions $\ef{2.1}_{3,4}$ state that $\Gamma(t)$ is the interface to be investigated; $\ef{2.1}_5$ indicates that the interface moves with the normal component of the fluid velocity; and $\ef{2.1}_6$ is the initial conditions for the density, velocity and domain.   Let $c=\sqrt{ p'(\rho)}=\sqrt{\gamma \rho^{ {\ga-1}  } }$ be the sound speed, the condition
\be\label{physical vacuum} -\infty<\nabla_{\bf n}(c^2)<0  \  \ {\rm on} \ \  \Gamma(t) \ee
defines a physical boundary (cf. \cite{7}, \cite{10'}, \cite{16'}, \cite{23}-\cite{25}).

The compressible Euler equations of isentropic flow  with damping are closely related to the
porous media equation (cf. \cite{HL, HMP, HPW, 23}):
\begin{equation}\label{pm}
 \rho_t =\Delta p(\rho),
\end{equation}
when $ \ef{2.1} _2$ is simplified to Darcy's law:
\begin{equation}\label{darcy}
\nabla_{\bf x} p(\rho)=- \rho {\bf u}.
\end{equation}
For \ef{pm}, basic understanding of the solution with finite mass $M>0$
is provided by Barenblatt (cf. \cite{ba}), which is spherically symmetric and given by
\be\label{1.6}
\bar\rho({\bf x} ,t)= \bar\rho(r,t)=(1+t)^{-\frac{3}{3\ga-1}}\lt( A - B ({1+ t})^{-\frac{2}{3\ga-1}}   r^2 \rt)^{\frac{1}{\ga-1}} \ \  {\rm with} \ \  r=|{\bf x}|,
   \ee
where 
\be\label{constAB}
B= \f{\ga-1}{2\ga(3\ga-1)} \ \ {\rm and} \ \  (\ga A)^{\frac{3\ga-1}{2(\ga-1)}}= M\ga^{\frac{1}{\ga-1}} (\ga B)^{\frac{3}{2}}
   \lt( \int_0^1 y^2 \lt(1-y^2\rt)^{\frac{1}{\ga-1}}dy \rt)^{-1}.
   \ee
Clearly,  
\be\label{barenblattbdry}
\int_0^{\bar R(t)} r^2 \bar\rho(r,t) dr=M  \ \ {\rm for} \ \ t\ge 0, \ \  {\rm where} \ \   \bar R(t)=
\sqrt{A/B}(1+t)^{ {1}/({3\ga-1})}.
\ee 
The corresponding  Barenblatt velocity $\bar {\bf u}$ is defined  by
\begin{align*}
\bar {\bf u}({\bf x}, t)= ( {{\bf x}}/ r ) \bar u(r,t)   \  \     \ \  \textrm{ in the region} \ \ \{(r,t): 0\le r\le \bar R(t), \ t >0\} ,
\end{align*}
where
\begin{align*}
\bar u(r,t)=- \f{  p\lt(\bar\rho\rt)_r}{\bar\rho}=    \f{r}{(3\ga-1)(1+t)}
\ \  {\rm satisfying} \ \ \bar u(0,t)=0 \ \ {\rm and} \ \ \dot{\bar R}(t)=\bar u \lt(\bar R(t), t\rt).
\end{align*}
So, $(\bar\rho, \bar {\bf u})$ defined in the region $\{(r,t): 0\le r\le \bar R(t), \ t >0\}$
 solves \ef{pm} and \ef{darcy}.

The vacuum boundary $r=\bar R(t)$  of Barenblatt's solution is physical. This is the major motivation to study  the physical vacuum free boundary problem of compressible Euler equations with damping. To this end, a class of explicit spherically symmetric  solutions to problem \ef{2.1} was constructed  in \cite{23}, which are of the following form: 
\begin{align}\label{liuexplicitsolution}  
 \Omega(t) =B_{R(t)}({\bf 0}) , \ \ 
 c^2({\bf x}, t)= c^2(r, t)={e(t)-b(t)r^2} \ \  {\rm and} \ \  {\bf u}({\bf x}, t)= ({{\bf x}}/{r} )u(r, t),
  \end{align}
where
$$ R(t)= \sqrt {e(t)/b(t)}   \ \   {\rm and} \ \    u(r, t)=a(t) r .$$
In \cite{23},  a system of ordinary differential equations for $(e,b,a)(t) $ was derived with $e(t), b(t)>0$ for $t\ge 0$ by substituting \ef {liuexplicitsolution} into $\ef{2.1}_{1, 2}$ and the time-asymptotic equivalence of this explicit solution and  Barenblatt's  solution with the same total mass was shown. Indeed, the Barenblatt solution  of \ef{pm} and \ef{darcy} can be obtained by the same ansatz as \ef{liuexplicitsolution}:
$$ 
 \bar c^2({\bf x}, t)=\bar e(t)-\bar b(t)r^2  \ \  {\rm and} \ \  {\bf u}({\bf x}, t)= \bar a(t) {{\bf x}}  . $$
 Substituting this into \ef{pm}, \ef{darcy} and \ef{barenblattbdry} with $\bar R (t)= \sqrt{{\bar e(t)}/{\bar b(t)}}$   gives
 $$
   \bar e(t)= \ga A (1+  t)^{-3({\gamma-1})/({3\ga-1})},\   \bar b(t)=\ga B(1+ t)^{-1}
  \  \ {\rm and} \  \ \bar a(t)=  {(3\ga-1)^{-1}(1+t)^{-1}},$$
  where $A$ and $B$ are determined by \ef{constAB}. Precisely, it was proved in \cite{23} the following  time-asymptotic equivalence:
  $$
  (a,\ b,\ e)(t)=(\bar a, \ \bar b, \ \bar e)(t)+ O(1)(1+t)^{-1}{\ln (1+t)} \ \ {\rm as}\  \ t\to\infty.$$

A question was raised in \cite{23} whether this equivalence is still true for general solutions to  problem \ef{2.1}.  For this purpose,   we seek  solutions to  problem \ef{2.1} of the form:
$$\Omega(t) =B_{R(t)}({\bf 0}), \ \  \rho({\bf x}, t)=\rho(r, t), \ \  {\bf u}({\bf x}, t)= ({{\bf x}}/{r}) u(r, t)  \ \   {\rm with} \ \  r=|{\bf x}|.$$ 
Then  problem \ef{2.1} reduces to
\begin{equation}\label{equation1}
\begin{split}
&    (r^2\rho)_t+ (r^2\rho u)_r=0    & {\rm in } & \  \  \lt(0, \  R(t)\rt), \\
&\rho( u_t +u  u_r)+ p_r=- \rho u  & {\rm in } & \  \  \lt(0, \  R(t)\rt),  \\
&\rho > 0  & {\rm in} &   \  \  \lt[0, \  R(t)\rt),\\
& \rho\lt(R(t),t\rt) =0,      \ \ u(0,t)=0,     \\
&\dot R(t)=u(R(t), t) \ \ {\rm with} \ \ R(0)=R_0,   \\
& (\rho,u)(r,t=0)=\lt(\rho_0, u_0\rt)(r) & {\rm  on } & \  \  \lt(0, \  R_0\rt),
\end{split}
\end{equation}
so that $R(t)$ is the  radius of the domain occupied by the gas at time $t$ and $r=R(t)$  represents  the vacuum free boundary which issues from $r=R_0$ and moves with the fluid velocity.  One of motivations to study the spherically symmetric solution  is that  the Barenblatt solution posses the same symmetry, and it is expected that spherically symmetric solutions will provide insights on the local and
long time behavior of solutions to the general three-dimensional problem \ef{2.1}. Locally,  at each point ${\bf x}$ in $\Omega(t)$, it might be plausible to rotate a solution in all possible ways about ${\bf x}$ and average all rotations in the spirit of spherical mean.
In long time, for a general three-dimensional problem, it is expected that the geometry of the boundary becomes more and more symmetric
due to the dissipation of damping which dissipates the total energy.

In the spherically symmetric setting, the physical vacuum boundary condition \ef{physical vacuum} reduces to
$$
 -\infty< ({  c^2})_r <0 $$
in a small neighborhood of the boundary.  To capture this  singularity, the initial domain is taken to be a  ball  $\{0\le r\le R_0\}$ and the initial density  is assumed to satisfy
\be\label{156}
\rho_0(r)>0 \ \ {\rm for} \ \  0\le r<R_0 , \ \ \rho_0(R_0)=0 \ \ {\rm and}  \ \
  -\iy<   \lt(\rho_0^{\ga-1}\rt)_r <0 \  \ {\rm at} \  \ r=R_0.
   \ee
We   require that the initial total mass is the same as that of the Barenblatt solution, that is,
\be\label{156'} \int_{0}^{R_0} r^2 \rho_0(r)dr =\int_{0}^{\bar{R}(0)} r^2 \bar \rho_0(r)dr= M.\ee
The conservation law of mass, $\ef{equation1}_1$, and  \ef{barenblattbdry}  give 
$$\int_{0}^{R(t)} r^2 \rho(r,t)dr = \int_{0}^{R_0} r^2 \rho_0(r)dr = M=\int_0^{\bar R(t)} r^2 \bar\rho(r,t) dr \ \ {\rm for} \ \ t\ge 0. $$
In the present work, we prove the global existence of smooth solutions to problem \ef{equation1} when initial data  are small spherically symmetric perturbations of Barenblatt  solutions and they have the same total masses. Moreover, we obtain the pointwise  convergence with a rate of density which gives the detailed behavior of the density, the convergence rate of velocity in supreme norm and the precise expanding rate of  physical vacuum boundaries.  The results obtained in this article  also prove the nonlinear asymptotic stability of  Barenblatt   solutions in the setting of physical vacuum free boundary problems.

The physical vacuum  with  the sound speed being $C^{ {1}/{2}}$-H$\ddot{\rm o}$lder continuous  across  vacuum boundaries  makes it
 challenging and interesting in the study of  free boundary problems in compressible fluids, even for  the local-in-time existence theory,  since  standard methods of symmetric hyperbolic systems (cf. \cite{17}) do not apply.
Indeed, characteristic speeds of the compressible isentropic Euler equations become singular with infinite spatial derivatives
at  vacuum boundaries which creates much severe  difficulties in analyzing the regularity near boundaries.
Recently, important progress has been made in the local-in-time well-posedness theory for the  compressible
Euler equations (cf. \cite{7}-\cite{10'}, \cite{16,16'}). On the other hand,  it  poses a great challenge to extend the local-in-time existence theory to the global one of smooth solutions,  due to the strong degeneracy near vacuum states caused by the singular behavior of  physical vacuum.    Obtaining  the
global-in-time regularity of solutions near  vacuum boundaries by establishing the uniform-in-time higher-order estimates is the key to  analyses. This is nontrivial due to the strong degenerate nonlinear hyperbolic nature. 
To obtain global-in-time  estimates, it is essential to show decay estimates, which are achieved in the present work by  introducing  time weights to quantify the long time behavior of solutions. This is in sharp contrast to the weighted estimates used in establishing the local-in-time well-posedness theory (cf. \cite{7}-\cite{10'}, \cite{16,  16', LXZ}), where only spatial weights are involved.

It should be noted that,  as the first step to understand  global solutions and their long time behavior for physical vacuum boundary problems of the isentropic compressible Euler equations with frictional damping, Luo and the author (cf. \cite{LZ}) proved  the global smooth solutions and convergence to  Barenblatt solutions as time goes to infinity in one-dimensional case,  based on a construction of  higher-order  weighted functionals with both space and time
weights capturing the behavior of solutions both near vacuum states and in large time,  an introduction of a new ansatz, higher-order nonlinear energy estimates and elliptic estimates. In general, much more obstacles appear in the study of   multi-dimensional problems of   compressible Euler equations as a prototype.   Compared with the one-dimensional case studied in \cite{LZ}, it is much more difficult and involved to solve
the three-dimensional spherically symmetric problem, \ef{equation1},   
 in  the construction of the nonlinear weighted functionals,  nonlinear weighted estimates and   elliptic estimates.  
Besides the difficulty of degeneracy of the equations at vacuum states, one of the difficulties in solving \ef{equation1} is the coordinates singularity at the origin, the center of symmetry, which carries the true three-dimensional nature. We succeed in constructing  suitable weights to resolve the coordinates singularity in this paper. As an intermediate step passing from one-dimensional case in \cite{LZ} to the general three-dimensional problem, \ef{2.1}, we believe the
ideas and  estimates including the nonlinear weighted functionals and  pointwise decay estimates developed in this paper will
contribute to a understanding of the behavior of solutions to  problem  \ef{2.1}.

To the best of our knowledge, the results obtained in this paper are among the first ones on the global existence of smooth solutions for   physical vacuum free boundary problems in  inviscid compressible fluids.
It should be pointed that the $L^p$-convergence of $L^{\infty}$-weak solutions for the Cauchy problem of the one-dimensional compressible Euler equations with damping  to Barenblatt  solutions of the porous media equations was given in \cite{HMP} with  $p=2$ if $1<\ga\le 2$ and $p=\ga$ if $\ga>2$ and in \cite{HPW} with $p=1$, respectively, using   entropy-type estimates for the solution itself without deriving  estimates for derivatives.
However, the interfaces separating gases and vacuum cannot be traced in the framework of $L^{\infty}$-weak solutions.  The aim of this work is to understand the behavior and long time dynamics of   physical vacuum boundaries, for which obtaining the global-in-time regularity of solutions is essential.

\section{Reformulation of the problem and main results}
\subsection{Fix the domain and Lagrangian variables}
We make the initial domain of the Barenblatt solution, $\lt(0, \bar R (0)\rt)$, as the reference domain and define a diffeomorphism $\eta_0: \lt(0, \bar R (0)\rt) \to \lt(0, R_0\rt)  $ by
\bee
\int_0^{\eta_0(r)} r^2 \rho_0(r) dr = \int_0^{r} r^2 \bar\rho_0(r)dr \ \  {\rm for} \ \
r\in  \lt(0, \bar R (0)\rt),
\eee
where
$\bar\rho_0(r): = \bar\rho(r,0) $ is the initial density of the Barenblatt solution. Clearly,
\be\label{lagrangiandensity}
 \eta_0^2(r)  \rho_0(\eta_0(r)) \eta_{0r}(r) = r^2 \bar\rho_0(r)    \ \  {\rm for} \ \
 r\in  \lt(0, \bar R (0)\rt) .
 \ee
Due to \ef{156}, \ef{1.6} and the fact that the total mass of the Barenblatt solution is the same as that of $\rho_0$, \ef{156'}, the diffeomorphism $\eta_0$ is well defined. For simplicity of presentation,  set
$$ \mathcal{I} : = \lt(0, \ \bar R(0)\rt)=\lt(0, \ \sqrt{A /B} \rt). $$
To fix the boundary, we transform system \ef{equation1} into Lagrangian variables.  For $r\in \mathcal{I}$, we define the Lagrangian variable  $\eta(r, t)$ by
\be\label{haz2.2}
  \eta_t(r, t)= u(\eta(r, t), t) \ \ {\rm for} \  \ t>0 \  \ {\rm and} \ \  \eta(r, 0)=\eta_0(r) ,
\ee
and set the Lagrangian density and velocity by
\be\label{o2.1}
f(r, t)=\rho(\eta(r, t), t) \ \ {\rm and} \ \  v(r, t)=u(\eta(r,t), t).
\ee
Then  the Lagrangian version of system \ef{equation1} can be written on the reference domain $\mathcal{I}$ as
\be\label{new419} \begin{split}
& (\eta^2 f)_t +r^2f {v_r}/{\eta_r}=0    & {\rm in}& \  \    \mathcal{I}\times (0, \iy),\\
& f  v_t+  { (f^{\gamma})_r}/{  \eta_r}=- f v\ \  &{\rm in}& \  \   \mathcal{I}\times (0, \iy), \\
& v(0,t)=0 & {\rm on} & \ \  (0,\iy),\\
 & (f, v) =\lt(\rho_0(\eta_0), u_0(\eta_0)\rt)  &  {\rm on}& \  \   \mathcal{I} \times \{t=0\}.
\end{split}
\ee
It should be noted that we need $\eta_r(r,t)>0$ for $r\in \mathcal{I}$ and $t\ge 0$ to make the Lagrangian transformation sensible, which will be verified in \ef{ht3.3}.  Indeed, $\eta_r>0$ implies $\eta(r,t)>0$ for $r\in \mathcal{I}$ and $t\ge 0$, due to the boundary condition that the center of the symmetry does not move, $v(0,t)=0$. The map $\eta(\cdot, t)$ defined above can be extended to $\bar {\mathcal{I}}= [0, \ \sqrt{A/ B }  ]$. In the setting, the vacuum free boundaries for problem \ef{equation1} are given by
\be\label{vbs}  R(t)=\eta\lt(\bar R(0), \  t\rt)= \eta\lt(\sqrt{A/B}, \ t\rt) \ \  {\rm for} \ \  t\ge 0.\ee
 It follows from solving $\eqref{new419}_1$ and using \ef{lagrangiandensity} that
\be\label{newlagrangiandensity}
f(r,t)\eta^2(r,t) \eta_r(r,t)= \rho_0(\eta_0(r)) \eta_0^2(r) \eta_{0r}(r) =r^2 \bar\rho_0(r), \ \  r\in \mathcal{I}  .\ee
So, the initial density of the Barenblatt solution, $\bar\rho_0$, can be viewed as a parameter and  system \ef{new419} can be rewritten as
\be\label{419} \begin{split}
& \bar\rho_0\eta_{tt} + \bar\rho_0 \eta_t  + \lt( \frac{\eta}{r}\rt)^2  \left[    \lt(\frac{r^2}{\eta^2}\frac{\bar\rho_0}{ \eta_r}\rt)^\ga    \right]_r =0   \ \  &{\rm in}& \  \    \mathcal{I}\times (0, \iy), \\
&  \eta(0, t)=0,    &  {\rm on}& \  \ (0,\iy),\\
& (\eta, \eta_t) = \lt(\eta_0, u_0(\eta_0)\rt)  &  {\rm on}& \  \    \mathcal{I}\times (0, \iy).
\end{split}
\ee

\subsection{Ansatz}
Define the Lagrangian variable $\bar\eta(r, t)$ for the Barenblatt flow in $\bar {\mathcal{I}}$ by
\be\label{haz2.7}
 \pl_t \bar \eta(r, t)= \bar u(\bar \eta (r, t), t)=\f{  \bar\eta(r,t)}{(3\ga-1)(1+ t)}  \ \ {\rm for} \  \ t>0 \  \ {\rm and} \ \  \bar \eta(r, 0)=r,
\ee
so that
\be\label{bareta}
\bar \eta(r,t)= r \lt(1 +  t \rt)^{{1}/({3\ga-1})}   \ \  {\rm for} \ \  (r,t)\in \bar{\mathcal{I}}\times [0,\iy)
\ee
and
\bee\label{barequ}\lt.\begin{split}
&    \bar\rho_0 \bar \eta_t + \lt(\frac{\bar\eta}{r}\rt)^2  \left[    \lt(\frac{r^2}{\bar \eta^2}\frac{\bar\rho_0}{\bar \eta_r}\rt)^\ga    \right]_r =0  \ \  {\rm in} \  \  {\mathcal{I}}\times (0,\iy).
\end{split}\rt.\eee
Since $\bar\eta$ does not solve $\ef{419}_1$ exactly, we introduce a correction $h(t)$ which is a solution of the following initial value problem of ordinary differential equations:
\be\label{pomt}\begin{split}
& h_{tt} + h_t -   (\bar \eta_r+h)^{2-3\ga} /(3\ga-1)  + \bar\eta_{rtt}  +\bar\eta_{rt}    =0, \ \ t >0,  \\
& h(t=0)=h_t(t=0)=0.
 \end{split}\ee
 (Notice that $\bar \eta_r$,  $\bar \eta_{rt}$ and $\bar \eta_{rtt}$ are independent of $r$.)
The new ansatz is then given by
\be\label{teeta}
\tilde \eta(r, t):=\bar \eta (r,t)  + r h(t)
,
\ee
so that
 \be\label{equeta} \begin{split}
&  \bar\rho_0  \tilde\eta_{tt}   +  \bar\rho_0  \tilde\eta_t + \lt(\frac{\tilde\eta}{r}\rt)^2 \left[    \lt(\frac{r^2}{\tilde\eta^2}\frac{\bar\rho_0}{ \tilde\eta_r}\rt)^\ga    \right]_r =0    \ \  {\rm in} \  \  {\mathcal{I}}\times (0,\iy).
\end{split}
\ee
It should be noted that $\tilde\eta_r$ is independent of $r$. We will prove in the Appendix that  $\tilde\eta$ behaves similar to $\bar\eta$, that is, there exist positive constants
$K$ and $C(n)$ independent of time $t$ such that for all $t\ge 0$,
\be\label{decay}\begin{split}
&\lt(1 +   t \rt)^{ {1}/({3\ga-1})} \le \tilde \eta_{r}(t) \le K \lt(1 +   t \rt)^{ {1}/({3\ga-1})},  \ \   \ \   \tilde\eta_{rt} \ge 0 ,  \\
&\lt|\f{d^k\tilde \eta_{r}(t)}{dt^k}\rt| \le C(n)\lt(1 +   t \rt)^{\frac{1}{3\ga-1}-k},   \ \ k=1, 2,  \cdots, n.
 \end{split}\ee
Moreover, there exists a certain constant $C$ independent of $t$ such that
 \be\label{decayforh}\begin{split}
 0\le h(t) \le C (1+t)^{\frac{1}{3\ga-1}-1}\ln(1+t) \ \  {\rm and} \ \ |h_t(t)| \le C (1+t)^{\frac{1}{3\ga-1}-2}\ln(1+t), \ \  t\ge 0.
 \end{split}\ee
The proof of \ef{decayforh} will also be given in the Appendix.
\subsection{Main results}
To state the main theorem, we write equation $\ef{419}_1$ in a perturbation form around the Barenblatt solution. Let
$$\zeta(r,t):= \eta(r,t)/r- {\tilde\eta(r,t)}/r .$$
Thus,
\be\label{hh2-12}
\eta(r,t) =  {\tilde \eta }(r,t) + r\zeta(r,t)    \ \  {\rm and} \ \  \eta_r (r,t)=  {\tilde \eta }_r(t) + \zeta(r,t) + r\zeta_r(r,t).
\ee
It follows from $\ef{419}_1$ and \ef{equeta} that
\be\label{pertb}\begin{split}
   &r \bar\rho_0 \zeta_{tt}   +
   r  \bar\rho_0  \zeta_t
   +   \lt(\tilde\eta_r+ \zeta\rt)^2\left[ \bar\rho_0^\ga     \lt(\tilde\eta_r+\zeta\rt)^{-2\ga}\lt(\tilde\eta_r+\zeta +r\zeta_r\rt)^{-\ga}   \right]_r -    \tilde\eta_r^{2-3\ga}  \lt(   {\bar\rho_0}  ^\ga    \right)_r
   =0 ,
\end{split}\ee
Denote
$$ \alpha := 1/(\ga-1), \  \ l:=3+ \min\lt\{m \in  \mathbb{N}: \ \  m> \alpha \rt\}=4 +[\aa] .$$
For $j=0,\cdots, l$ and  $i=0,\cdots, l-j$,  we set
\bee\label{}\begin{split}
&\mathcal{ E}_{j}(t) : = (1+ t)^{2j} \int_\mathcal{I}   \lt[r^4 \bar\rho_0  \lt(\pl_t^j \zeta\rt)^2 +r^2 \bar\rho_0^\ga  \lt|\pl_t^j \lt( \zeta, r\zeta_r \rt)  \rt|^2  + (1+ t)  r^4 \bar\rho_0    \lt(\pl_t^{j} \zeta_t\rt)^2 \rt]   dr,\\
&\mathcal{ E}_{j, i}(t): =  (1+  t)^{2j}  \int_\mathcal{I} \lt[r^2 \bar\rho_0^{1+(i-1)(\ga-1)}  \lt(\pl_t^j \pl_r^i \zeta\rt)^2 + r^4 \bar\rho_0^{1+(i+1)(\ga-1)}  \lt(\pl_t^j \pl_r^{i+1} \zeta\rt)^2 \rt] dr  .
\end{split}\eee
 The higher-order norm is defined by
$$
  \mathcal{E}(t) :=  \sum_{j=0}^l \lt(\mathcal{ E}_{j}(t) + \sum_{i=1}^{l-j} \mathcal{ E}_{j, i}(t)  \rt).$$
  It will be proved in Lemma \ref{lem34} that
$$
  \sup_{r\in \mathcal{I}} \lt\{\sum_{j=0}^2 (1+  t)^{2j} \lt |\pl_t^j \zeta(r, t)\rt |^2   +  \sum_{j=0}^1 (1+  t)^{2j} \lt |\pl_t^j \zeta_r(r, t)\rt|^2\rt\} \le C  \mathcal{E}(t)
$$
for some constant $C$ independent of $t$. So the boundedness of $\mathcal{E}(t)$ gives the uniform boundedness and decay of the perturbation $\zeta$ and its derivatives. In what follows, we state our main result.

{\begin{thm}\label{thm1} Suppose that \ef{156'} holds. There exists a  constant $\bar \delta >0$ such that if
$\mathcal{E}(0)\le \bar \da,$
then   the  problem \eqref{419}  admits a global unique smooth solution  in $\mathcal{I}\times[0, \iy)$ satisfying for all $t\ge 0$,
$$\mathcal{E}(t)\le C\mathcal{E}(0) $$
and
\be\label{thm1est2}\begin{split}
  \sup_{r\in \mathcal{I}} &  \lt\{   \sum_{j=0}^2 (1+  t)^{2j} \lt |\pl_t^j \zeta(r, t)\rt |^2   +  \sum_{j=0}^1 (1+  t)^{2j} \lt |\pl_t^j \zeta_r(r, t)\rt|^2
      +
  \sum_{
  i+j\le l-2,\   2i+j \ge 3} (1+  t)^{2j} \rt. \\
 & \lt. \times\lt |  \bar\rho_0^{\f{(\ga-1)(2i+j-3)}{2}}\pl_t^j \pl_r^i \zeta(r, t)\rt|^2  + \sum_{
  i+j=l-1 } (1+  t)^{2j}\lt | r \bar\rho_0^{\f{(\ga-1)(2i+j-3)}{2}}\pl_t^j \pl_r^i \zeta(r, t)\rt|^2  \rt.\\
&\lt.   + \sum_{
  i+j=l  } (1+  t)^{2j}\lt | r^2 \bar\rho_0^{\f{(\ga-1)(2i+j-3)}{2}}\pl_t^j \pl_r^i \zeta(r, t)\rt|^2 \rt\} \le C \mathcal{E}(0),
\end{split}\ee
where $C$ is a positive constant  independent of $t$.
\end{thm}
It should be noticed that  the time derivatives involved in the initial higher-order energy norm,
$\mathcal{E}(0)$, can be determined via the equation by the initial data $\rho_0$ and $u_0$  (see \cite{10,LXZ} for instance).

As a corollary of Theorem \ref{thm1}, we have the following theorem for solutions to the original vacuum free boundary problem
\ef{equation1}.
\begin{thm}\label{thm2} Suppose that \ef{156'} holds. There exists a  constant $\bar\delta >0$ such that if
 $\mathcal{E}(0)\le \bar\da, $
then  the  problem
\ef{equation1}  admits a global unique smooth solution $\lt(\rho(\eta,t), u(\eta,t),  R(t)\rt)$ for $t\in[0,\iy)$  satisfying
 \be\label{thm2est1}\begin{split}
 \lt|\rho\lt(\eta(r, t),t\rt)-\bar\rho\lt(\bar\eta(r, t), t\rt)\rt|
 \le & C\lt(A-B r ^2\rt)^{\frac{1}{\ga-1}}(1+t)^{-\frac{4}{3\ga-1}} \\
 & \times \lt(\sqrt{\mathcal{E}(0)}+(1+t)^{-\frac{3\ga-2}{3\ga-1}}\ln(1+t)\rt), \end{split}
 \ee
\be\label{thm2est2}
 \lt|u\lt(\eta(r, t),t\rt)-\bar u\lt(\bar\eta(r, t), t\rt)\rt|
\le C r (1+t)^{-1}\lt( \sqrt{\mathcal{E}(0)}+(1+t)^{-\frac{3\ga-2}{3\ga-1}}\ln(1+t)\rt),
\ee
for all $r \in \mathcal{I}$ and $t\ge 0 $; and for all $t\ge 0$,
\be\label{thm2est3}  c_1(1+t)^{\frac{1}{3\ga-1}}\le R(t)\le c_2 (1+t)^{\frac{1}{3\ga-1}} , \ee
\be\label{thm2est4} \lt|\frac{d^k R(t)}{dt^k}\rt|\le C(1+t)^{\frac{1}{3\ga-1}-k} , \ \   k=1, 2 ,3 ,\ee
\be\label{thm2est5} c_3 (1+t)^{-\frac{3\ga-2}{3\ga-1}}\le \lt|    \lt(\rho^{\ga-1}\rt)_\eta (\eta,t)  \rt| \le c_4(1+t)^{-\frac{3\ga-2}{3\ga-1}}  \ \ {\rm when} \  \   \frac{1}{2}R(t)\le \eta \le R(t).\ee
Here $C$, $c_1$, $c_2$, $c_3$ and $c_4$ are    positive  constants  independent of  $t$.
\end{thm}
The pointwise behavior of the density and  velocity for the vacuum free boundary problem
\ef{equation1}  to that of the Barenblatt  solution are given by  \ef{thm2est1} and \ef{thm2est2}, respectively. It is also shown in \ef{thm2est1} that
the difference of density  to problem \ef{equation1} and the corresponding Barenblatt density decays at the
rate of $(1+t)^{-{4}/(\ga+1)}$ in $L^\iy$, while the density of  the Barenblatt solution, $\bar\rho$, decays at the
rate of $(1 + t)^{-3/(\ga+1)}$ in $L^\iy$ (see \ef{1.6}).
\ef{thm2est3} gives the precise expanding rate
of the vacuum boundaries of the problem \ef{equation1}, which is the same as that for the Barenblatt solution shown in \ef{barenblattbdry}. Furthermore, it  verifies in \ef{thm2est5} that the vacuum boundary $R(t)$ is physical at any finite time.

\section{Proof of Theorem \ref{thm1}}
The proof  is based on the local existence of smooth solutions  (cf. \cite{LXZ,10, 16}) and continuation arguments. The uniqueness of the smooth solutions can be obtained as in section 11 of \cite{LXZ}.  In order to prove the global existence
 of  smooth solutions, we need to obtain the uniform-in-time {\it a priori} estimates  on any given time interval $[0, T]$ satisfying $\sup_{t\in [0, T]}\mathcal{E}(t)<\infty$. For this purpose, we use a bootstrap argument by making the following  {\it a priori} assumption: Let $\zeta$ be a smooth solution to \ef{pertb} on $[0 , T]$, there exists a suitably small fixed positive number $\ea_0\in (0,1)$ independent of $t$   such that for
 $t \in [ 0, T]$,
\be\label{apb}\begin{split}
&  \sum_{j=0}^2 (1+  t)^{2j} \lt \|\pl_t^j \zeta (\cdot,t) \rt \|_{L^\iy}^2   +  \sum_{j=0}^1 (1+  t)^{2j} \lt \|\pl_t^j \zeta_r (\cdot,t) \rt\|_{L^\iy}^2
     +
  \sum_{
  i+j\le l-2,\   2i+j \ge 3} (1+  t)^{2j}   \\
    & \times   \lt \|  \bar\rho_0^{\f{(\ga-1) (2i+j-3) }{2}}\pl_t^j \pl_r^i \zeta(\cdot,t) \rt\|_{L^\iy}^2 +
  \sum_{
  i+j=l-1} (1+  t)^{2j}   \lt \| r \bar\rho_0^{\f{(\ga-1) (2i+j-3) }{2}} \pl_t^j \pl_r^i \zeta  (\cdot,t) \rt\|_{L^\iy}^2    \le \epsilon_0^2 .
\end{split}\ee
This in particular implies,  noting  $\ef{decay}$, that for $0\le \ta_1, \ta_2\le 1$,
\be\label{basic}
\frac{1}{2}(1+t)^{\frac{1}{3\ga-1 }}\le (\tilde \eta_r+ \ta_1 \zeta + \ta_2 r\zeta_r )(r, t)\le 2K(1+t)^{\frac{1}{3\ga-1 }}, \ \ (r, t)\in \mathcal{I}\times [0, T].\ee
Moreover, it follows from \ef{hh2-12} and \ef{basic} that
\be\label{ht3.3}
\frac{1}{2}(1+t)^{\frac{1}{3\ga-1 }}\le \eta_r(r, t), \  r^{-1}\eta(r,t)\le 2K(1+t)^{\frac{1}{3\ga-1 }}, \ \ (r, t)\in \mathcal{I}\times [0, T].\ee
Here $K$ is the positive constant appearing in $\ef{decay}_1$.

Under this {\it a priori} assumption, we prove in section 3.2 the following elliptic estimates:
$$\mathcal{ E}_{j, i}(t) \le C   \sum_{\iota=0}^{i+j}\mathcal{ E}_{\iota}(t) ,    \ \ {\rm when}      \ \  j \ge 0,  \ \  i \ge 1, \ \  i+j\le l, $$
where $C$ is a positive constant independent of $t$.
With the {\it a priori} assumption and elliptic estimates,  we show in section 3.3 the following nonlinear weighted energy estimate: for some  positive constant $C$ independent of $t$,
$$\mathcal{E}_j(t)     \le C \sum_{\iota=0}^j \mathcal{ E}_{\iota}(0), \ \  j=0,1, \cdots, l. $$
Finally, the {\it a priori} assumption \ef{apb} can be verified in section 3.4 by proving
\begin{align*}\label{}
  &  \sum_{j=0}^2 (1+  t)^{2j} \lt \|\pl_t^j \zeta (\cdot,t) \rt \|_{L^\iy}^2   +  \sum_{j=0}^1 (1+  t)^{2j} \lt \|\pl_t^j \zeta_r (\cdot,t) \rt\|_{L^\iy}^2
     +
  \sum_{
  i+j\le l-2,\   2i+j \ge 3} (1+  t)^{2j}   \notag\\
    & \times   \lt \|  \bar\rho_0^{\f{(\ga-1) (2i+j-3) }{2}}\pl_t^j \pl_r^i \zeta(\cdot,t) \rt\|_{L^\iy}^2 +
  \sum_{
  i+j=l-1} (1+  t)^{2j}   \lt \| r \bar\rho_0^{\f{(\ga-1) (2i+j-3) }{2}} \pl_t^j \pl_r^i \zeta  (\cdot,t) \rt\|_{L^\iy}^2
   \notag \\
 & + \sum_{
  i+j=l} (1+  t)^{2j}   \lt \| r^2  \bar\rho_0^{\f{(\ga-1) (2i+j-3) }{2}} \pl_t^j \pl_r^i \zeta (\cdot,t) \rt\|_{L^\iy}^2 \le C\mathcal{E}(t).
\end{align*}
for some positive constant $C$  independent of $t$.
This closes the whole bootstrap argument for  small initial perturbations  and completes the proof of Theorem \ref{thm1}.

\subsection{Preliminaries}
In this subsection, we  list some embedding estimates for weighted Sobolev spaces which will be used later  and introduce some notations to simplify the presentation.
For any bounded interval $I$, set
$d(r)=dist(r, \partial I).$
  For any  $a>0$ and nonnegative integer $b$, the  weighted Sobolev space  $H^{a, b}( {I})$ is given by
$$ H^{a, b}({I}) := \lt\{   d^{a/2}F\in L^2( {I}): \ \  \int_ {I}     d^a|\pl_r ^k F|^2dr<\infty, \ \  0\le k\le b\rt\}$$
  with the norm
$$ \|F\|^2_{H^{a, b}({I})} := \sum_{k=0}^b \int_ {I}    d^a|\pl_r^k F|^2dr.$$
Then  for $b\ge  {a}/{2}$, it holds the following {\it embedding of weighted Sobolev spaces} (cf. \cite{18'}):
 $$ H^{a, b}( {I})\hookrightarrow H^{b- {a}/{2}}( {I})$$
    with the estimate
  \be\label{wsv} \|F\|_{H^{b- {a}/{2}}( {I})} \le C  \|F\|_{H^{a, b}( {I})} \ee
for some positive constant $C$ depending on $a$, $b$ and $I$.

The following general version of the {\it Hardy inequality} whose proof can be found  in \cite{18'} will also be used often in this paper.
 Let $k>1$ be a given real number and $F$ be a function satisfying
$$
\int_0^{\da} r^k\lt(F^2 + F_r^2\rt) dr < \iy,
$$
where $\da$ is a positive constant; then it holds that
$$
\int_0^{\da} r^{k-2} F^2 dr \le C(\delta, k) \int_0^{\da} r^k \lt( F^2 + F_r^2 \rt)  dr,
$$
where $C(\delta, k)$ is a constant depending only on $\da$ and $k$. As a consequence, one has
\be\label{hardy2'}
\int^{\sqrt{A/B}}_{\sqrt{A/(4B)} } \lt(\sqrt{A/B}-r\rt)^{k-2} F^2 dr \le C \int^{\sqrt{A/B}}_{\sqrt{A/(4B)} }  \lt(\sqrt{A/B}-r\rt)^k \lt( F^2 + F_r^2 \rt)  dr,
\ee
where $C$ is a constant depending on $A$, $B$ and $k$.

\vskip 0.5cm
{\bf Notations:}

1) Throughout the rest of paper,   $C$  will denote a positive constant which  only depend on the parameters of the problem,
$\ga$ and $M$, but does not depend on the data. They are referred as universal and can change
from one inequality to another one. Also we use $C(\beta)$ to denote  a certain positive constant
depending on quantity $\beta$.

2) We will employ the notation $a\lesssim b$ to denote $a\le C b$, $a \sim b$ to denote $C^{-1}b\le a\le Cb$ and $a \gtrsim b$ to denote $a \ge  C^{-1} b$,
where  $C$ is the universal constant  as defined
above.

3) In the rest of the paper, we will use the notations
$$ \int=:\int_{\mathcal{I}}  \ , \ \   \|\cdot\|=:\|\cdot\|_{L^2(\mathcal{I})}  \ \    {\rm and} \ \   \|\cdot\|_{L^{\infty}}=:\|\cdot\|_{L^{\infty}(\mathcal{I})}.$$

4)
We set
$$\sigma(r):=\bar\rho_0^{\ga-1}(r)=A-Br^2, \ \ r \in \mathcal{I}.$$
Then $\mathcal{ E}_{j}$ and $\mathcal{ E}_{j,i}$ can be rewritten as
\bee\label{}\begin{split}
&\mathcal{ E}_{j}(t)  = (1+ t)^{2j} \int    \lt[r^4 \sigma^\aa  \lt(\pl_t^j \zeta\rt)^2 +r^2 \sigma^{\aa+1}  \lt|\pl_t^j \lt( \zeta, r\zeta_r \rt)  \rt|^2  + (1+ t)  r^4 \sigma^\aa    \lt(\pl_t^{j} \zeta_t\rt)^2 \rt](r,t)   dr,\\
&\mathcal{ E}_{j, i}(t) =  (1+  t)^{2j}  \int \lt[r^2 \sigma^{\aa+i-1 }  \lt(\pl_t^j \pl_r^i \zeta\rt)^2 + r^4 \sigma^{\aa+i+1}  \lt(\pl_t^j \pl_r^{i+1} \zeta\rt)^2 \rt](r,t) dr  .
\end{split}\eee

5) We set
$$\mathcal{I}_o:= \lt( 0, \ {\sqrt{A/(4B)} } \rt) \ \  {\rm and} \ \  \mathcal{I}_b:= \lt(   {\sqrt{A/(4B)} } ,  \ {\sqrt{A/B} } \rt).$$
Then
$$\mathcal{I}=\mathcal{I}_o \cup \mathcal{I}_b.$$
Moreover, it gives from the Hardy  inequality \ef{hardy2'} that for $k>1$,
\be\label{hardy2}
\int_{\mathcal{I}_b } \sigma^{k-2}(r) F^2 dr \le C(A,B,k) \int_{\mathcal{I}_b}  \sigma^k(r)  \lt( F^2 + F_r^2 \rt)  dr,
\ee
provided that the right-hand side of \ef{hardy2} is finite.

\subsection{Elliptic estimates}
In this subsection, we prove the following elliptic estimates.
\begin{prop}\label{prop1} Suppose that \ef{apb} holds for suitably small positive number $\ea_0\in(0,1)$. Then it holds that for $t\in[0,T]$,
\be\label{prop1est}\begin{split}
   \mathcal{E}_{j,i}(t)\lesssim  \sum_{\iota=0}^{i+j} \mathcal{E}_\iota(t) \ \ {\rm when} \ \  j\ge 0, \ \ i\ge 1, \ \  i+j\le l.
\end{split}\ee
\end{prop}
The proof of this proposition consists of Lemma \ref{lem321} and Lemma \ref{lem322}.

\subsubsection{Lower-order elliptic estimates}
Dividing equation \ef{pertb} by $\bar\rho_0$, one has
\bee\label{}\begin{split}
   &r \zeta_{tt}   +
   r    \zeta_t
   +   \sigma\lt(\tilde\eta_r+ \zeta\rt)^2\left[      \lt(\tilde\eta_r+\zeta\rt)^{-2\ga}\lt(\tilde\eta_r+\zeta +r\zeta_r\rt)^{-\ga}   \right]_r \\
   & + \f{\ga}{\ga-1} \sigma_r \left[     \lt(\tilde\eta_r+\zeta\rt)^{2-2\ga}\lt(\tilde\eta_r+\zeta +r\zeta_r\rt)^{-\ga}  - \tilde\eta_r^{2-3\ga}    \right]
   =0 .
\end{split}\eee
Note that
\bee\label{}\begin{split}
  & \lt(\tilde\eta_r+ \zeta\rt)^2\left[      \lt(\tilde\eta_r+\zeta\rt)^{-2\ga}\lt(\tilde\eta_r+\zeta +r\zeta_r\rt)^{-\ga}   \right]_r \\
 &  =-2\ga \lt(\tilde\eta_r+\zeta\rt)^{1-2\ga}\lt(\tilde\eta_r+\zeta +r\zeta_r\rt)^{-\ga} \zeta_r
   -\ga \lt(\tilde\eta_r+\zeta\rt)^{2-2\ga}\lt(\tilde\eta_r+\zeta +r\zeta_r\rt)^{-\ga-1} \lt(2\zeta_r+r\zeta_{rr}\rt) \\
 &  =- \ga \tilde\eta_r^{1-3\ga}\lt(4\zeta_r+r\zeta_{rr}\rt) +\mathfrak{J}_1,
\end{split}\eee
\bee\label{}\begin{split}
  \lt(\tilde\eta_r+\zeta\rt)^{2-2\ga}\lt(\tilde\eta_r+\zeta +r\zeta_r\rt)^{-\ga}  - \tilde\eta_r^{2-3\ga}
= -\ga \tilde\eta_r^{1-3\ga}  \lt(r\zeta_r\rt) +(2-3\ga)\tilde\eta_r^{1-3\ga} \zeta + \mathfrak{J}_2 ,
\end{split}\eee
where
\be\label{oj1j2}\begin{split}
 \mathfrak{J}_1 :
   =&-2\ga \lt[\lt(\tilde\eta_r+\zeta\rt)^{1-2\ga}\lt(\tilde\eta_r+\zeta +r\zeta_r\rt)^{-\ga} -\tilde\eta_r^{1-3\ga} \rt]  \zeta_r
   \\
   &-\ga \lt[\lt(\tilde\eta_r+\zeta\rt)^{2-2\ga}\lt(\tilde\eta_r+\zeta +r\zeta_r\rt)^{-\ga-1} -\tilde\eta_r^{1-3\ga} \rt]\lt(2\zeta_r+r\zeta_{rr}\rt)
 ,\\
 \mathfrak{J}_2:=& \lt(\tilde\eta_r+\zeta\rt)^{2-2\ga}\lt(\tilde\eta_r+\zeta +r\zeta_r\rt)^{-\ga}  - \tilde\eta_r^{2-3\ga} +\ga \tilde\eta_r^{1-3\ga}  \lt(r\zeta_r\rt) -(2-3\ga)\tilde\eta_r^{1-3\ga} \zeta  .
\end{split}\ee
Then,
\be\label{o3-7}\begin{split}
&   \ga \tilde\eta_r^{1-3\ga} \lt[r\sigma  \zeta_{rr} +  4\sigma \zeta_r   +\f{\ga}{\ga-1} r\sigma_r  \zeta_r  \rt] \\
=&r \zeta_{tt}   +
   r    \zeta_t +\f{\ga (2-3\ga)}{\ga-1} \sigma_r\tilde\eta_r^{1-3\ga} \zeta  + \sigma \mathfrak{J}_1
     + \f{\ga}{\ga-1} \sigma_r \mathfrak{J}_2
 .
\end{split}\ee
\begin{lem}\label{lem321} Assume that \ef{apb} holds for suitably small positive number $\ea_0\in(0,1)$. Then,
$$\mathcal{E}_{0,1}(t) \lesssim  {\mathcal{\mathcal{E}}}_0(t) +  \mathcal{E}_1(t), \ \  0\le t\le T.$$
\end{lem}
{\bf Proof}. Multiply equation \ef{o3-7} by $\tilde\eta_r^{3\ga-1}r  \sigma^{\aa/2}$ and square the spatial $L^2$-norm of the product to obtain
\be\label{o3-8}\begin{split}
&   \lt\| r^2 \sigma^{1+\f{\aa}{2}}  \zeta_{rr} +  4 r \sigma^{1+\f{\aa}{2}}   \zeta_r   +\lt(1+\aa\rt)  r^2 \sigma^{ \f{\aa}{2}}\sigma_r \zeta_r  \rt\|^2 \\
\lesssim & (1+t)^2 \lt( \lt\| r^2 \sigma^{\f{\aa}{2}} \zeta_{tt} \rt\|^2  +
  \lt\| r^2 \sigma^{\f{\aa}{2}} \zeta_{t } \rt\|^2 +\lt\| r\sigma^{1+\f{\aa}{2}}  \mathfrak{J}_1 \rt\|^2
     +\lt\|r \sigma^{ \f{\aa}{2}} \sigma_r \mathfrak{J}_2\rt\|^2   \rt)+ \lt\| r \sigma^{\f{\aa}{2}} \zeta  \rt\|^2 \\
     \lesssim  &   \mathcal{E}_1 + (1+t)^2 \lt(  \lt\| r\sigma^{1+\f{\aa}{2}}  \mathfrak{J}_1 \rt\|^2
     +\lt\|r \sigma^{ \f{\aa}{2}} \sigma_r \mathfrak{J}_2\rt\|^2   \rt)+ \lt\| r \sigma^{\f{\aa}{2}} \zeta  \rt\|^2
\end{split}\ee
where we have used \ef{decay} and the definition of $\mathcal{E}_1$. It follows from the Taylor expansion,  \ef{basic} and \ef{apb} that
\bee\label{}\begin{split}
& |\mathfrak{J}_1| \lesssim (1+t)^{-\frac{3\ga}{3\ga-1}}\lt(\lt|r\zeta_r\rt| +\lt|\zeta\rt|\rt)\lt(\lt|r\zeta_{rr}\rt| + \lt|\zeta_r\rt|
 \rt) \lesssim (1+t)^{-\frac{3\ga}{3\ga-1}}\ea_0 \lt(\lt|r\zeta_{rr}\rt| + \lt|\zeta_r\rt|
 \rt),  \\
& \lt|\mathfrak{J}_2\rt|\lesssim (1+t)^{-\frac{3\ga}{3\ga-1}}\lt(\lt|r\zeta_r\rt|^2 +\lt|\zeta\rt|^2\rt) \lesssim (1+t)^{-\frac{3\ga}{3\ga-1}}\ea_0  \lt(\lt|r\zeta_r\rt| +\lt|\zeta\rt|\rt).
\end{split}\eee
Thus,
\be\label{o3-9}\begin{split}
& (1+t)^2 \lt(  \lt\| r\sigma^{1+\f{\aa}{2}}  \mathfrak{J}_1 \rt\|^2
     +\lt\|r \sigma^{ \f{\aa}{2}} \sigma_r \mathfrak{J}_2\rt\|^2   \rt) \\
     \lesssim  & \ea_0^2 \lt(\lt\| r^2 \sigma^{1+\f{\aa}{2}}  \zeta_{rr} \rt\|^2
     +\lt\| r  \sigma^{1+\f{\aa}{2}}  \zeta_{r } \rt\|^2+\lt\|r^2 \sigma^{ \f{\aa}{2}} \sigma_r \zeta_r\rt\|^2  +\lt\|r  \sigma^{ \f{\aa}{2}} \sigma_r \zeta \rt\|^2  \rt).
\end{split}\ee
Note that
\be\label{o3-10}\begin{split}
 \lt\| r \sigma^{\f{\aa}{2}} \zeta  \rt\|^2 =   \int_{\mathcal{I}_o} r^2 \sigma^{\aa} \zeta^2 dr
 +\int_{\mathcal{I}_b} r^2 \sigma^{\aa} \zeta^2 dr
 \lesssim  \int_{\mathcal{I}_o} r^2 \sigma^{1+ \aa} \zeta^2 dr
 +\int_{\mathcal{I}_b} r^4 \sigma^{\aa} \zeta^2 dr \lesssim {\mathcal{\mathcal{E}}}_0.
\end{split}\ee
Then, it yields from \ef{o3-8}, \ef{o3-9} and \ef{o3-10} that
\be\label{o3-11}\begin{split}
&   \lt\| r^2 \sigma^{1+\f{\aa}{2}}  \zeta_{rr} +  4 r \sigma^{1+\f{\aa}{2}}   \zeta_r   +\lt(1+\aa\rt)  r^2 \sigma^{ \f{\aa}{2}}\sigma_r \zeta_r  \rt\|^2 \\
     \lesssim  &  {\mathcal{\mathcal{E}}}_0 +  \mathcal{E}_1  +\ea_0^2 \lt(\lt\| r^2 \sigma^{1+\f{\aa}{2}}  \zeta_{rr} \rt\|^2
     +\lt\| r  \sigma^{1+\f{\aa}{2}}  \zeta_{r } \rt\|^2+\lt\|r^2 \sigma^{ \f{\aa}{2}} \sigma_r \zeta_r\rt\|^2   \rt).
\end{split}\ee

In what follows, we analyze the left-hand side of \ef{o3-11}, which  can be expanded as
\be\label{o3-12}\begin{split}
 &  \lt\| r^2 \sigma^{1+\f{\aa}{2}}  \zeta_{rr} +  4 r \sigma^{1+\f{\aa}{2}}   \zeta_r   +\lt(1+\aa\rt)  r^2 \sigma^{ \f{\aa}{2}}\sigma_r \zeta_r  \rt\|^2 \\
 =& \lt\| r^2 \sigma^{1+\f{\aa}{2}}  \zeta_{rr} \rt\|^2
     +16\lt\| r  \sigma^{1+\f{\aa}{2}}  \zeta_{r } \rt\|^2+(1+\aa)^2\lt\|r^2 \sigma^{ \f{\aa}{2}} \sigma_r \zeta_r\rt\|^2  \\
    & + \int \lt[4r^3\sigma^{2+\aa} +(1+\aa) r^4 \sigma^{1+\aa} \sigma_r \rt] \lt(\zeta_r^2\rt)_r dr
     +8(1+\aa)\int r^3\sigma^{1+\aa}\sigma_r \zeta_r^2 dr .
\end{split}\ee
With the help of the integration by parts and the fact $\sigma_r=-2Br$, one has
\bee\label{}\begin{split}
&\int \lt[4r^3\sigma^{2+\aa} +(1+\aa) r^4 \sigma^{1+\aa} \sigma_r \rt] \lt(\zeta_r^2\rt)_r dr \\
= & - \int \lt[4r^3\sigma^{2+\aa} +(1+\aa) r^4 \sigma^{1+\aa} \sigma_r \rt]_r  \zeta_r^2  dr\\
 \ge & -12 \int r^2 \sigma^{2+\aa} \zeta_r^2 dr-(1+\aa)^2 \int r^4 \sigma^\aa \sigma_r^2 \zeta_r^2dr
-C \int r^4 \sigma^{1+\aa} \zeta_r^2 dr.
\end{split}\eee
Substitute this into \ef{o3-12} and use  $\sigma_r=-2Br$  to give
\bee\label{}\begin{split}
 &  \lt\| r^2 \sigma^{1+\f{\aa}{2}}  \zeta_{rr} +  4 r \sigma^{1+\f{\aa}{2}}   \zeta_r   +\lt(1+\aa\rt)  r^2 \sigma^{ \f{\aa}{2}}\sigma_r \zeta_r  \rt\|^2 \\
 \ge & \lt\| r^2 \sigma^{1+\f{\aa}{2}}  \zeta_{rr} \rt\|^2
     +4\lt\| r  \sigma^{1+\f{\aa}{2}}  \zeta_{r } \rt\|^2 -C \int r^4 \sigma^{1+\aa} \zeta_r^2 dr.
\end{split}\eee
In view of \ef{o3-11}, we then see that
\be\label{o3-13}\begin{split}
&   \lt\| r^2 \sigma^{1+\f{\aa}{2}}  \zeta_{rr} \rt\|^2
     +4\lt\| r  \sigma^{1+\f{\aa}{2}}  \zeta_{r } \rt\|^2  \\
      \lesssim  &  {\mathcal{\mathcal{E}}}_0 +  \mathcal{E}_1  +\ea_0^2 \lt(\lt\| r^2 \sigma^{1+\f{\aa}{2}}  \zeta_{rr} \rt\|^2
     +\lt\| r  \sigma^{1+\f{\aa}{2}}  \zeta_{r } \rt\|^2+\lt\|r^2 \sigma^{ \f{\aa}{2}} \sigma_r \zeta_r\rt\|^2   \rt).
\end{split}\ee
On the other hand, it follows from \ef{o3-11} and \ef{o3-13} that
\bee\label{}\begin{split}
& \lt\|\lt(1+\aa\rt)  r^2 \sigma^{ \f{\aa}{2}}\sigma_r \zeta_r\rt\|^2 \\
 \le & 2 \lt\| r^2 \sigma^{1+\f{\aa}{2}}  \zeta_{rr} +  4 r \sigma^{1+\f{\aa}{2}}   \zeta_r   +\lt(1+\aa\rt)  r^2 \sigma^{ \f{\aa}{2}}\sigma_r \zeta_r  \rt\|^2
+2 \lt\| r^2 \sigma^{1+\f{\aa}{2}}  \zeta_{rr} +  4 r \sigma^{1+\f{\aa}{2}}   \zeta_r    \rt\|^2  \\
     \lesssim  &  {\mathcal{\mathcal{E}}}_0 +  \mathcal{E}_1  +\ea_0^2 \lt(\lt\| r^2 \sigma^{1+\f{\aa}{2}}  \zeta_{rr} \rt\|^2
     +\lt\| r  \sigma^{1+\f{\aa}{2}}  \zeta_{r } \rt\|^2+\lt\|r^2 \sigma^{ \f{\aa}{2}} \sigma_r \zeta_r\rt\|^2   \rt).
\end{split}\eee
This, together with \ef{o3-13}, gives
\be\label{o3-14}\begin{split}
&\lt\| r^2 \sigma^{1+\f{\aa}{2}}  \zeta_{rr} \rt\|^2
     +\lt\| r  \sigma^{1+\f{\aa}{2}}  \zeta_{r } \rt\|^2+\lt\|r^2 \sigma^{ \f{\aa}{2}} \sigma_r \zeta_r\rt\|^2\\
          \lesssim  &  {\mathcal{\mathcal{E}}}_0 +  \mathcal{E}_1  +\ea_0^2 \lt(\lt\| r^2 \sigma^{1+\f{\aa}{2}}  \zeta_{rr} \rt\|^2
     +\lt\| r  \sigma^{1+\f{\aa}{2}}  \zeta_{r } \rt\|^2+\lt\|r^2 \sigma^{ \f{\aa}{2}} \sigma_r \zeta_r\rt\|^2   \rt),
\end{split}\ee
which implies, with the aid of the smallness of $\ea_0$, that
\bee\label{}\begin{split}
 \lt\| r^2 \sigma^{1+\f{\aa}{2}}  \zeta_{rr} \rt\|^2
     +\lt\| r  \sigma^{1+\f{\aa}{2}}  \zeta_{r } \rt\|^2+\lt\|r^2 \sigma^{ \f{\aa}{2}} \sigma_r \zeta_r\rt\|^2
          \lesssim     {\mathcal{\mathcal{E}}}_0 +  \mathcal{E}_1.
\end{split}\eee
In view of  $\sigma_r=-2B r$, we then see that
\bee\label{o3-15}\begin{split}
 \lt\| r^2 \sigma^{1+\f{\aa}{2}}  \zeta_{rr} \rt\|^2
     +\lt\| r  \sigma^{1+\f{\aa}{2}}  \zeta_{r } \rt\|^2+\lt\|r^3 \sigma^{ \f{\aa}{2}}  \zeta_r\rt\|^2
          \lesssim     {\mathcal{\mathcal{E}}}_0 +  \mathcal{E}_1 ,
\end{split}\eee
which implies
\be\label{o3-16}\begin{split}
 \lt\| r  \sigma^{ \f{\aa}{2}}  \zeta_{r } \rt\|^2
  = & \int_{\mathcal{I}_o} r^2 \sigma^{\aa} \zeta_r ^2 dr
 +\int_{\mathcal{I}_b} r^2 \sigma^{\aa} \zeta_r^2 dr \\
 \lesssim & \int_{\mathcal{I}_o} r^2 \sigma^{2+ \aa} \zeta_r^2 dr
 +\int_{\mathcal{I}_b} r^6 \sigma^{\aa} \zeta_r^2 dr \lesssim {\mathcal{\mathcal{E}}}_0 +  \mathcal{E}_1.
\end{split}\ee
This finishes the proof of Lemma \ref{lem321}. $\Box$

\subsubsection{Higher-order elliptic estimates}
For $i\ge 1$ and $j\ge 0$, it yields from $\pl_t^j\pl_r^{i-1} \ef{o3-7}$ and $\sigma_r=-2Br$ that
\be\label{o3-22}\begin{split}
&   \ga \tilde\eta_r^{1-3\ga} \lt[r\sigma \pl_t^j \pl_r^{i+1} \zeta  +  (i+3) \sigma \pl_t^j \pl_r^i \zeta    +\lt(\aa+i\rt) r\sigma_r \pl_t^j \pl_r^i \zeta   \rt] \\
=&r \pl_t^{j+2}\pl_r^{i-1}\zeta   +r \pl_t^{j+1}\pl_r^{i-1}\zeta
      + \mathfrak{P}_1 + \mathfrak{P}_2
 .
\end{split}\ee
where
\be\label{3-23}\begin{split}
\mathfrak{P}_1:= &-\ga   \sum_{\iota=1}^j \lt[\pl_t^\iota \lt( \tilde\eta_r^{1-3\ga} \rt)\rt] \pl_t^{j-\iota} \lt[r\sigma \pl_r^{i+1} \zeta  +  (i+3) \sigma \pl_r^i \zeta    +\lt(\aa+i\rt) r\sigma_r  \pl_r^i \zeta   \rt] \\
&-\ga\pl_t^j\lt\{\tilde\eta_r^{1-3\ga} \lt[\sum_{\iota=2}^{i-1} C_{i-1}^\iota \lt[\pl_r^{\iota}(r \sigma )\rt] \pl_r^{i+1-\iota} \zeta  + 4\sum_{\iota=1}^{i-1} C_{i-1}^\iota \lt(\pl_r^{\iota}\sigma \rt) \pl_r^{i -\iota} \zeta\rt.\rt. \\
& \lt.\lt.  + (\aa+1)\sum_{\iota=1}^{i-1} C_{i-1}^\iota \lt[\pl_r^{\iota}\lt(r\sigma_r\rt) \rt] \pl_r^{i -\iota} \zeta \rt]  \rt\} + (i-1) \pl_r^{i-2} \lt(\pl_t^{j+2}\zeta   +  \pl_t^{j+1} \zeta\rt) \\
&-\f{2\ga (2-3\ga)B}{\ga-1} \pl_t^j \lt[\tilde\eta_r^{1-3\ga} \lt(r \pl_r^{i-1} \zeta +(i-1) \pl_r^{i-2} \zeta \rt)\rt],
\end{split}\ee
\be\label{3-24}\begin{split}
\mathfrak{P}_2:= &  \pl_r^{i-1} \lt(\sigma \pl_t^j \mathfrak{J}_1\rt)
     + (1+\aa) \pl_r^{i-1} \lt(  \sigma_r\pl_t^j  \mathfrak{J}_2 \rt).
\end{split}\ee
(Recall that $\mathfrak{J}_1$ and $\mathfrak{J}_2$ are defined in \ef{oj1j2}.)
Here and thereafter $C_m^j$ is used to denote the binomial coefficients for $0\le j\le m$,
$$C_m^j = \frac{m!}{j!(m-j)!}$$
and summations $\sum_{\iota=1}^{i-1}$ and $\sum_{\iota=2}^{i-1}$ should be understood as zero when $i=1$ and $i=1,2$, respectively.  Multiply equation \ef{o3-22} by $\tilde\eta_r^{3\ga-1}r  \sigma^{(\aa+i-1)/2}$, square the spatial $L^2$-norm of the product and use \ef{decay} to give
\bee\label{}\begin{split}
&    \lt\|r^2\sigma^{\frac{\aa+i+1}{2}} \pl_t^j \pl_r^{i+1} \zeta  +  (i+3) r\sigma^{\frac{\aa+i+1}{2}} \pl_t^j \pl_r^i \zeta    +\lt(\aa+i\rt) r^2  \sigma^{\frac{\aa+i-1}{2}} \sigma_r \pl_t^j \pl_r^i \zeta   \rt\|^2  \\
\lesssim & (1+t)^2 \lt(\lt\| r^2 \sigma ^{\frac{\aa+i-1}{2}} \pl_t^{j+2}\pl_r^{i-1}\zeta  \rt\|^2   +\lt\| r^2 \sigma ^{\frac{\aa+i-1}{2}} \pl_t^{j+1}\pl_r^{i-1}\zeta  \rt\|^2 \rt)
   \\&   +(1+t)^2 \lt( \lt\|r  \sigma ^{\frac{\aa+i-1}{2}}  \mathfrak{P}_1 \rt\|^2 + \lt\|r  \sigma ^{\frac{\aa+i-1}{2}}  \mathfrak{P}_2 \rt\|^2 \rt)
 .
\end{split}\eee
Similar to the derivation of \ef{o3-14} and \ef{o3-16}, we can then obtain
\be\label{3-25}\begin{split}
& (1+t)^{-2j} \mathcal{ E}_{j, i}(t) =  \lt\|r^2\sigma^{\frac{\aa+i+1}{2}} \pl_t^j \pl_r^{i+1} \zeta    \rt\|^2 + \lt\|  r\sigma^{\frac{\aa+i-1}{2}} \pl_t^j \pl_r^i \zeta       \rt\|^2   \\
\lesssim & \lt\|r^2\sigma^{\frac{\aa+i }{2}} \pl_t^j \pl_r^{i } \zeta    \rt\|^2 +  (1+t)^2 \lt(\lt\| r^2 \sigma ^{\frac{\aa+i-1}{2}} \pl_t^{j+2}\pl_r^{i-1}\zeta  \rt\|^2   +\lt\| r^2 \sigma ^{\frac{\aa+i-1}{2}} \pl_t^{j+1}\pl_r^{i-1}\zeta  \rt\|^2 \rt)
   \\&   +(1+t)^2 \lt( \lt\|r  \sigma ^{\frac{\aa+i-1}{2}}  \mathfrak{P}_1 \rt\|^2 + \lt\|r  \sigma ^{\frac{\aa+i-1}{2}}  \mathfrak{P}_2 \rt\|^2 \rt)
 .
\end{split}\ee
We will use this estimate to prove the following lemma by the mathematical induction.

\begin{lem}\label{lem322} Assume that \ef{apb} holds for suitably small positive number $\ea_0\in(0,1)$. Then for $j\ge 0$, $i\ge 1$ and $2\le i+j\le l$,
\be\label{lem322est}\begin{split}
\mathcal{E}_{j,i}(t)\lesssim   \sum_{\iota=0}^{i+j} \mathcal{E}_\iota(t), \ \  t\in [0,T].
\end{split}\ee
\end{lem}
{\bf Proof}. We use the induction for $i+j$ to prove this lemma. As shown in Lemma \ref{lem321} we know that  \ef{lem322est} holds for $i+j=1$. For $1\le k\le l- 1$, we make the induction hypothesis that  \ef{lem322est} holds for all
$ j\ge 0$, $i\ge 1$ and  $i+j\le k$, that is,
\be\label{3-26}
\mathcal{ E}_{j, i}(t) \lesssim \sum_{\iota=0}^{i+j}\mathcal{ E}_{\iota}(t),  \ \  j\ge 0, \ \ i\ge 1, \ \  i+j \le k.
\ee
It then suffices to prove  \ef{lem322est}  for $ j\ge 0$, $i\ge 1$ and  $i+j=k+1$.
(Indeed, there exists an order of $(i,j)$ for the proof. For example,  when $i+j=k+1$ we will bound $\mathcal{ E}_{k+1-\iota, \iota}$ from $\iota=1$ to $k+1$ step by step.)

Before going to the estimate, we  notice a fact that $\mathcal{E}_{j,0}\lesssim \mathcal{E}_j$ for $j=0,\cdots,l$. Indeed, it follows from  \ef{hardy2} that
\bee\label{}\begin{split}
 \int_{\mathcal{I}_b}   \sigma^{\aa -1 }  \lt(\pl_t^j  \zeta\rt)^2 dr
 \lesssim & \int_{\mathcal{I}_b}   \sigma^{\aa +1 }  \lt[\lt(\pl_t^j  \zeta\rt)^2 + \lt(\pl_t^j  \zeta_r\rt)^2  \rt]dr \\
 \lesssim & \int_{\mathcal{I}_b}   \sigma^{\aa +1 }  \lt[r^2 \lt(\pl_t^j  \zeta\rt)^2 + r^4 \lt(\pl_t^j  \zeta_r\rt)^2  \rt]dr \le  (1+t)^{-2j} \mathcal{E}_j(t),
  \end{split}\eee
which implies
\be\label{3-27}\begin{split}
 \mathcal{E}_{j,0}(t) = & (1+t)^{2j} \int \lt[r^2 \sigma^{\aa -1 }  \lt(\pl_t^j  \zeta\rt)^2 + r^4 \sigma^{\aa +1}  \lt(\pl_t^j    \zeta_r\rt)^2 \rt](r,t)  dr \\
\le&  (1+t)^{2j}  \lt[ \int_{\mathcal{I}_o}  r^2 \sigma^{\aa -1 }  \lt(\pl_t^j  \zeta\rt)^2 (r,t) dr + \int_{\mathcal{I}_b}  r^2 \sigma^{\aa -1 }  \lt(\pl_t^j  \zeta\rt)^2(r,t) dr \rt] + \mathcal{E}_j(t)
\\
\lesssim  & (1+t)^{2j}  \lt[ \int_{\mathcal{I}_o}  r^2 \sigma^{\aa + 1 }  \lt(\pl_t^j  \zeta\rt)^2 (r,t) dr + \int_{\mathcal{I}_b}  \sigma^{\aa -1 }  \lt(\pl_t^j  \zeta\rt)^2 (r,t) dr \rt] + \mathcal{E}_j(t) \\
 \lesssim & \mathcal{E}_j(t), \ \  j=0,1,\cdots,l.
  \end{split}\ee
This, together with the induction hypothesis  \ef{3-26}, gives
\be\label{3-28}
\mathcal{ E}_{j, i}(t) \lesssim \sum_{\iota=0}^{i+j}\mathcal{ E}_{\iota}(t),  \ \  j\ge 0, \ \ i\ge 0, \ \  i+j \le k.
\ee

In what follows, we assume $ j\ge 0$, $i\ge 1$ and  $i+j=k+1\le l$. First, We estimate $\mathfrak{P}_1$ and $\mathfrak{P}_2$ given by \ef{3-23} and \ef{3-24}, respectively.
For $\mathfrak{P}_1$, it follows from \ef{decay} and $\sigma_r=-2B r$ that
\bee\label{}\begin{split}
|\mathfrak{P}_1| \lesssim  & \sum_{\iota=1}^j (1+t)^{-1-\iota}\lt(\lt| rw \pl_t^{j-\iota} \pl_r^{i+1} \zeta\rt| + \lt|\pl_t^{j-\iota} \pl_r^{i } \zeta\rt|\rt)
+\sum_{\iota=0}^j \sum_{m=1}^{i-1} (1+t)^{-1-\iota} \lt|\pl_t^{j-\iota}\pl_r^m \zeta \rt|\\
& + (i-1) \lt(\lt| \pl_t^{j+2} \pl_r^{i-2} \zeta \rt| +\lt| \pl_t^{j+1} \pl_r^{i-2} \zeta \rt| + \sum_{\iota=0}^j (1+t)^{-1-\iota} \lt| \pl_t^{j-\iota} \pl_r^{i-2} \zeta \rt|\rt)\\
&+ \sum_{\iota=0}^j (1+t)^{-1-\iota}  \lt| r  \pl_t^{j-\iota} \pl_r^{i-1} \zeta \rt|,
 \end{split}\eee
which implies
\bee\label{}\begin{split}
 \lt\|r  \sigma ^{\frac{\aa+i-1}{2}}  \mathfrak{P}_1 \rt\|^2  \lesssim   \sum_{\iota=1}^j (1+t)^{-2-2\iota}\lt(\lt\| r^2 \sigma ^{\frac{\aa+i+ 1}{2}}  \pl_t^{j-\iota} \pl_r^{i+1} \zeta\rt\|^2 + \lt\|r \sigma ^{\frac{\aa+i-1}{2}} \pl_t^{j-\iota} \pl_r^{i } \zeta\rt\|^2 \rt)
\\
 +\sum_{\iota=0}^j  (1+t)^{-2-2\iota} \lt( \sum_{m=1}^{i-1} \lt\|r \sigma ^{\frac{\aa+i-1}{2}} \pl_t^{j-\iota}\pl_r^m \zeta \rt\|^2  +   \lt\| r^2   \sigma ^{\frac{\aa+i-1}{2}}\pl_t^{j-\iota} \pl_r^{i-1} \zeta \rt\|^2 \rt) \\
  + (i-1)^2 \lt(\sum_{\iota=j+1}^{j+2}\lt\|r \sigma ^{\frac{\aa+i-1}{2}}  \pl_t^{\iota} \pl_r^{i-2} \zeta \rt\|^2    + \sum_{\iota=0}^j (1+t)^{-2-2\iota} \lt\| r \sigma ^{\frac{\aa+i-1}{2}}  \pl_t^{j-\iota} \pl_r^{i-2} \zeta \rt\|^2 \rt).
 \end{split}\eee
So,
\be\label{3-29}
 \lt\|r  \sigma ^{\frac{\aa+i-1}{2}}  \mathfrak{P}_1 \rt\|^2  \lesssim  \lt\{    \begin{split}  & (1+t)^{-2-2j} \lt(\sum_{\iota=0}^{j-1} \mathcal{E}_{ \iota,1} + \sum_{\iota=0}^j\mathcal{E}_{ \iota}\rt)(t),   &i=1, \\
&   (1+t)^{-2-2j} \lt(\sum_{\iota=0}^{j-1} \mathcal{E}_{ \iota,i} + \sum_{\iota=0}^j \sum_{m=1}^{i-1} \mathcal{E}_{ \iota,m } + \sum_{\iota=0}^{j+2} \mathcal{E}_{ \iota, i-2} \rt) (t), & i\ge 2.
 \end{split} \rt. \ee
For $\mathfrak{P}_2$, it follows from \ef{decay}, \ef{apb}, \ef{basic} and $\sigma_r=-2B r$ that
\bee\label{} \begin{split}
 |\mathfrak{P}_2| \lesssim & \sum_{n=0}^j \sum_{m=0}^{i-1} K_{nm} \lt(\lt|\pl_t^{j-n}\pl_r^{i-1-m}\lt( \sigma  r\zeta_{rr}  \rt)\rt|+
 \lt|\pl_t^{j-n}\pl_r^{i-1-m}\lt( \sigma  \zeta_{r }  \rt)\rt|  \rt.\\
& \qquad \qquad  \qquad   \lt.  +
 \lt|\pl_t^{j-n}\pl_r^{i-1-m}\lt(  \sigma_r  r \zeta_{r }  \rt)\rt|+
 \lt|\pl_t^{j-n}\pl_r^{i-1-m}\lt(  \sigma_r \zeta   \rt)\rt|\rt)\\
 \lesssim & \sum_{n=0}^j \sum_{m=0}^{i-1} K_{nm} \lt(\lt| \sigma  r \pl_t^{j-n}\pl_r^{i-m +1}  \zeta  \rt|+   \sum_{\iota=0}^{i-m }  \lt|\pl_t^{j-n}\pl_r^{\iota}  \zeta  \rt|\rt)
 =:  \sum_{n=0}^j \sum_{m=0}^{i-1} \mathfrak{P}_{2nm} ,
 \end{split}   \eee
where
\bee\label{}\begin{split}
&  K_{00} = \ea_0(1+t)^{-1-\frac{1}{3\ga-1}}; \\
&  K_{10}= \ea_0(1+t)^{-2-\frac{1}{3\ga-1}}, \ \
  K_{01} = (1+t)^{-1-\frac{1}{3\ga-1}}\lt(\ea_0 + |r\pl_r^2 \zeta |\rt);\\
&  K_{20}=\ea_0(1+t)^{-3-\frac{1}{3\ga-1}} + (1+t)^{-1-\frac{1}{3\ga-1}}\lt|r\pl_t^2\pl_r \zeta \rt|,\\
&K_{11}= (1+t)^{-2-\frac{1}{3\ga-1}} \lt( \ea_0 + \lt|r\pl_r^2 \zeta  \rt|\rt) + (1+t)^{-1-\frac{1}{3\ga -1}}\lt| r \pl_t\pl_r^2 \zeta  \rt| , \\
&K_{02}=(1+t)^{-1-\frac{1}{3\ga-1}}\lt(\lt|   \pl_r^2 \zeta  \rt| + \lt| r \pl_r^3 \zeta  \rt|\rt) + (1+t)^{-1-\frac{2}{3\ga-1}}\lt( \ea_0^2 + \lt|r\pl_r^2 \zeta  \rt|^2 \rt).
\end{split}\eee
We do not list here $K_{nm}$ for $n+m\ge 3$ since we can use the same method to estimate $\mathfrak{P}_{2nm}$ for $n+m\ge 3$ as that for $n+m\le 2$. Easily, $\mathfrak{P}_{200}$
and $\mathfrak{P}_{210}$ can be bounded by
\bee\label{}  \begin{split}
 \lt\|r  \sigma ^{\frac{\aa+i-1}{2}}  \mathfrak{P}_{200} \rt\|^2 \lesssim  & \ea_0^2 (1+t)^{-2}
\lt(\lt\|r^2 \sigma ^{\frac{\aa+i+1}{2}}  \pl_t^{j} \pl_r^{i+1} \zeta \rt\|^2 + \sum_{\iota=0}^i\lt\|r \sigma ^{\frac{\aa+i-1}{2}}  \pl_t^{j} \pl_r^{\iota} \zeta \rt\|^2\rt)\\
\lesssim & \ea_0^2(1+t)^{-2-2j} \lt(\mathcal{E}_{j,i} +\sum_{\iota=0}^{i-1} \mathcal{E}_{j,\iota} \rt)(t),
 \end{split}   \eee
\bee\label{}  \begin{split}
 \lt\|r  \sigma ^{\frac{\aa+i-1}{2}}  \mathfrak{P}_{210} \rt\|^2 \lesssim  & \ea_0^2 (1+t)^{-4}
\lt(\lt\|r^2 \sigma ^{\frac{\aa+i+1}{2}}  \pl_t^{j-1} \pl_r^{i+1} \zeta \rt\|^2 + \sum_{\iota=0}^i\lt\|r \sigma ^{\frac{\aa+i-1}{2}}  \pl_t^{j-1} \pl_r^{\iota} \zeta \rt\|^2\rt)\\
\lesssim & \ea_0^2(1+t)^{-2-2j}  \sum_{\iota=0}^{i} \mathcal{E}_{j-1,\iota}(t).
 \end{split}   \eee
For $\mathfrak{P}_{201}$, we use \ef{apb} to get $|\sigma^{{1}/{2}}\pl_r^2 \zeta |\lesssim \ea_0$
and then obtain
\bee\label{}  \begin{split}
 \lt\|r  \sigma ^{\frac{\aa+i-1}{2}}  \mathfrak{P}_{201} \rt\|^2 \lesssim  & \ea_0^2 (1+t)^{-2}
\lt(\lt\|r^2 \sigma ^{\frac{\aa+i }{2}}  \pl_t^{j } \pl_r^{i } \zeta \rt\|^2 + \sum_{\iota=0}^{i-1}\lt\|r \sigma ^{\frac{\aa+i-2}{2}}  \pl_t^{j } \pl_r^{\iota} \zeta \rt\|^2\rt)\\
\lesssim & \ea_0^2(1+t)^{-2-2j}  \sum_{\iota=0}^{i-1} \mathcal{E}_{j ,\iota}(t),
 \end{split}   \eee
For $\mathfrak{P}_{220}$, we   use \ef{apb} again to get
 $|r\sigma^{{1}/{2}}\pl_t^2 \pl_r \zeta |\lesssim \ea_0 (1+t)^{-2}$ and then achieve
 \bee\label{}  \begin{split}
 \lt\|r  \sigma ^{\frac{\aa+i-1}{2}}  \mathfrak{P}_{220} \rt\|^2 \lesssim  & \ea_0^2 (1+t)^{-6}
\lt(\lt\|r^2 \sigma ^{\frac{\aa+i }{2}}  \pl_t^{j-2 } \pl_r^{i+1 } \zeta \rt\|^2 + \sum_{\iota=0}^{i }\lt\|r \sigma ^{\frac{\aa+i-2}{2}}  \pl_t^{j-2 } \pl_r^{\iota} \zeta \rt\|^2\rt)\\
\lesssim & \ea_0^2(1+t)^{-2-2j}  \sum_{\iota=0}^{i+1} \mathcal{E}_{j-2 ,\iota}(t),
 \end{split}   \eee
because it can be derived from \ef{hardy2} that
 \bee\label{}  \begin{split}
& \lt\|r \sigma ^{\frac{\aa+i-2}{2}}  \pl_t^{j-2 } \pl_r^{i} \zeta \rt\|^2
 =\int_{\mathcal{I}_o}  r^2 \sigma ^{ {\aa+i-2} } \lt|  \pl_t^{j-2 } \pl_r^{i} \zeta \rt |^2 dr
 +\int_{\mathcal{I}_b}  r^2 \sigma ^{ {\aa+i-2} } \lt|  \pl_t^{j-2 } \pl_r^{i} \zeta \rt |^2 dr\\
 \lesssim  &\int_{\mathcal{I}_o}  r^2 \sigma ^{ {\aa+i-1} } \lt|  \pl_t^{j-2 } \pl_r^{i} \zeta \rt |^2 dr
 +\int_{\mathcal{I}_b}   \sigma ^{ {\aa+i-2} } \lt|  \pl_t^{j-2 } \pl_r^{i} \zeta \rt |^2dr \\
 \lesssim & \int_{\mathcal{I}_o}  r^2 \sigma ^{ {\aa+i-1} } \lt|  \pl_t^{j-2 } \pl_r^{i} \zeta \rt |^2 dr
 +\int_{\mathcal{I}_b}   \sigma ^{ {\aa+i } }\lt( \lt|  \pl_t^{j-2 } \pl_r^{i} \zeta \rt |^2 + \lt|  \pl_t^{j-2 } \pl_r^{i+1} \zeta \rt |^2\rt)dr\\
  \lesssim & \int_{\mathcal{I}_o}  r^2 \sigma ^{ {\aa+i-1} } \lt|  \pl_t^{j-2 } \pl_r^{i} \zeta \rt |^2 dr
 +\int_{\mathcal{I}_b}   \sigma ^{ {\aa+i } }\lt(r^2  \lt|  \pl_t^{j-2 } \pl_r^{i} \zeta \rt |^2 + r^2 \lt|  \pl_t^{j-2 } \pl_r^{i+1} \zeta \rt |^2\rt)dr\\
 \lesssim & \int   r^2 \sigma ^{ {\aa+i-1} } \lt|  \pl_t^{j-2 } \pl_r^{i} \zeta \rt |^2 dr
 +\int    r^2  \sigma ^{ {\aa+i } } \lt|  \pl_t^{j-2 } \pl_r^{i+1} \zeta \rt |^2dr.
 \end{split}   \eee
Similar to the estimate for $\mathfrak{P}_{220}$, we can obtain
 \bee\label{}  \begin{split}
 \lt\|r  \sigma ^{\frac{\aa+i-1}{2}}  \mathfrak{P}_{211} \rt\|^2 \lesssim  & \ea_0^2 (1+t)^{-4}
\lt(\lt\|r^2 \sigma ^{\frac{\aa+i-1 }{2}}  \pl_t^{j-1 } \pl_r^{i  } \zeta \rt\|^2 + \sum_{\iota=0}^{i-1 }\lt\|r \sigma ^{\frac{\aa+i-3}{2}}  \pl_t^{j-1 } \pl_r^{\iota} \zeta \rt\|^2\rt)\\
\lesssim & \ea_0^2(1+t)^{-2-2j}  \sum_{\iota=0}^{i } \mathcal{E}_{j-1 ,\iota}(t),
 \end{split}   \eee
 \bee\label{}  \begin{split}
 \lt\|r  \sigma ^{\frac{\aa+i-1}{2}}  \mathfrak{P}_{202} \rt\|^2 \lesssim  & \ea_0^2 (1+t)^{-2}
\lt(\lt\|r^2 \sigma ^{\frac{\aa+i-2 }{2}}  \pl_t^{j } \pl_r^{i-1 } \zeta \rt\|^2 + \sum_{\iota=0}^{i-2 }\lt\|r \sigma ^{\frac{\aa+i-4}{2}}  \pl_t^{j   } \pl_r^{\iota} \zeta \rt\|^2\rt)\\
\lesssim & \ea_0^2(1+t)^{-2-2j}  \sum_{\iota=0}^{i-1} \mathcal{E}_{j  ,\iota}(t) .
 \end{split}   \eee
It should be noted that $\mathfrak{P}_{211}$ and $\mathfrak{P}_{202}$ appear when  $i\ge 2$ and $i\ge 3$, respectively. This  ensures the application of the Hardy inequality \ef{hardy2}. Other cases can be done similarly, since the leading term of $K_{nm}$ is
$$\sum_{q=0}^n (1+t)^{-1-\frac{1}{3\ga-1}-q}\lt(\lt|r \pl_t^{n-q} \pl_r^{m+1}\zeta\rt| + \lt|  \pl_t^{n-q} \pl_r^{m }\zeta\rt| \rt)$$
and
\be\label{zuihou}\begin{split}
\sum_{n=0}^j \sum_{m=0}^{i-1}  \sum_{q=0}^n (1+t)^{-2 -2q}  \lt\|  r  \sigma ^{\frac{\aa+i-1}{2}} \lt(\lt|r \pl_t^{n-q} \pl_r^{m+1}\zeta\rt| + \lt|  \pl_t^{n-q} \pl_r^{m }\zeta\rt| \rt) \lt(\lt| \sigma  r \pl_t^{j-n}\pl_r^{i-m +1}  \zeta  \rt|  \rt.\rt. \\
 + \lt.\lt.   \sum_{\iota=0}^{i-m }  \lt|\pl_t^{j-n}\pl_r^{\iota}  \zeta  \rt|\rt)\rt\|^2
\lesssim  \ea_0^2(1+t)^{-2-2j}   \lt(\mathcal{E}_{j,i} +\sum_{0\le \iota \le j, \ p \ge 0, \ \iota+p \le i+j-1 } \mathcal{E}_{ \iota,p } \rt)(t).
\end{split}\ee
(Estimate \ef{zuihou} will be verified in the Appendix.) Now, we may conclude that
\be\label{3-30}\begin{split}
 \lt\|r  \sigma ^{\frac{\aa+i-1}{2}}  \mathfrak{P}_{2 } \rt\|^2
\lesssim & \ea_0^2(1+t)^{-2-2j}   \lt(\mathcal{E}_{j,i} +\sum_{0\le \iota \le j, \ p \ge 0, \ \iota+p \le i+j-1 } \mathcal{E}_{ \iota,p } \rt)(t) .
 \end{split}   \ee
Substitute \ef{3-29} and \ef{3-30} into \ef{3-25} gives, for suitably small $\ea_0$, that
\be\label{3-31}
  \mathcal{ E}_{j, i}(t)
\lesssim  \lt\{ \begin{split} & \mathcal{ E}_{j }(t)   + \mathcal{ E}_{j+1 }(t) + \sum_{ \iota \ge 0, \ p \ge 0, \ \iota+p \le  j  } \mathcal{E}_{ \iota,p } (t)
+  \sum_{\iota=0}^j\mathcal{E}_{ \iota} (t) , & i=1, \\
&\mathcal{ E}_{j, i-1}(t) +   \mathcal{ E}_{j+2, i-2}(t) + \mathcal{ E}_{j+1, i-2}(t) +\sum_{0\le \iota \le j, \ p \ge 0, \ \iota+p \le i+j-1 } \mathcal{E}_{ \iota,p } (t)  , \ \  &i\ge 2.
\end{split} \rt. \ee

Now, we use  estimate \ef{3-28}, derived from the induction hypothesis  \ef{3-26}, and \ef{3-31} to show that \ef{lem322est} holds for $i+j=k+1$. First, choosing $j=k$ and $i=1$ in \ef{3-31} gives
\be\label{3-32}
\mathcal{ E}_{k, 1}(t) \lesssim \sum_{\iota=0}^{k+1}\mathcal{E}_{ \iota} (t)   + \sum_{ \iota \ge 0, \ p \ge 0, \ \iota+p \le  k  } \mathcal{E}_{ \iota,p } (t) \lesssim \sum_{\iota=0}^{k+1}\mathcal{E}_{ \iota} (t).
\ee
We choose   $j=k-1$ and $i=2$ in \ef{3-31} and use \ef{3-27}-\ef{3-28} to show
\bee\label{}
\mathcal{ E}_{k-1, 2}(t) \lesssim   \mathcal{ E}_{k-1,  1}(t) +   \mathcal{ E}_{k+1, 0}(t) + \mathcal{ E}_{k, 0}(t) +\sum_{0\le \iota \le k-1, \ p \ge 0, \ \iota+p \le k } \mathcal{E}_{ \iota,p } (t) \lesssim \sum_{\iota=0}^{k+1}\mathcal{E}_{ \iota} (t).
\eee
For $\mathcal{ E}_{k-2, 3}$, it follows from \ef{3-31}, \ef{3-28} and \ef{3-32} to obtain
\bee\label{}
\mathcal{ E}_{k-2, 3}(t) \lesssim   \mathcal{ E}_{k-2,  2}(t) +   \mathcal{ E}_{k , 1}(t) + \mathcal{ E}_{k-1, 1}(t) +\sum_{0\le \iota \le k-2, \ p \ge 0, \ \iota+p \le k } \mathcal{E}_{ \iota,p } (t) \lesssim \sum_{\iota=0}^{k+1}\mathcal{E}_{ \iota} (t).
\eee
The other cases can be handled similarly. So we have proved \ef{lem322est} when $i+j=k+1$. This finishes the proof of Lemma \ref{lem322}. $\Box$

\subsection{Nonlinear weighted energy estimates}
In this section, we show that the weighted energy $\mathcal{E}_j(t)$ can be bounded by the initial date for all $t\in[0,T]$.
\begin{prop}\label{prop2} Suppose that \ef{apb} holds for suitably small positive number $\ea_0\in(0,1)$. Then it holds that for $t\in[0,T]$,
\be\label{prop2est}\begin{split}
  \mathcal{E}_j(t)    \lesssim   \sum_{\iota=0}^j \mathcal{E}_\iota (0) , \ \  j=0,1,\cdots, l.
\end{split}\ee
\end{prop}
The proof of this proposition consists of Lemma \ref{lem331} and Lemma \ref{lem332}.

\subsubsection{Basic energy estimates}
\begin{lem}\label{lem331}  Assume that \ef{apb} holds for suitably small positive number $\ea_0\in(0,1)$. Then,
\be\label{E0t}\begin{split}
  \mathcal{E}_0(t)   +
 \int_0^t  \int \lt[(1+s)^{-1} r^2 \bar\rho_0^\ga \lt(\zeta^2+ \lt(r\zeta_r\rt)^2\rt) +(1+s)  r^4 \bar\rho_0  \zeta_s^2  \rt] drds
 \lesssim    \mathcal{E}_0(0), \  \  t\in[0,T] .
\end{split}\ee
\end{lem}
{\bf Proof}. Multiplying \ef{pertb} by $r^3 \zeta_t$, and integrating the product with respect to the spatial variable, we obtain, using the integration by parts, that
\be\label{4.1}\begin{split}
  \f{d}{dt}\int \f{1}{2}r^4\bar\rho_0 \zeta_{t}^2 dr   +
  \int r^4 \bar\rho_0  \zeta_t^2 dr + \int  \bar\rho_0^\ga \mathfrak{L}_1 dr
   =0
\end{split}\ee
where
$$\mathfrak{L}_1:=-\lt(\tilde\eta_r+\zeta\rt)^{-2\ga}\lt(\tilde\eta_r+\zeta +r\zeta_r\rt)^{-\ga} \lt[r^3\lt(\tilde\eta_r+\zeta\rt)^2\zeta_t \rt]_r
+  \tilde\eta_r^{2-3\ga} \lt(r^3 \zeta_t\rt)_r=:-\mathfrak{L}_{11}+\mathfrak{L}_{12}.$$
For  $\mathfrak{L}_{11}$, note that
\bee\label{}\begin{split}
 &\lt[r^3\lt(\tilde\eta_r+\zeta\rt)^2\zeta_t \rt]_r \\
 =&3 r^2 \lt(\tilde\eta_r+\zeta\rt)^2\zeta_t + 2 r^2 \lt(\tilde\eta_r+\zeta\rt) \lt(r\zeta_r\rt) \zeta_t + r^2\lt(\tilde\eta_r+\zeta\rt)^2 \lt(r\zeta_{rt}\rt) \\
 =&
  2 r^2 \lt(\tilde\eta_r+\zeta\rt) \lt(\tilde\eta_r+\zeta +r\zeta_r\rt)  \zeta_t +r^2 \lt(\tilde\eta_r+\zeta\rt)^2 \lt(\zeta+ r\zeta_r\rt)_t,
\end{split}\eee
thus,
\bee\label{}\begin{split}
  \mathfrak{L}_{11}=& 2r^2\lt(\tilde\eta_r+\zeta\rt)^{1-2\ga }\lt(\tilde\eta_r+\zeta +r\zeta_r\rt)^{1-\ga }\zeta_t  +  r^2 \lt(\tilde\eta_r+\zeta\rt)^{2-2\ga }\lt(\tilde\eta_r+\zeta +r\zeta_r\rt)^{-\ga}\lt(\zeta+ r\zeta_r\rt)_t \\
 =& \f{r^2}{1-\ga} \lt[\lt(\tilde\eta_r+\zeta\rt)^{2-2\ga }\lt(\tilde\eta_r+\zeta +r\zeta_r\rt)^{1-\ga }\rt]_t \\
 &- r^2\lt[2\lt(\tilde\eta_r+\zeta\rt)^{1-2\ga }\lt(\tilde\eta_r+\zeta +r\zeta_r\rt)^{1-\ga }  +  \lt(\tilde\eta_r+\zeta\rt)^{2-2\ga }\lt(\tilde\eta_r+\zeta +r\zeta_r\rt)^{-\ga} \rt] \tilde\eta_{rt}.
\end{split}\eee
Clearly, $ \mathfrak{L}_{12}$ can be rewritten as
\bee\label{}\begin{split}
   \mathfrak{L}_{12} =   r^2  \lt( 3\zeta + r\zeta_r\rt)_t \tilde\eta_r^{2-3\ga}
  =r^2 \lt[  \lt( 3\zeta + r\zeta_r\rt) \tilde\eta_r^{2-3\ga}\rt]_t -(2 -3\ga)r^2  \lt( 3\zeta + r\zeta_r\rt) \tilde\eta_r^{1-3\ga}\tilde\eta_{rt}.
\end{split}\eee
Substitute these calculations into \ef{4.1} to give
\be\label{4.2}\begin{split}
 \f{d}{dt}\int \lt( \f{1}{2}r^4\bar\rho_0 \zeta_{t}^2 + r^2 \bar\rho_0^\ga \widetilde{\mathfrak{E}}_0 \rt) dr   +
  \int r^4 \bar\rho_0  \zeta_t^2 dr
 +\int r^2\bar\rho_0^\ga  \tilde\eta_{rt}   \mathfrak{F}   dr
   =0 ,
\end{split}\ee
where
\bee\label{}\begin{split}
  \widetilde{\mathfrak{E}}_0:=& \f{1 }{ \ga-1} \lt[\lt(\tilde\eta_r+\zeta\rt)^{2-2\ga }\lt(\tilde\eta_r+\zeta +r\zeta_r\rt)^{1-\ga } -\tilde\eta_r^{3-3\ga }  + (\ga-1) \lt( 3\zeta + r\zeta_r\rt) \tilde\eta_r^{2-3\ga}\rt]
,\\
  \mathfrak{F}: = & 2\lt(\tilde\eta_r+\zeta\rt)^{1-2\ga }\lt(\tilde\eta_r+\zeta +r\zeta_r\rt)^{1-\ga }  +  \lt(\tilde\eta_r+\zeta\rt)^{2-2\ga }\lt(\tilde\eta_r+\zeta +r\zeta_r\rt)^{-\ga} \\
  &-3 \tilde\eta_r^{2-3\ga}
   -(2-3\ga)  \lt( 3\zeta + r\zeta_r\rt) \tilde\eta_r^{1-3\ga}.
\end{split}\eee
It follows from the Taylor expansion, the smallness of $\zeta$ and $r\zeta_r$  which is a consequence of \ef{apb}, and \ef{decay}   that
\bee\label{}\begin{split}
  \widetilde{\mathfrak{E}}_0 = & \tilde\eta_r^{1-3\ga}
 \lt[\f{3}{2}(3\ga-2) \zeta^2 +  (3\ga-2)\zeta r \zeta_{r}
 + \f{\ga}{2} \lt( r   \zeta_r \rt)^2 \rt] \\
   & +O(1)  \tilde\eta_r^{ -3\ga} \lt(|\zeta|+|r\zeta_{r}|\rt)
   \lt( \zeta ^2 +\lt( r   \zeta_r \rt)^2\rt)   \\
    \sim &  \tilde\eta_r^{1-3\ga}  \lt( \zeta ^2 +\lt( r   \zeta_r \rt)^2\rt)
   \sim     (1+t)^{-1}\lt(\zeta^2+ \lt(r\zeta_r\rt)^2\rt),
 \end{split}\eee
\bee\label{}\begin{split}
 \mathfrak{F} \ge  & (3\ga-1) \tilde\eta_r^{ -3\ga}
 \lt[\f{3}{2}(3\ga-2) \zeta^2 +  (3\ga-2)\zeta r \zeta_{r}
 + \f{\ga}{2} \lt( r   \zeta_r \rt)^2 \rt] \\
  & -C  \tilde\eta_r^{ -3\ga-1} \lt(|\zeta|+|r\zeta_{r}|\rt)
   \lt( \zeta ^2 +\lt( r   \zeta_r \rt)^2\rt)   \\
    \ge&     (1+t)^{-\frac{3\ga }{3\ga-1}}\lt(\zeta^2+ \lt(r\zeta_r\rt)^2\rt) \ge 0.
 \end{split}\eee
Here and thereafter the notation $O(1)$ represents a finite number could be positive or negative.
We then  have, by integrating \ef{4.2} with respect to the temporal variable, that
\bee\label{}\begin{split}
\lt. \int \lt( \f{1}{2}r^4\bar\rho_0 \zeta_{t}^2 +r^2 \bar\rho_0^\ga  \widetilde{\mathfrak{E}}_0 \rt)(r,s) dr \rt|_{s=0}^{t}  +
  \int_0^t \int r^4 \bar\rho_0  \zeta_s^2 drds \le  0
\end{split}\eee
and
\be\label{4.4}\begin{split}
 & \int \lt[  r^4\bar\rho_0 \zeta_{t}^2 + (1+t)^{-1}r^2 \bar\rho_0^\ga \lt(\zeta^2+ \lt(r\zeta_r\rt)^2\rt) \rt](r,t) dr   +
  \int_0^t \int r^4 \bar\rho_0  \zeta_s^2 drds \\
  \lesssim &  \int \lt[  r^4\bar\rho_0 \zeta_{t}^2 +  r^2 \bar\rho_0^\ga \lt(\zeta^2+ \lt(r\zeta_r\rt)^2\rt) \rt](r,0) dr.
\end{split}\ee

Multiplying \ef{pertb} by $r^3 \zeta$, and integrating the product with respect to the spatial variable, we have, using the integration by parts, that
\be\label{4.5}\begin{split}
  \f{d}{dt}\int  r^4\bar\rho_0 \lt( \frac{1}{2}\zeta^2 + \zeta \zeta_t \rt) dr
   + \int  \bar\rho_0^\ga  \mathfrak{L}_2 dr
   = \int r^4 \bar\rho_0  \zeta_t^2 dr,
\end{split}\ee
where
$$ \mathfrak{L}_2:=-\lt(\tilde\eta_r+\zeta\rt)^{-2\ga}\lt(\tilde\eta_r+\zeta +r\zeta_r\rt)^{-\ga} \lt[r^3\lt(\tilde\eta_r+\zeta\rt)^2\zeta  \rt]_r
+  \tilde\eta_r^{2-3\ga} \lt(r^3 \zeta \rt)_r ,$$
which can be rewritten as
\bee\label{}\begin{split}
   \mathfrak{L}_2= & - 2r^2\lt(\tilde\eta_r+\zeta\rt)^{1-2\ga }\lt(\tilde\eta_r+\zeta +r\zeta_r\rt)^{1-\ga }\zeta  -  r^2 \lt(\tilde\eta_r+\zeta\rt)^{2-2\ga }\lt(\tilde\eta_r+\zeta +r\zeta_r\rt)^{-\ga}\\
  &\times \lt(\zeta+ r\zeta_r\rt)+ r^2 \tilde\eta_r^{2-3\ga}  \lt(3\zeta+ r\zeta_r\rt)\\
  =& r^2 \lt[3\tilde\eta_r^{2-3\ga}-2 \lt(\tilde\eta_r+\zeta\rt)^{1-2\ga }\lt(\tilde\eta_r+\zeta +r\zeta_r\rt)^{1-\ga }
  -\lt(\tilde\eta_r+\zeta\rt)^{2-2\ga }\lt(\tilde\eta_r+\zeta +r\zeta_r\rt)^{-\ga}\rt]\zeta\\
  &+r^2 \lt[\tilde\eta_r^{2-3\ga}-\lt(\tilde\eta_r+\zeta\rt)^{2-2\ga }\lt(\tilde\eta_r+\zeta +r\zeta_r\rt)^{-\ga}\rt]r\zeta_r .
\end{split}\eee
Again,  we use the Taylor expansion, \ef{basic} and  \ef{apb} to obtain
\bee\label{}\begin{split}
   \mathfrak{L}_2  \gtrsim  &   r^2  (1+t)^{-1}
 \lt[3(3\ga-2) \zeta^2 + 2 (3\ga-2)\zeta r \zeta_{r}
 + \ga \lt( r   \zeta_r \rt)^2 \rt] \\
 &  - C r^2  (1+t)^{-\frac{3\ga}{3\ga-1}} \lt(|\zeta|+|r\zeta_{r}|\rt)
   \lt( \zeta ^2 +\lt( r   \zeta_r \rt)^2\rt) \\
   \gtrsim  & r^2  (1+t)^{-1}
 \lt[3(3\ga-2) \zeta^2 + 2 (3\ga-2)\zeta r \zeta_{r}
 + \ga \lt( r   \zeta_r \rt)^2 -C \ea_0 \lt( \zeta ^2 +\lt( r   \zeta_r \rt)^2\rt) \rt] \\
 \gtrsim  & r^2 (1+t)^{-1} \lt(\zeta^2+ \lt(r\zeta_r\rt)^2\rt) ,
\end{split}\eee
provide that $\ea_0$ is suitably small.
It then follows from \ef{4.5}, the Cauchy inequality  and \ef{4.4} that
\be\label{4.6}\begin{split}
 & \int   \lt( r^4\bar\rho_0 \zeta ^2\rt)  (r,t) dr   +
  \int_0^t  \int (1+s)^{-1} r^2 \bar\rho_0^\ga \lt(\zeta^2+ \lt(r\zeta_r\rt)^2\rt)  drds \\
  \lesssim &  \int \lt( r^4\bar\rho_0 \lt(\zeta^2 + \zeta_{t}^2\rt) \rt)  (r,0) dr + \int   \lt( r^4\bar\rho_0 \zeta_t^2\rt)  (r,t) dr
  +
  \int_0^t \int r^4 \bar\rho_0  \zeta_s^2 drds \\
  \lesssim &  \int \lt[  r^4\bar\rho_0 \lt(\zeta^2 + \zeta_{t}^2\rt) +  r^2 \bar\rho_0^\ga \lt(\zeta^2+ \lt(r\zeta_r\rt)^2\rt) \rt](r,0) dr =\mathcal{E}_0(0).
\end{split}\ee

Next, we show the time decay of the energy norm. Multiply  equation \ef{4.2} by (1+t) and integrate the product with respect to the temporal variable to get
\bee\label{4.7}\begin{split}
&(1+t)\int \lt( \f{1}{2}r^4\bar\rho_0 \zeta_{t}^2 +r^2 \bar\rho_0^\ga  \widetilde{\mathfrak{E}}_0 \rt)(r,t) dr   +
 \int_0^t (1+s) \int r^4 \bar\rho_0  \zeta_s^2 drds \\
 \le& \int \lt( \f{1}{2}r^4\bar\rho_0 \zeta_{t}^2 +r^2 \bar\rho_0^\ga  \widetilde{\mathfrak{E}}_0 \rt)(r,0) dr + \int_0^t \int \lt( \f{1}{2}r^4\bar\rho_0 \zeta_{s}^2 + \widetilde{\mathfrak{E}}_0 \rt)  drds \\
 \lesssim & \int \lt( \f{1}{2}r^4\bar\rho_0 \zeta_{t}^2 + r^2 \bar\rho_0^\ga \widetilde{\mathfrak{E}}_0 \rt)(r,0) dr + \int_0^t \int \lt[ r^4\bar\rho_0 \zeta_{s}^2
 + (1+s)^{-1} r^2 \bar\rho_0^\ga \lt(\zeta^2+ \lt(r\zeta_r\rt)^2\rt)  \rt] drds
 \\
 \lesssim &  \int \lt[  r^4\bar\rho_0 \lt(\zeta^2 + \zeta_{t}^2\rt) +  r^2 \bar\rho_0^\ga \lt(\zeta^2+ \lt(r\zeta_r\rt)^2\rt) \rt](r,0) dr =\mathcal{E}_0(0),
\end{split}\eee
where estimates \ef{4.4} and \ef{4.6} have been used to derive the last inequality. This means
\bee\label{}\begin{split}
& \int \lt[(1+t) r^4\bar\rho_0 \zeta_{t}^2 +  r^2 \bar\rho_0^\ga \lt(\zeta^2+ \lt(r\zeta_r\rt)^2\rt) \rt](r,t) dr   +
 \int_0^t (1+s) \int r^4 \bar\rho_0  \zeta_s^2 drds \lesssim \mathcal{E}_0(0),
\end{split}\eee
which, together with \ef{4.6}, gives \ef{E0t}. This finishes the proof of Lemma \ref{lem331}. $\Box$

\subsubsection{Higher-order energy estimates}
Equation \ef{pertb} reads
\bee\label{}\begin{split}
   &r \bar\rho_0 \zeta_{tt}   +
   r  \bar\rho_0  \zeta_t
   +    \left[ \bar\rho_0^\ga     \lt(\tilde\eta_r+\zeta\rt)^{2-2\ga}\lt(\tilde\eta_r+\zeta +r\zeta_r\rt)^{-\ga}   \right]_r
   -    \tilde\eta_r^{2-3\ga}  \lt(   {\bar\rho_0}  ^\ga    \right)_r
   \\
   &-2 \bar\rho_0^\ga     \lt(\tilde\eta_r+\zeta\rt)^{1-2\ga}\lt(\tilde\eta_r+\zeta +r\zeta_r\rt)^{-\ga}  \zeta_r  =0 .
\end{split}\eee
Let $k\ge 1$ be an integer and take the $k$-th time derivative of the equation above. One has
 \be\label{4-21}\begin{split}
    r \bar\rho_0 \pl_t^k \zeta_{tt}   +
   r  \bar\rho_0 \pl_t^k  \zeta_{t}
  +  \lt[{\bar\rho_0}^\ga  \lt(w_1 \pl_t^k \zeta + w_2 r\pl_t^k \zeta_r + K_1 \rt) \rt]_r  + \bar\rho_0^\ga \lt[ (3w_2-w_1)  \pl_t^k \zeta_r + K_2 \rt]   \\
  -2 \bar\rho_0^\ga   \lt(w_3  \zeta_r  \pl_t^k \zeta +K_3\rt)  +\pl_t^{k-1} \lt\{ {\bar\rho_0}^\ga  \tilde\eta_{rt} \lt[ w_1 -  (2-3\ga)  \tilde\eta_r^{1-3\ga} \rt]\rt\}_r
 -2 \bar\rho_0^\ga  \pl_t^{k-1} \lt(\tilde\eta_{rt}   w_3     \zeta_r \rt)
    =0.
\end{split}\ee
Here
\bee\label{}\begin{split}
 w_1=&(2-2\ga) \lt(\tilde\eta_r+\zeta\rt)^{1-2\ga}\lt(\tilde\eta_r+\zeta +r\zeta_r\rt)^{-\ga}  -\ga  \lt(\tilde\eta_r+\zeta\rt)^{2-2\ga}\lt(\tilde\eta_r+\zeta +r\zeta_r\rt)^{-\ga-1}, \\
 w_2=&-\ga  \lt(\tilde\eta_r+\zeta\rt)^{2-2\ga}\lt(\tilde\eta_r+\zeta +r\zeta_r\rt)^{-\ga-1},
\\
 w_3=&(1-2\ga) \lt(\tilde\eta_r+\zeta\rt)^{ -2\ga}\lt(\tilde\eta_r+\zeta +r\zeta_r\rt)^{-\ga}  -\ga  \lt(\tilde\eta_r+\zeta\rt)^{1-2\ga}\lt(\tilde\eta_r+\zeta +r\zeta_r\rt)^{-\ga-1};
\end{split}\eee
and
\bee\label{}\begin{split}
 K_1=  & \pl_t^{k-1}  \lt(w_1 \zeta_t + w_2 r\zeta_{tr} \rt) - \lt(w_1 \pl_t^k \zeta + w_2 r\pl_t^k \zeta_r \rt),  \\
 K_2=  & \pl_t^{k-1} \lt[ (3w_2-w_1) \zeta_{tr}\rt] -  (3w_2-w_1) \pl_t^k \zeta_r
   ,\\
 K_3= &  \pl_t^{k-1}( w_3  \zeta_r  \zeta_t )-  w_3  \zeta_r  \pl_t^k \zeta
   .
\end{split}\eee
It should be noted that $K_1$, $K_2$ and $K_3$ contain lower-order terms involving  $\pl_t^{\iota} (\zeta, \zeta_r)  $ with $\iota=0,\cdots,k-1$; and $w_1$, $w_2$ and $w_3$ can be expanded, according to the Taylor expansion and the smallness of $\zeta$ and $r\zeta_r$ which is a consequence of \ef{apb},  as follows
\be\label{w123}\begin{split}
 w_1=& (2-3\ga) \tilde\eta_r^{1-3\ga } +(3\ga-1) \tilde\eta_r^{-3\ga }\lt[( 3\ga-2)\zeta +\ga r\zeta_r \rt] +  \bar{w}_1 \\
 w_2=&-\ga\tilde\eta_r^{1-3\ga } +\ga \tilde\eta_r^{-3\ga }\lt[( 3\ga-1)\zeta +(\ga+1) r\zeta_r \rt] + \bar{w}_2  ,
\\
 w_3=& (1-3\ga) \tilde\eta_r^{ -3\ga } + \bar{w}_3.
\end{split}\ee
Here $\bar w_i$ satisfies
\be\label{barw123}\begin{split}
|\bar{w}_1|+|\bar{w}_2| \lesssim \tilde\eta_r^{-3\ga-1} \lt(|\zeta|^2+|r\zeta_r|^2\rt)  ,
\ \ {\rm and} \ \  |\bar{w}_3|\lesssim \tilde\eta_r^{-3\ga-1} \lt(|\zeta| +|r\zeta_r| \rt) .
\end{split}\ee
In particular, $K_1=K_2=K_3=0$ when $k=1$.

\begin{lem}\label{lem332} Assume that \ef{apb} holds for suitably small positive number $\ea_0\in(0,1)$. Then
 for all $j=  1,   \cdots,l $,
\be\label{lemekt}\begin{split}
&  \mathcal{E}_j(t)   +
 \int_0^t  \int \lt[(1+s)^{2j-1}  r^2 \bar\rho_0^\ga \lt| \pl_s^j\lt( \zeta, r \zeta_r\rt)\rt|^2 +(1+s)^{2j+1}  r^4 \bar\rho_0 \lt(\pl_s^{j } \zeta_s \rt)^2  \rt] drds \\
 \lesssim  & \sum_{\iota=0}^j \mathcal{E}_\iota (0), \ \  t\in [0,T] .
\end{split}\ee
\end{lem}
{\bf Proof}. We use induction to prove \ef{lemekt}. As shown in Lemma \ref{lem331}, we know that \ef{lemekt} holds for $j=0$. For $1\le k\le l$, we make the induction hypothesis that  \ef{lemekt} holds for all $j=0,\cdots,k-1$, that is,
\be\label{t411}\begin{split}
  &  \mathcal{E}_j(t)   +
 \int_0^t  \int \lt[(1+s)^{2j-1}  r^2 \bar\rho_0^\ga \lt| \pl_s^j\lt( \zeta, r \zeta_r\rt)\rt|^2 +(1+s)^{2j+1}  r^4 \bar\rho_0 \lt(\pl_s^{j } \zeta_s \rt)^2  \rt] drds
\\
  \lesssim  &  \sum_{\iota=0}^j \mathcal{E}_\iota (0)  \ \ {\rm for} \ \ {\rm all} \ \  j=0, 1, \cdots, k-1.
\end{split}\ee
It suffices to prove \ef{lemekt} holds for $j=k$ under the induction hypothesis \ef{t411}.

{\bf Step 1}. In this step, we prove that
\be\label{3-62}\begin{split}
   & \frac{d}{dt}\lt[\int  \frac{1}{2}  r^4   \bar\rho_0
       \lt(\pl_t^k \zeta_t \rt)^2  dr + \overline{\mathfrak{E}}_k \rt] + \int  r^4  \bar\rho_0  \lt(\pl_t^k \zeta_t\rt)^2 dr
  \\     \lesssim &     (\ea_0+\da) (1+t)^{  -2k-2} \mathcal{E}_k(t) +   \lt(\ea_0+\da^{-1}\rt) (1+t)^{  -2k-2} \sum_{\iota=0}^{k-1}  \mathcal{E}_\iota(t) ,
\end{split}\ee
for any positive number $\da>0$ which will be specified later,
where
$\overline{\mathfrak{E}}_k:= \int r^2 \bar\rho_0^\ga \widetilde{{\mathfrak{E}}}_k   dr + M_2$. Here  $\widetilde{{\mathfrak{E}}}_k$ and $M_2$ are defined by \ef{defn331} and \ef{m2}, respectively. Moreover, we show that  $\overline{\mathfrak{E}}_k$ satisfies the following estimates:
\be\label{3-63}\begin{split}
    \overline{\mathfrak{E}}_k \ge C^{-1} (1+t)^{-1}\int r^2 \bar\rho_0^\ga   \lt| \pl_t^k \lt(\zeta  , r  \zeta_r\rt) \rt|^2  dr - C (1+t)^{-2k-1} \sum_{\iota=0}^{k-1}  \mathcal{E}_\iota(t) ,
\end{split}\ee
\be\label{3-64}\begin{split}
   \overline{\mathfrak{E}}_k \lesssim   (1+t)^{-1}\int r^2 \bar\rho_0^\ga   \lt| \pl_t^k \lt(\zeta  , r  \zeta_r\rt) \rt|^2  dr  +  (1+t)^{-2k-1} \sum_{\iota=0}^{k-1}  \mathcal{E}_\iota(t)  .
\end{split}\ee

We start with integrating the production of
 \ef{4-21} and $r^3 \pl_t^k \zeta_t$   with respect to the spatial variable which gives
\be\label{t414}\begin{split}
   & \frac{d}{dt}\int \frac{1}{2}  r^4   \bar\rho_0
       \lt(\pl_t^k \zeta_t \rt)^2    dr + \int  r^4  \bar\rho_0  \lt(\pl_t^k \zeta_t\rt)^2 dr + N_1+ N_2 =0 ,
\end{split}\ee
where
\bee\label{}\begin{split}
    N_1: = &  - \int {\bar\rho_0}^\ga \lt(w_1 \pl_t^k \zeta + w_2 r\pl_t^k \zeta_r \rt)  \lt(r^3 \pl_t^k \zeta_t \rt)_r dr + \int r^2 \bar\rho_0^\ga  (3w_2-w_1)\lt( r \pl_t^k \zeta_r \rt)\pl_t^k \zeta_t  dr \\
 &-2 \int r^2 \bar\rho_0^\ga    w_3  \lt(r\zeta_r \rt) \lt(\pl_t^k \zeta\rt) \pl_t^k \zeta_t   dr, \end{split}\eee
 \bee\label{}\begin{split}
      N_2 :=  & - \int  {\bar\rho_0}^\ga \lt\{K_1+ \pl_t^{k-1}\lt[ \tilde\eta_{rt} \lt( w_1 -  (2-3\ga)  \tilde\eta_r^{1-3\ga} \rt ) \rt] \rt\}\lt(r^3  \pl_t^k \zeta_t  \rt)_r dr \\
    &
    + \int r^2 \bar\rho_0^\ga \lt(\pl_t^k \zeta_t \rt) \lt[ r(K_2-2K_3) -2 \pl_t^{k-1} \lt(\tilde\eta_{rt}   w_3    r \zeta_r \rt)\rt]     dr.
\end{split}\eee
Note that $N_1$ and $N_2$ can be rewritten as
\bee\label{}\begin{split}
  N_1 =  -\frac{1}{2}\int r^2 \bar\rho_0^\ga  \lt[\lt(3w_1+2w_3r\zeta_r\rt)\lt[\lt(\pl_t^k \zeta  \rt)^2\rt]_t +2 w_1 \lt[\lt(\pl_t^k \zeta  \rt) r\pl_t^k \zeta_{r}  \rt]_t
 +  w_2\lt[ \lt( r \pl_t^k \zeta_r \rt)^2\rt]_t \rt] dr \\
 =  -\frac{1}{2}\frac{d}{dt}\int r^2 \bar\rho_0^\ga  \lt[\lt(3w_1+2w_3r\zeta_r\rt)\lt(\pl_t^k \zeta  \rt)^2 + 2w_1 \lt(\pl_t^k \zeta  \rt) r\pl_t^k \zeta_{r}
 +  w_2 \lt( r \pl_t^k \zeta_r \rt)^2 \rt] dr + \widetilde{N}_1,
\end{split}\eee
\bee\label{}\begin{split}
  N_2=  & \frac{d}{dt} M_2 + \int  {\bar\rho_0}^\ga \lt\{K_{1t}+ \pl_t^{k}\lt[ \tilde\eta_{rt} \lt( w_1 -  (2-3\ga)  \tilde\eta_r^{1-3\ga} \rt ) \rt] \rt\}  \lt(r^3  \pl_t^k \zeta   \rt)_r dr \\
    &
    - \int r^2 \bar\rho_0^\ga \lt(\pl_t^k \zeta  \rt) \lt[ r(K_2-2K_3)_t -2 \pl_t^{k } \lt(\tilde\eta_{rt}   w_3    r \zeta_r \rt)\rt]     dr
    = : \frac{d}{dt} M_2 + \widetilde{N}_2,
\end{split}\eee
where
\bee\label{}\begin{split}
  \widetilde{N}_1  : =&  \frac{1}{2} \int r^2 \bar\rho_0^\ga  \lt[\lt(3w_1+ 2w_3r\zeta_r\rt)_t \lt(\pl_t^k \zeta  \rt)^2 + 2w_{1t} \lt(\pl_t^k \zeta  \rt) r\pl_t^k \zeta_{r}
 +  w_{2t} \lt( r \pl_t^k \zeta_r \rt)^2 \rt] dr ,
\end{split}\eee
\be\label{m2}\begin{split}
  M_2 : = &   -\int  {\bar\rho_0}^\ga \lt\{K_1+ \pl_t^{k-1}\lt[ \tilde\eta_{rt} \lt( w_1 -  (2-3\ga)  \tilde\eta_r^{1-3\ga} \rt ) \rt] \rt\}\lt(r^3  \pl_t^k \zeta  \rt)_r dr \\
    &
    + \int r^2 \bar\rho_0^\ga \lt(\pl_t^k \zeta\rt) \lt[ r(K_2-2K_3) -2 \pl_t^{k-1} \lt(\tilde\eta_{rt}   w_3    r \zeta_r \rt)\rt]     dr.
\end{split}\ee
It then follows from equation \ef{t414} that
\be\label{t415}\begin{split}
   & \frac{d}{dt}\lt[\int \lt( \frac{1}{2}  r^4   \bar\rho_0
       \lt(\pl_t^k \zeta_t \rt)^2 + r^2 \bar\rho_0^\ga  \widetilde{\mathfrak{E}}_k   \rt)dr + M_2\rt] + \int  r^4  \bar\rho_0  \lt(\pl_t^k \zeta_t\rt)^2 dr = -\widetilde{N}_1 -  \widetilde{N}_2,
\end{split}\ee
where
\be\label{defn331}\begin{split}
 \widetilde{\mathfrak{E}}_k : = -\frac{1}{2}\lt[\lt(3w_1+2w_3r\zeta_r\rt)\lt(\pl_t^k \zeta  \rt)^2 + 2w_1 \lt(\pl_t^k \zeta  \rt) r\pl_t^k \zeta_{r}
 +  w_2 \lt( r \pl_t^k \zeta_r \rt)^2 \rt]
\end{split}\ee
which satisfies
\be\label{t416}\begin{split}
 \widetilde{\mathfrak{E}}_k = & \tilde\eta_r^{1-3\ga}
 \lt[\frac{3}{2}(3\ga-2)\lt(\pl_t^k \zeta  \rt)^2 +  (3\ga-2)\lt(\pl_t^k \zeta  \rt) r\pl_t^k \zeta_{r}
 + \frac{\ga}{2}\lt( r \pl_t^k \zeta_r \rt)^2 \rt] \\
 &   +O(1) \tilde\eta_r^{ -3\ga}\lt(|\zeta|+|r\zeta_{r}|\rt)
   \lt(\lt(\pl_t^k \zeta  \rt)^2 +\lt( r \pl_t^k \zeta_r \rt)^2\rt)  \\
 \sim & \tilde\eta_r^{1-3\ga}  \lt[\lt(\pl_t^k \zeta  \rt)^2 +\lt( r \pl_t^k \zeta_r \rt)^2\rt]
 \sim  (1+t)^{-1}\lt[\lt(\pl_t^k \zeta  \rt)^2 +\lt( r \pl_t^k \zeta_r \rt)^2 \rt].
\end{split}\ee
 Here we have used \ef{w123}, \ef{decay} and the smallness of $\zeta$ and $r\zeta_r$  which is a consequence of \ef{apb} to derive the above equivalence.  We will show later that $M_2$ can be bounded by the integral of $\widetilde{\mathfrak{E}}_k$ and lower-order terms, see \ef{3-67}.

In what follows, we  analyze the terms on the right-hand side of \ef{t415}.    Clearly, $-\widetilde{N}_1$ can be bounded  by
\bee\label{}\begin{split}
 - \widetilde{N}_1  \le  &  (1- 3\ga)\int r^2 \bar\rho_0^\ga \tilde\eta_r^{-3\ga}\tilde\eta_{rt}
   \lt[\lt(\frac{9}{2}\ga-3\rt)\lt(\pl_t^k \zeta  \rt)^2  +  (3\ga-2)\lt(\pl_t^k \zeta  \rt) \lt(r\pl_t^k \zeta_r  \rt) \rt.\\
&\lt. + \frac{\ga}{2}\lt(r\pl_t^k \zeta_r  \rt)^2 \rt]  dr
  + C\int r^2 \bar\rho_0^\ga  \tilde\eta_r^{-3\ga }\lt[\tilde\eta_r^{ -1}\tilde\eta_{rt} \lt(|\zeta|+|r\zeta_{r}|\rt) \rt.\\
&\lt.   +\lt(1+\tilde\eta_r^{ -1} \lt(|\zeta|+|r\zeta_{r}|\rt) \rt) \lt(|\zeta_t|+|r\zeta_{rt}|\rt)\rt]  \lt(\lt(\pl_t^k \zeta  \rt)^2 +\lt( r \pl_t^k \zeta_r \rt)^2\rt)  dr.
\end{split}\eee
It should be noted that the first integral on the right-hand side of the inequality above is non-positive due to $\tilde\eta_{rt} \ge 0$. Thus, we have by use of \ef{decay} and \ef{apb} that
\be\label{t418}\begin{split}
 - \widetilde{N}_1  \lesssim   \ea_0 (1+t)^{-2-\frac{1}{3\ga-1}} \int r^2 \bar\rho_0^\ga  \lt(\lt(\pl_t^k \zeta  \rt)^2 +\lt( r \pl_t^k \zeta_r \rt)^2\rt)  dr.
\end{split}\ee
To control $\widetilde{N}_2$, we may rewrite it as
\bee\label{}\begin{split}
  \widetilde{N}_2=  &
    \int r^2 {\bar\rho_0}^\ga \lt\{  \lt(3 \pl_t^k \zeta +   r\pl_t^k \zeta_r  \rt) \pl_t^{k}\lt[ \tilde\eta_{rt} \lt( w_1 -  (2-3\ga)  \tilde\eta_r^{1-3\ga} \rt ) \rt]  +  2 \lt(\pl_t^k \zeta  \rt)  \pl_t^{k } \lt(\tilde\eta_{rt}   w_3    r \zeta_r \rt)  \rt\}  dr \\
 &  +  \int r^2  {\bar\rho_0}^\ga  \lt[ K_{1t}   \lt(3 \pl_t^k \zeta +   r\pl_t^k \zeta_r  \rt)   -  r(K_2-2K_3)_t \lt(\pl_t^k \zeta  \rt) \rt]dr
    =:  \widetilde{N}_{21} + \widetilde{N}_{22}.
\end{split}\eee
For $\widetilde{N}_{21}$, note that
 \bee\label{}\begin{split}
    & \pl_t^{k}\lt[ \tilde\eta_{rt} \lt( w_1 -  (2-3\ga)  \tilde\eta_r^{1-3\ga} \rt ) \rt]
     \\
     =  & (3\ga-1)  \pl_t^{k}\lt[ \tilde\eta_{rt} \tilde\eta_r^{-3\ga }\lt(( 3\ga-2)\zeta +\ga r\zeta_r \rt)    \rt] + \pl_t^{k}\lt(\tilde\eta_{rt}\bar{w}_1\rt)\\
  = & (3\ga-1)\tilde\eta_{rt} \tilde\eta_r^{-3\ga }\lt(( 3\ga-2)\pl_t^{k}\zeta +\ga r\pl_t^{k}\zeta_r \rt) \\
  & + O(1) \sum_{\iota=1}^{k}\lt|\pl_t^{\iota}\lt(\tilde\eta_{rt} \tilde\eta_r^{-3\ga }\rt) \rt|\lt|\pl_t^{k-\iota}\lt(\zeta, r \zeta_r\rt)\rt| +  \pl_t^{k}\lt(\tilde\eta_{rt}\bar{w}_1\rt).
\end{split}\eee
and
 \bee\label{}\begin{split}
  &  \pl_t^{k } \lt(\tilde\eta_{rt}   w_3    r \zeta_r \rt)=   (1-3\ga)\pl_t^{k } \lt(\tilde\eta_{rt} \tilde\eta_r^{ -3\ga }    r \zeta_r \rt)  + \pl_t^{k } \lt(\tilde\eta_{rt}  \bar{w}_3    r \zeta_r \rt)\\
    &= (1-3\ga)\tilde\eta_{rt} \tilde\eta_r^{-3\ga }\lt( r\pl_t^{k}\zeta_r \rt)  + O(1) \sum_{\iota=1}^{k}\lt|\pl_t^{\iota}\lt(\tilde\eta_{rt} \tilde\eta_r^{-3\ga }\rt) \rt|\lt|r \pl_t^{k-\iota}   \zeta_r \rt|  + \pl_t^{k } \lt(\tilde\eta_{rt}  \bar{w}_3    r \zeta_r \rt).
\end{split}\eee
 Thus,
\be\label{t419}\begin{split}
 -\widetilde{N}_{21}  \le  &    (1- 3\ga)\int r^2 \bar\rho_0^\ga \tilde\eta_r^{-3\ga}\tilde\eta_{rt}
   \lt[\lt( {9} \ga- 6 \rt)\lt(\pl_t^k \zeta  \rt)^2  +  (6\ga-4)\lt(\pl_t^k \zeta  \rt) \lt(r\pl_t^k \zeta_r  \rt) \rt.\\
&  \lt.
 +  {\ga} \lt(r\pl_t^k \zeta_r  \rt)^2 \rt]dr
  +C\int r^2 \bar\rho_0^\ga   \lt(\lt|\pl_t^k \zeta  \rt| +\lt| r \pl_t^k \zeta_r \rt| \rt) \\
 &  \times \lt[\sum_{\iota=1}^{k}\lt|\pl_t^{\iota}\lt(\tilde\eta_{rt} \tilde\eta_r^{-3\ga }\rt) \rt|\lt|\pl_t^{k-\iota}\lt(\zeta, r \zeta_r\rt)\rt| + \lt| \pl_t^{k}\lt(\tilde\eta_{rt}\bar{w}_1, \tilde\eta_{rt}  \bar{w}_3    r \zeta_r \rt)\rt|  \rt]dr.
\end{split}\ee
For $\widetilde{N}_{22}$,  note that
 \bee\label{}\begin{split}
     K_{1t} = &(k-1)\lt(w_{1t}\pl_t^k\zeta+ w_{2t} r\pl_t^k \zeta_r \rt)+ O(1) \sum_{\iota=2}^{k}\lt|\pl_t^{\iota}\lt(w_1,w_2\rt) \rt|\lt|\pl_t^{k+1-\iota}\lt(\zeta, r \zeta_r\rt)\rt|,\\
    r K_{2t} = & (k-1)\lt(3w_2-w_1  \rt)_t \lt(r\pl_t^k \zeta_r \rt) + O(1) \sum_{\iota=2}^{k}\lt|\pl_t^{\iota}\lt(w_1,w_2\rt) \rt|\lt| \pl_t^{k+1-\iota} \lt( r \zeta_r \rt) \rt|,\\
    r K_{3t} = & (k-1)\lt( w_3 r\zeta_r  \rt)_t \lt( \pl_t^k \zeta  \rt) + O(1) \sum_{\iota=2}^{k}\lt|\pl_t^{\iota}\lt( w_3 r\zeta_r \rt) \rt|\lt| \pl_t^{k+1-\iota} \zeta  \rt|.
\end{split}\eee
Thus,
\be\label{t420}\begin{split}
   \widetilde{N}_{22}  \ge  & 2(k-1) \widetilde{N}_{1}
    -C \int r^2 \bar\rho_0^\ga   \lt(\lt|\pl_t^k \zeta  \rt| +\lt| r \pl_t^k \zeta_r \rt| \rt)  \\
        & \times \sum_{\iota=2}^{k}\lt|\pl_t^{\iota}\lt(w_1,w_2, w_3 r\zeta_r\rt) \rt|\lt|\pl_t^{k+1-\iota}\lt(\zeta, r \zeta_r\rt)\rt|   dr .
\end{split}\ee
In a similar way to dealing with  $\widetilde{N}_{1}$ shown in \ef{t418}, we have, with the aid of \ef{t419} and \ef{t420}, that
\be\label{t421}\begin{split}
  - \widetilde{N}_2  \lesssim   &  \ea_0 (1+t)^{-2-\frac{1}{3\ga-1}} \int r^2 \bar\rho_0^\ga  \lt(\lt(\pl_t^k \zeta  \rt)^2 +\lt( r \pl_t^k \zeta_r \rt)^2\rt)  dr \\
   &+ \int r^2 \bar\rho_0^\ga   \lt(\lt|\pl_t^k \zeta  \rt| +\lt| r \pl_t^k \zeta_r \rt| \rt) Q  dr.
\end{split}\ee
where
\be\label{defnQ}\begin{split}
Q =  &  \sum_{\iota=1}^{k}\lt|\pl_t^{\iota}\lt(\tilde\eta_{rt} \tilde\eta_r^{-3\ga }\rt) \rt|\lt|\pl_t^{k-\iota}\lt(\zeta, r \zeta_r\rt)\rt|      + \lt| \pl_t^{k}\lt(\tilde\eta_{rt}\bar{w}_1, \tilde\eta_{rt}  \bar{w}_3    r \zeta_r \rt)\rt|\\
& + \sum_{\iota=2}^{k}\lt|\pl_t^{\iota}\lt(w_1,w_2, w_3 r\zeta_r\rt) \rt|\lt|\pl_t^{k+1-\iota}\lt(\zeta, r \zeta_r\rt)\rt|  .
\end{split}\ee
Therefore, it produces from \ef{t415}, \ef{t418} and \ef{t421} that
\be\label{t422}\begin{split}
   & \frac{d}{dt}\lt[\int \lt( \frac{1}{2}  r^4   \bar\rho_0
       \lt(\pl_t^k \zeta_t \rt)^2 + r^2 \bar\rho_0^\ga \widetilde{\mathfrak{E}}_k   \rt)dr + M_2\rt] + \int  r^4  \bar\rho_0  \lt(\pl_t^k \zeta_t\rt)^2 dr
  \\     \lesssim &  \ea_0 (1+t)^{-2-\frac{1}{3\ga-1}} \int r^2 \bar\rho_0^\ga  \lt| \pl_t^k \lt(\zeta  , r  \zeta_r\rt) \rt|^2  dr + \int r^2 \bar\rho_0^\ga   \lt| \pl_t^k \lt(\zeta  , r  \zeta_r\rt) \rt| Q dr.
\end{split}\ee

We are to bound the last term on the right-hand side of \ef{t422}. It follows from \ef{decay} and \ef{apb} that
\be\label{qbounds}\begin{split}
Q \lesssim  & \ea_0 (1+t)^{-2-\frac{1}{3\ga-1}}  \lt|\pl_t^{k}\lt(\zeta, r \zeta_r\rt)\rt|
+\widetilde{Q},
\end{split}\ee
where
\begin{align}\label{3-57}
&\widetilde{Q} : =      \sum_{\iota=1}^{k} (1+t)^{-2-\iota} \lt|\pl_t^{k-\iota}\lt(\zeta, r \zeta_r\rt)\rt|
 + (1+t)^{-1-\frac{1}{3\ga-1}}\lt|\pl_t^{2}\lt(\zeta, r \zeta_r\rt)\rt| \lt|\pl_t^{k-1}\lt(\zeta, r \zeta_r\rt)\rt| \notag \\
&  + (1+t)^{-2-\frac{1}{3\ga-1}}\lt|\pl_t^{2}\lt(\zeta, r \zeta_r\rt)\rt|\lt|\pl_t^{k-2}\lt(\zeta, r \zeta_r\rt)\rt| + (1+t)^{-1-\frac{1}{3\ga-1}}\lt|\pl_t^{3}\lt(\zeta, r \zeta_r\rt)\rt|\lt|\pl_t^{k-2}\lt(\zeta, r \zeta_r\rt)\rt|
\notag \\
& + \lt[(1+t)^{-1-\frac{1}{3\ga-1}}\lt|\pl_t^{4}\lt(\zeta, r \zeta_r\rt)\rt|
+(1+t)^{-2-\frac{1}{3\ga-1}}\lt|\pl_t^{3}\lt(\zeta, r \zeta_r\rt)\rt|
 + (1+t)^{-3-\frac{1}{3\ga-1}}\rt. \notag\\
& \lt. \times \lt|\pl_t^{2}\lt(\zeta, r \zeta_r\rt)\rt|
+(1+t)^{-1-\frac{2}{3\ga-1}}\lt|\pl_t^{2}\lt(\zeta, r \zeta_r\rt)\rt|^2 \rt]\lt|\pl_t^{k-3}\lt(\zeta, r \zeta_r\rt)\rt| + {\rm l.o.t.}.
\end{align}
Here and thereafter the notation ${\rm l.o.t.}$ is used to represent the lower-order terms involving $\pl_t^{\iota}\lt(\zeta, r \zeta_r\rt) $ with $\iota=2,\cdots, k-4$. It should be noticed that the second term on the right-hand side of \ef{3-57} only appears as $k-1\ge 2$, the third term as $k-2\ge 2$, the fourth term as $k-2\ge 3$, and so on. Clearly, we use \ef{apb} again to obtain
\begin{align}\label{}
&\widetilde{Q} \lesssim     \sum_{\iota=1}^{k} (1+t)^{-2-\iota} \lt|\pl_t^{k-\iota}\lt(\zeta, r \zeta_r\rt)\rt|
 + \ea_0 \sigma^{-\frac{1}{2}}(1+t)^{-3-\frac{1}{3\ga-1}} \lt|\pl_t^{k-1}\lt(\zeta, r \zeta_r\rt)\rt| \notag \\
&    +\ea_0 \sigma^{-1} (1+t)^{-4-\frac{1}{3\ga-1}} \lt|\pl_t^{k-2}\lt(\zeta, r \zeta_r\rt)\rt|
  +   \ea_0 \sigma^{-\frac{3}{2}} (1+t)^{-5-\frac{1}{3\ga-1}} \lt|\pl_t^{k-3}\lt(\zeta, r \zeta_r\rt)\rt| + {\rm l.o.t.},\notag
\end{align}
if $k \ge 7$. Similarly, we can bound ${\rm l.o.t.}$ and achieve
\begin{align}\label{}
&\widetilde{Q} \lesssim     \sum_{\iota=1}^{k} (1+t)^{-2-\iota} \lt|\pl_t^{k-\iota}\lt(\zeta, r \zeta_r\rt)\rt|
 + \ea_0 \sum_{\iota=1}^{[(k-1)/2]}\sigma^{-\frac{\iota}{2}}(1+t)^{-2-\iota-\frac{1}{3\ga-1}} \lt|\pl_t^{k-\iota}\lt(\zeta, r \zeta_r\rt)\rt|,\notag
\end{align}
which  implies
\begin{align}\label{new3-58}
 \int r^2 \bar\rho_0^\ga   \lt| \pl_t^k \lt(\zeta  , r  \zeta_r\rt) \rt| \widetilde{Q}  dr
\lesssim  \sum_{\iota=1}^{k} (1+t)^{-2-\iota} \int r^2 \bar\rho_0^\ga   \lt| \pl_t^k \lt(\zeta  , r  \zeta_r\rt) \rt| \lt|\pl_t^{k-\iota}\lt(\zeta, r \zeta_r\rt)\rt| dr \notag
 \\+ \ea_0 \sum_{\iota=1}^{[(k-1)/2]} (1+t)^{-2-\iota }  \int r^2 \bar\rho_0^\ga  \sigma^{-\frac{\iota}{2}}   \lt| \pl_t^k \lt(\zeta  , r  \zeta_r\rt) \rt| \lt|\pl_t^{k-\iota}\lt(\zeta, r \zeta_r\rt)\rt| dr =:   \widetilde{\mathcal{Q}}_1 +  \widetilde{\mathcal{Q}}_2.
\end{align}
Easily, it follows from the Cauchy-Schwarz inequality that for any $\delta>0$,
\begin{align}\label{3-58}
 \widetilde{\mathcal{Q}}_1 \lesssim & \da (1+t)^{-2 } \int r^2 \bar\rho_0^\ga   \lt| \pl_t^k \lt(\zeta  , r  \zeta_r\rt) \rt|^2  dr  +
   \da^{-1} \sum_{\iota=1}^{k} (1+t)^{-2-2\iota} \int r^2 \bar\rho_0^\ga     \lt|\pl_t^{k-\iota}\lt(\zeta, r \zeta_r\rt)\rt|^2  dr
 \notag\\
 \lesssim & \da (1+t)^{-2 -2k }  \mathcal{E}_k(t) + \da^{-1} (1+t)^{-2-2k}\sum_{\iota=0}^{k-1}  \mathcal{E}_\iota (t)
\end{align}
and
\begin{align}\label{3-59}
 \widetilde{\mathcal{Q}}_2 \lesssim &  \ea_0 (1+t)^{-2 } \int r^2 \bar\rho_0^\ga   \lt| \pl_t^k \lt(\zeta  , r  \zeta_r\rt) \rt|^2  dr \notag\\
 & +
   \ea_0 \sum_{\iota=1}^{[(k-1)/2]} (1+t)^{-2-2\iota} \int r^2 \bar\rho_0^\ga \sigma^{-\iota}    \lt|\pl_t^{k-\iota}\lt(\zeta, r \zeta_r\rt)\rt|^2  dr   .
\end{align}
In view of  the Hardy inequality \ef{hardy2}, we see that for $\iota=1,\cdots,[(k-1)/2]$,
\begin{align}\label{}
 &  \int_{\mathcal{I}_b}   \sigma^{\aa+1 -\iota}    \lt|\pl_t^{k-\iota}\lt(\zeta,   \zeta_r\rt)\rt|^2  dr
 \lesssim \int_{\mathcal{I}_b}   \sigma^{\aa+3  -\iota}     \lt|\pl_t^{k-\iota}\lt(\zeta,   \zeta_r, \zeta_{rr}\rt)\rt|^2  dr \lesssim \cdots  \notag\\
 \lesssim &\sum_{i=0}^{\iota+1} \int_{\mathcal{I}_b}   \sigma^{\aa+1 +\iota}     \lt|\pl_t^{k-\iota} \pl_r^{i} \zeta \rt|^2  dr \lesssim  \sum_{i=0}^{\iota+1} \int_{\mathcal{I}_b}  r^4 \sigma^{\aa+1 +\iota}     \lt|\pl_t^{k-\iota} \pl_r^{i} \zeta \rt|^2  dr\notag\\
 \lesssim & (1+t)^{2\iota-2k} \lt( \mathcal{E}_{k-\iota } + \sum_{i=1}^{\iota}  \mathcal{E}_{k-\iota, i}\rt)(t)   \lesssim (1+t)^{2\iota-2k}  \sum_{\iota=0}^{k}  \mathcal{E}_\iota(t), \notag
 \end{align}
due to   $\aa+1 -\iota \ge \aa - [([\aa]+1)/2] \ge 0$ for $k\le l$, which ensures the application of the Hardy inequality. Here the last inequality follows from the elliptic estimate \ef{prop1est}. Thus, we can obtain for $\iota=1,\cdots,[(k-1)/2]$,
\begin{align}\label{new3-64}
 &\int r^2 \bar\rho_0^\ga \sigma^{-\iota}    \lt|\pl_t^{k-\iota}\lt(\zeta, r \zeta_r\rt)\rt|^2  dr =  \int r^2  \sigma^{\aa+1 -\iota}    \lt|\pl_t^{k-\iota}\lt(\zeta, r \zeta_r\rt)\rt|^2  dr \notag\\
 = & \int_{\mathcal{I}_o} r^2  \sigma^{\aa+1 -\iota}    \lt|\pl_t^{k-\iota}\lt(\zeta, r \zeta_r\rt)\rt|^2  dr
 + \int_{\mathcal{I}_b} r^2  \sigma^{\aa+1 -\iota}    \lt|\pl_t^{k-\iota}\lt(\zeta, r \zeta_r\rt)\rt|^2  dr \notag\\
 \lesssim &  \int_{\mathcal{I}_o} r^2  \sigma^{\aa+1 }    \lt|\pl_t^{k-\iota}\lt(\zeta, r \zeta_r\rt)\rt|^2  dr
 + \int_{\mathcal{I}_b}   \sigma^{\aa+1 -\iota}    \lt|\pl_t^{k-\iota}\lt(\zeta,  \zeta_r\rt)\rt|^2  dr  \notag\\
  \lesssim & (1+t)^{2\iota-2k}  \sum_{\iota=0}^{k}  \mathcal{E}_\iota(t) .
 \end{align}
This, together with \ef{3-59}, implies
\begin{align}\label{3-60}
 \widetilde{\mathcal{Q}}_2
 \lesssim   \ea_0 (1+t)^{-2-2k}\sum_{\iota=0}^{k }  \mathcal{E}_\iota (t)  .
\end{align}
So, it yields from \ef{qbounds}, \ef{new3-58}, \ef{3-58} and \ef{3-60} that for $\da>0$,
\begin{align}\label{}
 \int r^2 \bar\rho_0^\ga   \lt| \pl_t^k \lt(\zeta  , r  \zeta_r\rt) \rt|  {Q}  dr
\lesssim    (1+t)^{  -2k-2} \lt[  (\ea_0+\da) \mathcal{E}_k(t) +   \lt(\ea_0+\da^{-1}\rt)  \sum_{\iota=0}^{k-1}  \mathcal{E}_\iota(t) \rt]. \notag
\end{align}
Substitute this into \ef{t422} to give \ef{3-62}.

To prove \ef{3-63} and \ef{3-64}, we adopt a similar but much easier way to dealing with $\widetilde{N}_2$ as shown in \ef{t421} to show
\be\label{3-65}\begin{split}
  |M_2| \lesssim  &   \int r^2 \bar\rho_0^\ga   \lt(\lt|\pl_t^k \zeta  \rt| +\lt| r \pl_t^k \zeta_r \rt| \rt) P dr ,
\end{split}\ee
where
\bee\label{}\begin{split}
  P=  &   \sum_{\iota=0}^{k-1}\lt|\pl_t^{\iota}\lt(\tilde\eta_{rt} \tilde\eta_r^{-3\ga }\rt) \rt|\lt|\pl_t^{k-1-\iota}\lt(\zeta, r \zeta_r\rt)\rt|    + \lt| \pl_t^{k-1}\lt(\tilde\eta_{rt}\bar{w}_1, \tilde\eta_{rt}  \bar{w}_3    r \zeta_r \rt)\rt| \\
   &+ \sum_{\iota=1}^{k-1}\lt|\pl_t^{\iota}\lt(w_1,w_2, w_3 r\zeta_r\rt) \rt|\lt|\pl_t^{k -\iota}\lt(\zeta, r \zeta_r\rt)\rt|.
\end{split}\eee
In view of \ef{decay} and \ef{apb}, we have
\bee\label{}\begin{split}
& P \lesssim     \sum_{\iota=0}^{k-1} (1+t)^{-2-\iota} \lt|\pl_t^{k-1-\iota}\lt(\zeta, r \zeta_r\rt)\rt|
 + \lt|\pl_t^{k-2}\lt(\zeta, r \zeta_r\rt)\rt|(1+t)^{-1-\frac{1}{3\ga-1}}\lt|\pl_t^{2}\lt(\zeta, r \zeta_r\rt)\rt| \\
& + \lt|\pl_t^{k-3}\lt(\zeta, r \zeta_r\rt)\rt|\lt[(1+t)^{-1-\frac{1}{3\ga-1}}\lt|\pl_t^{3}\lt(\zeta, r \zeta_r\rt)\rt|
+(1+t)^{-2-\frac{1}{3\ga-1}}\lt|\pl_t^{2}\lt(\zeta, r \zeta_r\rt)\rt|\rt] + {\rm l.o.t.},
\end{split}\eee
which implies
\begin{align}\label{}
& P  \lesssim     \sum_{\iota=0}^{k-1} (1+t)^{-2-\iota} \lt|\pl_t^{k-1-\iota}\lt(\zeta, r \zeta_r\rt)\rt|
 + \ea_0 \sum_{\iota=2}^{[k/2]}\sigma^{ \frac{1-\iota}{2}}(1+t)^{-1-\iota-\frac{1}{3\ga-1}} \lt|\pl_t^{k-\iota}\lt(\zeta, r \zeta_r\rt)\rt|.\notag
\end{align}
Similar to the derivation of \ef{new3-64}, we can use the Hardy inequality \ef{hardy2} and  elliptic estimate \ef{prop1est} to obtain
\begin{align}\label{}
   \int r^2 \bar\rho_0^\ga   P^2  dr \lesssim   &  (1+t)^{-2-2k} \sum_{\iota=0}^{k-1}  \mathcal{E}_\iota(t)
 + \ea_0^2 \sum_{\iota=2}^{[k/2]} (1+t)^{-2-2\iota } \int r^2 \bar\rho_0^\ga \sigma^{  {1-\iota} } \lt|\pl_t^{k-\iota}\lt(\zeta, r \zeta_r\rt)\rt|^2 dr\notag \\
 \lesssim  &   (1+t)^{-2-2k} \sum_{\iota=0}^{k-1}  \mathcal{E}_\iota(t)\notag .
\end{align}
It then gives from   the Cauchy inequality  and \ef{3-65} that for any $\da>0$,
\be\label{3-67}\begin{split}
  |M_2| \lesssim  &   \da (1+t)^{-1}\int r^2 \bar\rho_0^\ga   \lt| \pl_t^k \lt(\zeta  , r  \zeta_r\rt) \rt|^2  dr + \da^{-1} (1+t) \int r^2 \bar\rho_0^\ga   P^2  dr \\
  \lesssim &  \da (1+t)^{-1}\int r^2 \bar\rho_0^\ga   \lt| \pl_t^k \lt(\zeta  , r  \zeta_r\rt) \rt|^2  dr + \da^{-1} (1+t)^{-1-2k} \sum_{\iota=0}^{k-1}  \mathcal{E}_\iota(t) .
\end{split}\ee
This, together with \ef{t416}, proves \ef{3-63} (by choosing suitably small $\da$)  and \ef{3-64}.

\vskip 0.5cm
{\bf Step 2}. To control the fist term on the right-hand side of \ef{3-62}, we will prove   that
\be\label{3-73}\begin{split}
    &  \frac{d}{dt} \mathfrak{E}_k +  \int    \lt[(1+ t)^{-1} r^2 \bar\rho_0^\ga  \lt|\pl_t^k \lt( \zeta, r\zeta_r \rt)  \rt|^2  +    r^4 \bar\rho_0    \lt(\pl_t^{k} \zeta_t\rt)^2 \rt]   dr \\
  \lesssim &    (1+t)^{-1-2k} \sum_{\iota=0}^{k-1}  (1+ t)^{2\iota} \int    \lt[ r^2 \bar\rho_0^\ga  \lt|\pl_t^\iota \lt( \zeta, r\zeta_r \rt)  \rt|^2  + (1+ t)  r^4 \bar\rho_0    \lt(\pl_t^{\iota} \zeta_t\rt)^2 \rt]   dr,
\end{split}\ee
where
\bee\label{defnek}
  \mathfrak{E}_k : = \int r^4  \bar\rho_0 \lt[ \lt(\pl_t^k \zeta_t \rt)^2  +    \lt(\pl_t^k \zeta \rt)  \pl_t^k \zeta_t    +
   \frac{1}{2}  \lt(\pl_t^k \zeta  \rt)^2   \rt] dr + 2  \overline{\mathfrak{E}}_k.
\eee
We start with   integrating the product of \ef{4-21} and $r^3 \pl_t^k \zeta$ with respect to $r$ to give
\be\label{t412}\begin{split}
   & \frac{d}{dt}\int  r^4  \bar\rho_0 \lt(  \lt(\pl_t^k \zeta\rt) \pl_t^{k } \zeta_t   +
   \frac{1}{2}  \lt(\pl_t^k \zeta\rt)^2    \rt) dr - \int  r^4  \bar\rho_0 \lt(\pl_t^{k } \zeta_t \rt)^2 dr+  M_1 + M_2=0,
\end{split}\ee
where
\bee\label{}\begin{split}
    M_1= &  - \int {\bar\rho_0}^\ga \lt(w_1 \pl_t^k \zeta + w_2 r\pl_t^k \zeta_r \rt)  \lt(r^3 \pl_t^k \zeta \rt)_r dr + \int r^2 \bar\rho_0^\ga  (3w_2-w_1)\lt( r \pl_t^k \zeta_r \rt)\pl_t^k \zeta  dr \\
 &-2 \int r^2 \bar\rho_0^\ga    w_3  \lt(r\zeta_r \rt) \lt(\pl_t^k \zeta\rt)^2 dr
\end{split}\eee
and $M_2$ is defined in \ef{m2}.
 A direct calculation shows that $M_1$ is positive and can be bounded from below as follows
\bee\label{}\begin{split}
 M_1 =& -\int r^2 \bar\rho_0^\ga \lt[\lt(3w_1+2w_3r\zeta_{r}\rt) \lt(\pl_t^k \zeta\rt)^2 + 2w_1 \lt(\pl_t^k \zeta\rt)\lt( r \pl_t^k \zeta_r \rt)
 + w_2 \lt( r \pl_t^k \zeta_r \rt)^2 \rt] dr \\
  \ge & \int r^2 \bar\rho_0^\ga \tilde\eta_r^{1-3\ga}
 \lt\{(9\ga-6) \lt(\pl_t^k \zeta\rt)^2 +  (6\ga-4) \lt(\pl_t^k \zeta\rt)\lt( r \pl_t^k \zeta_r \rt)
 + \ga\lt( r \pl_t^k \zeta_r \rt)^2  \rt.\\
 &\lt. - C\lt(|\zeta|+|r\zeta_{r}|\rt)
   \lt[ \lt(\pl_t^k \zeta\rt)^2+\lt( r \pl_t^k \zeta_r \rt)^2 \rt]\rt\} dr
 \\
  \gtrsim  & \int r^2 \bar\rho_0^\ga \tilde\eta_r^{1-3\ga}
 \lt[ \lt(\pl_t^k \zeta\rt)^2+\lt( r \pl_t^k \zeta_r \rt)^2 \rt]dr
  \gtrsim   (1+t)^{-1} \int r^2 \bar\rho_0^\ga
 \lt[ \lt(\pl_t^k \zeta\rt)^2+\lt( r \pl_t^k \zeta_r \rt)^2 \rt]  dr,
\end{split}\eee
due to  \ef{w123}, the smallness of $\zeta_r$ and $r\zeta_r$ and \ef{decay}.
We then obtain,
by making a summation of $2\times\ef{3-62}$ and \ef{t412}, that
\be\label{t417}\begin{split}
    &  \frac{d}{dt} \mathfrak{E}_k + \int  r^4  \bar\rho_0  \lt(\pl_t^k \zeta_t\rt)^2 dr
     + (1+t)^{-1} \int r^2 \bar\rho_0^\ga
 \lt[ \lt(\pl_t^k \zeta\rt)^2+\lt( r \pl_t^k \zeta_r \rt)^2 \rt]  dr \\
  \lesssim &  \da (1+t)^{-1}\int r^2 \bar\rho_0^\ga   \lt| \pl_t^k \lt(\zeta  , r  \zeta_r\rt) \rt|^2  dr +  (\ea_0+\da) (1+t)^{  -2k-2} \mathcal{E}_k(t) \\
   &+   \da^{-1} (1+t)^{-1-2k} \sum_{\iota=0}^{k-1}  \mathcal{E}_\iota(t) + \lt(\ea_0+\da^{-1}\rt) (1+t)^{  -2k-2} \sum_{\iota=0}^{k-1}  \mathcal{E}_\iota(t) ,
\end{split}\ee
because of \ef{3-67}.
Notice from the Hardy inequality \ef{hardy2} that for $j=0,1,\cdots, l$,
\bee\label{}\begin{split}
&  \int r^4\bar\rho_0 \lt(\pl_t^j \zeta\rt)^2 dr
  =\int_{\mathcal{I}_o} r^4 \sigma^\aa \lt(\pl_t^j \zeta\rt)^2 dr +
  \int_{\mathcal{I}_b} r^4 \sigma^\aa \lt(\pl_t^j \zeta\rt)^2 dr
  \lesssim  \int_{\mathcal{I}_o} r^2 \sigma^{\aa+1} \lt(\pl_t^j \zeta\rt)^2 dr \\
&  \qquad \qquad + \int_{\mathcal{I}_b}   \sigma^\aa \lt(\pl_t^j \zeta\rt)^2 dr
  \lesssim  \int_{\mathcal{I}_o} r^2 \sigma^{\aa+1} \lt(\pl_t^j \zeta\rt)^2 dr
  + \int_{\mathcal{I}_b}   \sigma^{\aa +2}  \lt[\lt(\pl_t^j \zeta\rt)^2 + \lt(\pl_t^j \zeta_r\rt)^2 \rt]dr \\
&  \lesssim \int  r^2 \sigma^{\aa+1} \lt(\pl_t^j \zeta\rt)^2 dr
  + \int_{\mathcal{I}_b}   \sigma^{\aa +1}  \lt[r^2 \lt(\pl_t^j \zeta\rt)^2 + r^4 \lt(  \pl_t^j \zeta_r\rt)^2 \rt]dr
   \lesssim \int r^2 \bar\rho_0^\ga   \lt| \pl_t^j \lt(\zeta  , r  \zeta_r\rt) \rt|^2  dr.
\end{split}\eee
Thus,
\be\label{3-72}\begin{split}
\mathcal{ E}_{j}(t) \lesssim  (1+ t)^{2j} \int    \lt[ r^2 \bar\rho_0^\ga  \lt|\pl_t^j \lt( \zeta, r\zeta_r \rt)  \rt|^2  + (1+ t)  r^4 \bar\rho_0    \lt(\pl_t^{j} \zeta_t\rt)^2 \rt]   dr, \ \   j=0,\cdots,l.\end{split}\ee
This finishes the proof of \ef{3-73}, by using  \ef{t417} and \ef{3-72}, choosing suitably small $\da$ and noting the smallness of $\ea_0$. Moreover, it follows from   \ef{3-63} and \ef{3-64}  that
 \be\label{3-70}\begin{split}
  \mathfrak{E}_k  \ge &C^{-1} \int r^4  \bar\rho_0 \lt| \pl_t^k \lt(\zeta  ,  \zeta_t\rt) \rt|^2 dr + C^{-1} (1+t)^{-1}\int r^2 \bar\rho_0^\ga   \lt| \pl_t^k \lt(\zeta  , r  \zeta_r\rt) \rt|^2  dr \\
  &- C (1+t)^{-2k-1} \sum_{\iota=0}^{k-1}  \mathcal{E}_\iota(t),
\end{split}\ee
\be\label{3-71}\begin{split}
 \mathfrak{E}_k  \lesssim &  \int r^4  \bar\rho_0 \lt| \pl_t^k \lt(\zeta  ,  \zeta_t\rt) \rt|^2 dr +  (1+t)^{-1}\int r^2 \bar\rho_0^\ga   \lt| \pl_t^k \lt(\zeta  , r  \zeta_r\rt) \rt|^2  dr \\
 & +  (1+t)^{-2k-1} \sum_{\iota=0}^{k-1}  \mathcal{E}_\iota(t)
.
\end{split}\ee

\vskip 0.5cm
{\bf Step 3}. In this step, we show the time decay of the norm. We integrate \ef{3-73} and use the induction hypothesis \ef{t411} to show, noting \ef{3-70} and \ef{3-71}, that
 \bee\label{}\begin{split}
 & \int \lt[ r^4  \bar\rho_0 \lt| \pl_t^k \lt(\zeta  ,  \zeta_t\rt) \rt|^2  + (1+t)^{-1} r^2 \bar\rho_0^\ga   \lt| \pl_t^k \lt(\zeta  , r  \zeta_r\rt) \rt|^2  \rt](r,t) dr \\
 & +\int_0^t \int    \lt[(1+ s)^{-1} r^2 \bar\rho_0^\ga  \lt|\pl_s^k \lt( \zeta, r\zeta_r \rt)  \rt|^2  +    r^4 \bar\rho_0    \lt(\pl_s^{k} \zeta_s\rt)^2 \rt]   drds \\
 \lesssim  & \sum_{\iota=0}^k \mathcal{E}_\iota (0) +
 \sum_{\iota=0}^{k-1}  \int_0^t    (1+ s)^{ 2\iota -1-2k} \int    \lt[ r^2 \bar\rho_0^\ga  \lt|\pl_s^\iota \lt( \zeta, r\zeta_r \rt)  \rt|^2  + (1+ s)  r^4 \bar\rho_0    \lt(\pl_s^{\iota} \zeta_s\rt)^2 \rt]   dr \\
 \lesssim & \sum_{\iota=0}^k \mathcal{E}_\iota (0).
\end{split}\eee
Multiply  \ef{3-73} by $(1+t)^{p} $ and integrate the product with respect to the temporal variable from $p=1$ to $p=2k$  step by step to get
\be\label{3-74}\begin{split}
 & (1+t)^{2k}\int \lt[ r^4  \bar\rho_0 \lt| \pl_t^k \lt(\zeta  ,  \zeta_t\rt) \rt|^2  + (1+t)^{-1} r^2 \bar\rho_0^\ga   \lt| \pl_t^k \lt(\zeta  , r  \zeta_r\rt) \rt|^2  \rt](r,t) dr \\
 & +\int_0^t (1+s)^{2k} \int    \lt[(1+ s)^{-1} r^2 \bar\rho_0^\ga  \lt|\pl_s^k \lt( \zeta, r\zeta_r \rt)  \rt|^2  +    r^4 \bar\rho_0    \lt(\pl_s^{k} \zeta_s\rt)^2 \rt]   drds \\
 \lesssim  & \sum_{\iota=0}^k \mathcal{E}_\iota (0) +
 \sum_{\iota=0}^{k-1}  \int_0^t    (1+ s)^{ 2\iota -1 } \int    \lt[ r^2 \bar\rho_0^\ga  \lt|\pl_s^\iota \lt( \zeta, r\zeta_r \rt)  \rt|^2  + (1+ s)  r^4 \bar\rho_0    \lt(\pl_s^{\iota} \zeta_s\rt)^2 \rt]   dr \\
 \lesssim & \sum_{\iota=0}^k \mathcal{E}_\iota (0).
\end{split}\ee
With this estimate at hand, we finally integrate $(1+t)^{2k+1}\ef{3-62}$ with respect to the temporal variable and use \ef{3-72}, \ef{t411} and \ef{3-74}  to show
\be\label{3-75}\begin{split}
 & (1+t)^{2k}\int \lt[ (1+t) r^4  \bar\rho_0 \lt| \pl_t^k    \zeta_t  \rt|^2  +   r^2 \bar\rho_0^\ga   \lt| \pl_t^k \lt(\zeta  , r  \zeta_r\rt) \rt|^2  \rt](r,t) dr \\
 & +\int_0^t (1+s)^{2k+1} \int     r^4 \bar\rho_0    \lt(\pl_s^{k} \zeta_s\rt)^2     drds \\
 \lesssim  & \sum_{\iota=0}^k \mathcal{E}_\iota (0) +
 \sum_{\iota=0}^{k }  \int_0^t    (1+ s)^{ 2\iota -1 } \int    \lt[ r^2 \bar\rho_0^\ga  \lt|\pl_s^\iota \lt( \zeta, r\zeta_r \rt)  \rt|^2  + (1+ s)  r^4 \bar\rho_0    \lt(\pl_s^{\iota} \zeta_s\rt)^2 \rt]   dr \\
 \lesssim & \sum_{\iota=0}^k \mathcal{E}_\iota (0).
\end{split}\ee
It finally follows from \ef{3-74} and \ef{3-75} that
\bee\label{}\begin{split}
      \mathcal{E}_k(t)   +
 \int_0^t  \int  (1+s)^{2k-1}\lt[  r^2 \bar\rho_0^\ga \lt| \pl_s^k\lt( \zeta, r \zeta_r\rt)\rt|^2 +(1+s)^{2}  r^4 \bar\rho_0 \lt(\pl_s^{k } \zeta_s \rt)^2  \rt] drds
  \lesssim     \sum_{\iota=0}^k \mathcal{E}_\iota (0) .
\end{split}\eee
This completes the proof of Lemma \ref{lem332}. $\Box$

\subsection{Verification of the a priori assumption}
In this subsection, we prove the following lemma.
\begin{lem}\label{lem34} Suppose that $\mathcal{E}(t)$ is finite, then it holds that
\begin{align}\label{lem34est}
  &  \sum_{j=0}^2 (1+  t)^{2j} \lt \|\pl_t^j \zeta (\cdot,t) \rt \|_{L^\iy}^2   +  \sum_{j=0}^1 (1+  t)^{2j} \lt \|\pl_t^j \zeta_r (\cdot,t) \rt\|_{L^\iy}^2
     +
  \sum_{
  i+j\le l-2,\   2i+j \ge 3} (1+  t)^{2j}   \notag\\
    & \times   \lt \|  \sigma^{\f{ 2i+j-3 }{2}}\pl_t^j \pl_r^i \zeta(\cdot,t) \rt\|_{L^\iy}^2 +
  \sum_{
  i+j=l-1} (1+  t)^{2j}   \lt \| r \sigma^{\f{ 2i+j-3 }{2}}\pl_t^j \pl_r^i \zeta  (\cdot,t) \rt\|_{L^\iy}^2
   \notag \\
 & + \sum_{
  i+j=l} (1+  t)^{2j}   \lt \| r^2 \sigma^{\f{ 2i+j-3 }{2}}\pl_t^j \pl_r^i \zeta (\cdot,t) \rt\|_{L^\iy}^2 \lesssim \mathcal{E}(t).
\end{align}
\end{lem}
Once this lemma is proved, the a priori assumption \ef{apb} is then verified and the proof of Theorem \ref{thm1} is finished, since it follows from the elliptic estimate \ef{prop1est}  and the nonlinear weighted energy estimate \ef{prop2est} that
$$\mathcal{E}(t)\lesssim \mathcal{E}(0), \ \  t\in [0,T].$$

\vskip 0.5cm
\noindent {\bf Proof}. The proof consists of two steps. In Step 1, we derive the $L^\iy$-bounds away from the boundary, that is,
  \begin{align}\label{3-83}
 & \sum_{i+j \le l-2} \lt\|\pl_t^j \pl_r^i \zeta\rt\|_{L^\iy(\mathcal{I}_o)}^2
  +\sum_{i+j =l-1} \lt\|r\pl_t^j \pl_r^i \zeta\rt\|_{L^\iy(\mathcal{I}_o)}^2
  +\sum_{i+j =l-2} \lt\|r^2\pl_t^j \pl_r^i \zeta\rt\|_{L^\iy(\mathcal{I}_o)}^2
\notag\\
& \lesssim   (1+t)^{-2j} \mathcal{E}(t).
\end{align}
Away from the origin, we show in Step 2 the following $L^\iy$-estimates:
\be\label{3-85}\begin{split}
\sum_{j=0}^3 (1+  t)^{2j} \lt\|\pl_t^j \zeta\rt\|_{L^\iy(\mathcal{I}_b)}^2   +  \sum_{j=0}^1 (1+  t)^{2j} \lt\|\pl_t^j \zeta_r\rt\|_{L^\iy(\mathcal{I}_b)}^2 \lesssim \mathcal{E}(t),
\end{split}\ee
\be\label{fak3-77}\begin{split}
\lt\| \sigma^{\f{2i+j-3}{2}}\pl_t^j \pl_r^i \zeta \rt\|_{L^\iy(\mathcal{I}_b)}^2  \lesssim     (1+t)^{-2j}\mathcal{E}(t)  \  \ {\rm when} \  \  2i+j \ge 4.
 \end{split}\ee
We obtain \ef{lem34est} by using \ef{3-83}-\ef{fak3-77} and noting the facts $l\ge 4$ and $\mathcal{I}=\mathcal{I}_o \cup \mathcal{I}_b$. It suffices to  show \ef{3-83}-\ef{fak3-77}.

To this end, we first notice some facts. It follows from \ef{3-27} that $\mathcal{E}_{j,0}\lesssim \mathcal{E}_j$ for $j=0,\cdots, l$, which implies
 \be\label{newnorm}\begin{split}
\sum_{j=0}^l \lt(\mathcal{ E}_{j}(t) + \sum_{i=0}^{l-j} \mathcal{ E}_{j, i}(t)  \rt)\lesssim   \mathcal{E}(t).
\end{split}\ee
The following embedding (cf. \cite{adams}): $H^{1/2+\da}(\mathcal{I})\hookrightarrow L^\iy(\mathcal{I}) $
with the estimate
\be\label{half} \|F\|_{L^\iy(\mathcal{I})} \le C(\delta) \|F\|_{H^{1/2+\da}(\mathcal{I})},\ee
for $\da>0$ will be used in the rest of the proof.

{\bf Step 1} (away from the boundary). It follows from \ef{newnorm} that for $j=0,1,\cdots, l$,
\begin{align}\label{3-78}
(1+t)^{2j} \lt[\sum_{i=0}^{l-j}\int_{\mathcal{I}_o}   r^2   \lt(\pl_t^j \pl_r^i \zeta\rt)^2 dr + \int_{\mathcal{I}_o}   r^4   \lt(\pl_t^j \pl_r^{l-j+1} \zeta\rt)^2 dr \rt]    \lesssim \sum_{i=0}^{l-j} \mathcal{ E}_{j, i}(t) \le  \mathcal{E}(t),
\end{align}
which implies, using \ef{wsv}, that for  $j=0,1,\cdots, l-1$,
\begin{align}\label{3-79}
\lt\|\pl_t^j \zeta\rt\|_{H^{l-j-1}(\mathcal{I}_o)}^2
\lesssim & \lt\|\pl_t^j \zeta\rt\|_{H^{2, l-j}(\mathcal{I}_o)}^2 =  \sum_{i=0}^{l-j}\int_{\mathcal{I}_o}   {dist^2(r, \pl\mathcal{I}_o )}    \lt(\pl_t^j \pl_r^i \zeta\rt)^2 dr
\notag\\
\le  & \sum_{i=0}^{l-j}\int_{\mathcal{I}_o}   r^2    \lt(\pl_t^j \pl_r^i \zeta\rt)^2 dr
   \le (1+t)^{-2j} \mathcal{E}(t)  .
\end{align}
In view of \ef{half} and \ef{3-79}, we  see that  for $j=0,1,\cdots, l-2$,
\begin{align}\label{3-80}
\sum_{i=0}^{l-j-2}\lt\|\pl_t^j \pl_r^i \zeta\rt\|_{L^\iy(\mathcal{I}_o)}^2
\lesssim \sum_{i=0}^{l-j-2}\lt\|\pl_t^j \pl_r^i \zeta\rt\|_{H^1(\mathcal{I}_o)}^2
\lesssim & \lt\|\pl_t^j \zeta\rt\|_{H^{ l-j-1}(\mathcal{I}_o)}^2 \lesssim (1+t)^{-2j} \mathcal{E}(t).
\end{align}
It gives from \ef{half}, \ef{3-78} and \ef{3-79} that for
$j=0,1,\cdots, l-1$,
 \begin{align}\label{3-81}
 \lt\|r\pl_t^j \pl_r^{l-j-1} \zeta\rt\|_{L^\iy(\mathcal{I}_o)}^2 \lesssim  &
 \lt\|r\pl_t^j \pl_r^{l-j-1} \zeta\rt\|_{H^1(\mathcal{I}_o)}^2 \lesssim \lt\| \pl_t^j \pl_r^{l-j-1} \zeta\rt\|_{L^2(\mathcal{I}_o)}^2 + \lt\|r\pl_t^j \pl_r^{l-j } \zeta\rt\|_{L^2(\mathcal{I}_o)}^2 \notag\\
 \lesssim & (1+t)^{-2j} \mathcal{E}(t)
\end{align}
and for
$j=0,1,\cdots, l$,
  \begin{align}\label{3-82}
 \lt\|r^2 \pl_t^j \pl_r^{l-j} \zeta\rt\|_{L^\iy(\mathcal{I}_o)}^2 \lesssim  &
 \lt\|r^2 \pl_t^j \pl_r^{l-j} \zeta\rt\|_{H^1(\mathcal{I}_o)}^2 \lesssim \lt\|r \pl_t^j \pl_r^{l-j } \zeta\rt\|_{L^2(\mathcal{I}_o)}^2 + \lt\|r^2 \pl_t^j \pl_r^{l-j +1 } \zeta\rt\|_{L^2(\mathcal{I}_o)}^2 \notag\\
 \lesssim & (1+t)^{-2j} \mathcal{E}(t).
\end{align}
So that we can derive  \ef{3-83}  from   \ef{3-80}-\ef{3-82}.

{\bf Step 2} (away from the origin).  We set
\begin{align}\label{3-84}
 d_b(r): = dist (r, \pl \mathcal{I}_b) \le \sqrt{A/B}-r \lesssim \sa(r), \ \  r\in \mathcal{I}_b.\end{align}
It follows from \ef{wsv} and \ef{3-84} that for $j\le 5+[\aa]-\aa$,
\begin{align}\label{}
&\lt\|\pl_t^j \zeta \rt\|_{H^{\frac{5-j+[\aa]-\aa}{2}}(\mathcal{I}_b)}^2
=  \lt\|\pl_t^j \zeta \rt\|_{H^{l-j+1-\frac{l-j+1+\aa}{2}}(\mathcal{I}_b)}^2
\lesssim\lt\|\pl_t^j \zeta \rt\|_{H^{ {l-j+1+\aa} , l-j+1}(\mathcal{I}_b)}^2\notag \\
= &\sum_{k=0}^{l-j+1} \int_{\mathcal{I}_b}    d_b^{\aa+1+l-j}(r) |\pl_r^k \pl_t^j \zeta|^2dr
\lesssim \sum_{k=0}^{l-j+1} \int_{\mathcal{I}_b}   \sigma^{\aa+1+l-j}|\pl_r^k \pl_t^j \zeta|^2dr
\notag\\
\lesssim & \sum_{k=0}^{l-j+1} \int_{\mathcal{I}_b}   \sigma^{\aa+k}|\pl_r^k \pl_t^j \zeta|^2dr \lesssim \sum_{k=0}^{l-j+1} \int_{\mathcal{I}_b} r^4  \sigma^{\aa+k}|\pl_r^k \pl_t^j \zeta|^2dr \notag\\
\le &    (1+t)^{-2j}\lt(\mathcal{ E}_{j}(t) + \sum_{k=1}^{l-j } \mathcal{ E}_{j,k}(t)\rt)\le (1+t)^{-2j}  \mathcal{E}(t).\notag
\end{align}
This, together with \ef{half}, gives \ef{3-85}.

To prove  \ef{fak3-77},   we denote
$$\psi:=  \sigma^{\f{2i+j-3}{2}}\pl_t^j \pl_r^i \zeta.$$
In what follows, we assume $2i+j\ge 4$ and $i+j\le l$ and show that
\be\label{hahahaha}
\lt\|\psi\rt\|_{L^\iy(\mathcal{I}_b)}^2
\lesssim     (1+t)^{-2j}\mathcal{E}(t) .
\ee
The estimate \ef{hahahaha} will be proved  by separating the cases when $\alpha$ is or is not an integer.

{\bf Case 1} ($\aa\neq [\aa]$). When $\aa$ is not an integer,  we choose  $  \sigma^{2(l-i-j)+ \aa-[\aa]}$ as the spatial weight.
A simple calculation yields
\bee\label{}\begin{split}
 \lt|\pl_r   \psi \rt| \lesssim\lt|  \sigma^{\f{2i+j-3}{2}}\pl_t^j \pl_r^{i+1} \zeta\rt|
 + \lt|  \sigma^{\f{2i+j-3}{2}-1}\pl_t^j \pl_r^i \zeta\rt|,
 \end{split}\eee
\bee\label{}\begin{split}
 \lt|\pl_r^2 \psi \rt| \lesssim\lt|  \sigma^{\f{2i+j-3}{2}}\pl_t^j \pl_r^{i+2}\zeta\rt|
 + \lt|  \sigma^{\f{2i+j-3}{2}-1}\pl_t^j \pl_r^{i+1} \zeta \rt| + \lt|  \sigma^{\f{2i+j-3}{2}-2}\pl_t^j \pl_r^i \zeta\rt|,
\end{split}\eee
$$\cdots\cdots$$
\be\label{l313}\begin{split}
 &\lt|\pl_r^k \psi \rt| \lesssim  \sum_{p=0}^k \lt|  \sigma^{\f{2i+j-3}{2}-p}\pl_t^j \pl_r^{i+k-p} \zeta\rt|   \ \ {\rm for} \ \   k=1, 2, \cdots, l+1-j-i.
\end{split}\ee
It  follows from \ef{l313} that for $1\le k \le l+1-i-j$,
\bee\label{}\begin{split}
& \int _{\mathcal{I}_b}  \sigma^{2(l-i-j)+ \aa-[\aa]  } \lt|\pl_r^k \psi \rt|^2 dr \lesssim   \sum_{p=0}^k \int_{\mathcal{I}_b}      \sigma^{\aa+l-j+1-2p}\lt|\pl_t^j  \pl_r^{i+k-p} \zeta\rt|^2 dr \\
 \lesssim &   \int_{\mathcal{I}_b}   \sigma^{l-i-j+1-k}  \sum_{p=0}^1     \sigma^{\aa+i+k-2p} \lt|\pl_t^j  \pl_r^{i+k-p} \zeta\rt|^2 dr +  \sum_{p=2}^k \int_{\mathcal{I}_b}  \sigma^{\aa+l-j+1-2p}\lt|\pl_t^j  \pl_r^{i+k-p} \zeta\rt|^2 dr
 \\ \lesssim &   \sum_{p=0}^1 \int_{\mathcal{I}_b}    r^4    \sigma^{\aa+i+k-2p} \lt|\pl_t^j  \pl_r^{i+k-p} \zeta\rt|^2 dr + \sum_{p=2}^k  \int_{\mathcal{I}_b}     \sigma^{\aa+l-j+1-2p}\lt|\pl_t^j  \pl_r^{i+k-p} \zeta\rt|^2 dr
 \\
 \lesssim & (1+t)^{-2j} \mathcal{E}_{j,i+k-1}+  \sum_{p=2}^k \int_{\mathcal{I}_b}  \sigma^{\aa+l-j+1-2p}\lt|\pl_t^j  \pl_r^{i+k-p} \zeta\rt|^2 dr.
\end{split}\eee
To bound the 2nd term on the right-hand side of the inequality above, notice that
\be\label{3.14}\begin{split}
&\aa+l-j+1-2p \\
= &2( l+1-i-j- k)   +2(k-p) + (\aa-[\aa])+(2i+j-4) -1  > -1
\end{split}\ee
for $p\in [2, k]$, due to $\aa> [\aa]$ and $2i+j\ge 4$. We then have, with the aid of the Hardy inequality \ef{hardy2}, that for $p\in [2, k]$,
\bee\label{}\begin{split}
  & \int_{\mathcal{I}_b}  \sigma^{\aa+l-j+1-2p}\lt|\pl_t^j  \pl_r^{i+k-p} \zeta\rt|^2 dr
  \lesssim  \int_{\mathcal{I}_b}  \sigma^{\aa+l-j+1-2p+2} \sum_{\iota=0}^1 \lt|\pl_t^j  \pl_r^{i+k-p+\iota} \zeta\rt|^2 dr  \lesssim \cdots \\
  & \lesssim  \int_{\mathcal{I}_b}  \sigma^{\aa+l-j+1} \sum_{\iota=0}^p \lt|\pl_t^j  \pl_r^{i+k-p+\iota} \zeta\rt|^2 dr
  = \sum_{\iota=0}^p  \int_{\mathcal{I}_b}   \sigma^{(l+1-i-j-k)+(p-\iota)}     \sigma^{\aa+i+k-p+\iota }
 \lt|\pl_t^j  \pl_r^{i+k-p+\iota} \zeta\rt|^2 dr \\
& \lesssim \sum_{\iota=0}^p  \int_{\mathcal{I}_b}   r^4   \sigma^{\aa+i+k-p+\iota }
 \lt|\pl_t^j  \pl_r^{i+k-p+\iota} \zeta\rt|^2 dr
 \le   \sum_{\iota=i+k-p}^{i+k-1} (1+t)^{-2j} \mathcal{E}_{j,\iota}.
\end{split}\eee
That yields,  for $1\le k \le l+1-i-j$,
\bee\label{}\begin{split}
  \int _{\mathcal{I}_b}  \sigma^{2(l-i-j)+ \aa-[\aa]  } \lt|\pl_r^k \psi \rt|^2 dr
 \lesssim & (1+t)^{-2j} \mathcal{E}_{j,i+k-1}+ \sum_{p=2}^k  \sum_{\iota=i+k-p}^{i+k-1} (1+t)^{-2j} \mathcal{E}_{j,\iota} \\
 \lesssim &(1+t)^{-2j} \sum_{\iota=i}^{i+k-1} \mathcal{E}_{j,\iota}.
\end{split}\eee
Therefore,  it follows from  \ef{3-84}  and \ef{newnorm} that
\bee\label{}\begin{split}
&\lt\|\psi\rt\|_{H^{2(l-i-j)+\aa-[\aa] , \ l+1-i-j }(\mathcal{I}_b)}^2 = \sum_{k=0}^{l+1-i-j} \int_{\mathcal{I}_b}    d_b^{2(l-i-j)+ \aa-[\aa]  } \lt|\pl_r^k \psi \rt|^2dr
\\
\lesssim & \sum_{k=0}^{l+1-i-j} \int_{\mathcal{I}_b}   \sigma^{2(l-i-j)+ \aa-[\aa]  } \lt|\pl_r^k \psi \rt|^2dr  \lesssim    \int_{\mathcal{I}_b}      \sigma^{\aa+l-j+1 }\lt|\pl_t^j  \pl_r^{i  } \zeta\rt|^2 dr  +  (1+t)^{-2j} \sum_{\iota=i}^{l-j} \mathcal{E}_{j,\iota} \\
\lesssim &  \int_{\mathcal{I}_b}   r^4   \sigma^{\aa+i+1 }\lt|\pl_t^j  \pl_r^{i  } \zeta\rt|^2 dr  +  (1+t)^{-2j} \sum_{\iota=i}^{l-j} \mathcal{E}_{j,\iota}
 \lesssim    (1+t)^{-2j} \sum_{\iota=i}^{l-j} \mathcal{E}_{j,\iota}
 \le  (1+t)^{-2j}\mathcal{E}(t).
\end{split}\eee
When $\aa$ is not an integer, $\aa-[\aa]\in(0,1)$. So, it follows from  \ef{half} and \ef{wsv}  that
\be\label{3-89}\begin{split}
\lt\|\psi\rt\|_{L^\iy(\mathcal{I}_b)}^2 \lesssim \lt\|\psi\rt\|_{H^{1-\frac{\aa-[\aa]}{2}}(\mathcal{I}_b)}^2
\lesssim \lt\|\psi\rt\|_{H^{2(l-i-j)+\aa-[\aa] , \ l+1-i-j }(\mathcal{I}_b)}^2
\lesssim     (1+t)^{-2j}\mathcal{E}(t).
 \end{split}\ee
{\bf Case 2} ($\aa=[\aa]$). In this case $\aa$ is an integer,  we choose  $  \sigma^{2(l-i-j)+ 1/2 }$ as the spatial weight. As shown in Case 1, we
have   for $1\le k \le l+1-i-j$,
\bee\label{}\begin{split}
  \int_{\mathcal{I}_b}   \sigma^{2(l-i-j)+ 1/2} \lt|\pl_r^k \psi \rt|^2dr
 \lesssim   (1+t)^{-2j} \mathcal{E}_{j,i+k-1}+ \sum_{p=2}^k  \int_{\mathcal{I}_b}     \sigma^{\aa+l-j+1-2p+ 1/2}\lt|\pl_t^j  \pl_r^{i+k-p} \zeta\rt|^2 dr.
\end{split}\eee
Note that for $1\le k \le l+1-i-j$ and  $2\le p\le k$,
\bee\label{}\begin{split}
 \aa+l-j+1-2p + \frac{1}{2}
=  2( l+1-i-j- k)   +2(k-p)  +(2i+j-4) - \frac{1}{2} \ge -\frac{1}{2}.
\end{split}\eee
We can then use the Hardy inequality \ef{hardy2} to obtain
\bee\label{}\begin{split}
  \int_{\mathcal{I}_b}   \sigma^{2(l-i-j)+ 1/2 } \lt|\pl_r^k \psi \rt|^2 dr
 \lesssim  (1+t)^{-2j} \sum_{\iota=i}^{i+k-1} \mathcal{E}_{j,\iota}, \ \ k=1, 2, \cdots, l-j+1-i,
\end{split}\eee
which, together with  \ef{3-84}  and \ef{newnorm}, implies that
\bee\label{}\begin{split}
 \lt\|\psi\rt\|_{H^{2(l-i-j)+1/2 , \ l+1-i-j }(\mathcal{I}_b)}^2
 \lesssim    (1+t)^{-2j} \sum_{\iota=i}^{l-j} \mathcal{E}_{j,\iota} \le  (1+t)^{-2j}\mathcal{E}(t).
\end{split}\eee
Therefore, it follows from  \ef{half} and \ef{wsv}  that
\be\label{3-90}\begin{split}
\lt\|\psi\rt\|_{L^\iy(\mathcal{I}_b)}^2 \lesssim \lt\|\psi\rt\|_{H^{3/4}(\mathcal{I}_b)}^2
\lesssim \lt\|\psi\rt\|_{H^{2(l-i-j)+1/2 , \ l+1-i-j }(\mathcal{I}_b)}^2
\lesssim     (1+t)^{-2j}\mathcal{E}(t).
 \end{split}\ee
In view of \ef{3-89} and \ef{3-90}, we obtain \ef{hahahaha} or equivalently \ef{fak3-77}.
This completes the proof of   Lemma \ref{lem34}. $\Box$

\section{Proof of Theorem \ref{thm2}}
In this section, we prove Theorem \ref{thm2}. First, it follows from \ef{o2.1}, \ef{newlagrangiandensity}, \ef{1.6}, \ef{bareta}, \ef{haz2.2} and \ef{haz2.7} that
for $(r,t)\in \mathcal{I}\times[0,\iy)$,
$$\rho(\eta(r,t),t)-\bar\rho(\bar\eta(r,t),t)=\frac{r^2 \bar\rho_0(r)}{\eta^2(r,t)\eta_r(r,t)}-\frac{r^2 \bar\rho_0(r)}{\bar\eta^2(r,t)\bar\eta_r(r,t)}$$
and
$$u(\eta(r,t),t)-\bar u(\bar\eta(r,t),t)= \eta_t(r,t)- \bar\eta_t(r,t).$$
Then, we have, using \ef{hh2-12}, \ef{teeta}, \ef{bareta}, \ef{basic}, \ef{decayforh} and \ef{thm1est2} that
$$\lt|\rho(\eta(r,t),t)-\bar\rho(\bar\eta(r,t),t)\rt|\lesssim (A-Br^2)^{\frac{1}{\ga-1}}(1+t)^{-\frac{4}{3\ga-1}}\lt[\sqrt{\mathcal{E}(0)} + (1+t)^{-\frac{3\ga-2}{3\ga-1}}\ln(1+t)\rt]$$
and
$$\lt|u(\eta(r,t),t)-\bar u(\bar\eta(r,t),t)\rt|\lesssim r (1+t)^{-1}\lt[\sqrt{\mathcal{E}(0)} + (1+t)^{-\frac{3\ga-2}{3\ga-1}}\ln(1+t)\rt] .$$
This gives the proof of \ef{thm2est1} and \ef{thm2est2}.

For the boundary behavior, it follows from \ef{vbs}, \ef{hh2-12}, \ef{teeta} and  \ef{bareta} that
\bee\label{}\begin{split}
 R(t)=&\eta\lt(\sqrt{A/B},t\rt)= \lt(\tilde\eta+r \zeta\rt)\lt(\sqrt{A/B},t\rt)
 =\lt(\bar\eta+rh+r \zeta\rt)\lt(\sqrt{A/B},t\rt) \\
 =&
 \lt[ r\lt(\bar\eta_r+h+\zeta\rt)\rt]\lt(\sqrt{A/B},t\rt)= \sqrt{A/B} \lt[(1+t)^{1/(3\ga-1)} + h(t)+ \zeta(t)\rt]
  \end{split}\eee
which, together with \ef{thm1est2} and \ef{decayforh} that
$$R(t)\ge \sqrt{A/B} \lt[(1+t)^{\f{1}{3\ga-1}} - C \sqrt{\mathcal{E}(0)}\rt]$$
and
$$R(t)\le \sqrt{A/B} \lt[(1+t)^{\f{1}{3\ga-1}} + C(1+t)^{-\frac{3\ga-2}{3\ga-1}}\ln(1+t) + C \sqrt{\mathcal{E}(0)}\rt].$$
Thus, \ef{thm2est3} follows from the smallness of $\mathcal{E}(0)$. Notice that for  $k=1,2,3$,
\bee\label{}\begin{split}
 \frac{d^k R(t)}{dt^k}= \pl_t^{k}\tilde\eta \lt(\sqrt{A/B},t\rt)  +\lt( r \pl_t^k \zeta\rt)\lt(\sqrt{A/B},t\rt).
  \end{split}\eee
Therefore, \ef{thm2est4} follows from \ef{decay} and \ef{thm1est2}.

We are to verify the physical vacuum condition, \ef{thm2est5}. It follows from \ef{o2.1}, \ef{newlagrangiandensity}, \ef{hh2-12} that
\bee\label{}\begin{split}
 &\lt(\rho^{\ga-1}\rt)_\eta(\eta,t) = \frac{\lt(f^{\ga-1}\rt)_r(r,t)}{\eta_r(r,t)}
 =\frac{1}{\eta_r} \lt[\bar\rho_0^{\ga-1}\lt(\frac{r^2}{\eta^2 \eta_r}\rt)^{\ga-1}\rt]_r
 \\
 =&\frac{1}{\eta_r} \lt\{  \bar\rho_0^{\ga-1}\lt[\bar\rho_0^{\ga-1}\lt(\frac{r^2}{\eta^2 \eta_r}\rt)^{\ga-1}\rt]_r-2Br\lt(\frac{r^2}{\eta^2 \eta_r}\rt)^{\ga-1}  \rt\}\\
 =&  (1-\ga)\bar\rho_0^{\ga-1}\lt[2 \lt(\frac{\eta}{r} \rt)^{1-2\ga }\eta_r^{ -\ga} \zeta_r
  + \lt(\frac{\eta}{r} \rt)^{2-2\ga}\eta_r^{-\ga-1}\lt(2\zeta_r + r\zeta_{rr}\rt)\rt]  -2Br\lt(\frac{\eta}{r} \rt)^{2-2\ga}\eta_r^{ -\ga},
  \end{split}\eee
which implies, with the aid of \ef{ht3.3} and \ef{thm1est2}, that
  \be\label{4-.1}\begin{split}
  \lt|\lt(\rho^{\ga-1}\rt)_\eta(\eta,t)\rt| \lesssim   &
 (1+t)^{-1}\sqrt{\mathcal{E}(0)} + r (1+t)^{-1 + \frac{1}{3\ga-1}}
  \end{split}\ee
  and
    \be\label{4-.2}\begin{split}
 \lt| \lt(\rho^{\ga-1}\rt)_\eta(\eta,t) \rt| \ge & 2Br\lt(\frac{\eta}{r} \rt)^{2-2\ga}\eta_r^{ -\ga} -
 C (1+t)^{-1}\sqrt{\mathcal{E}(0)}  \\
 \ge  & C^{-1} r (1+t)^{-1 + \frac{1}{3\ga-1}} -
 C (1+t)^{-1}\sqrt{\mathcal{E}(0)} .
  \end{split}\ee
In view of  \ef{ht3.3}, we see that
$$\eta(r,t) \sim (1+t)^{\frac{1}{3\ga-1}} r,$$
which, together with \ef{thm2est3}, \ef{4-.1} and \ef{4-.2}, gives for $R(t)/2 \le \eta \le R(t)$,
  \bee\label{}\begin{split}
 C^{-1}   (1+t)^{-\frac{3\ga-2}{3\ga-1}} -
 C (1+t)^{-1}\sqrt{\mathcal{E}(0)} \le \lt|\lt(\rho^{\ga-1}\rt)_\eta(\eta,t)\rt| \lesssim
 (1+t)^{-1}\sqrt{\mathcal{E}(0)} +   (1+t)^{-\frac{3\ga-2}{3\ga-1}}.
  \end{split}\eee
Thus,  \ef{thm2est5} follows from the smallness of $\mathcal{E}(0)$.
This finishes the proof of Theorem \ref{thm2}. $\Box$

\section{Appendix}
{\bf Proof of $\ef{decay}_1$}. We may write \ef{pomt} as the following system
\be\label{s5.1}\lt\{\begin{split}
 &h_t= z,\\
 &z_t= -z - \lt[ \bar\eta_r^{2-3\ga} -  (\bar \eta_r+h)^{2-3\ga} \rt] /( {3\ga - 1} )- \bar\eta_{rtt}  \\
 &(h,z)(t=0)=(0,0).
\end{split}\rt.\ee
Recalling that $\bar\eta_r (t)=(1+t)^{{1}/({3\ga-1})}$, thus  $\bar\eta_{rtt}<0$. A simple phase plane analysis shows that there exist $0<t_0<t_1<t_2$ such that, starting from $(h, z)=(0, 0)$ at $t=0$, $h$ and $z$ increases in the interval $[0, t_0]$ and $z$ reaches its positive maxima at $t_0$; in the interval $[t_0, t_1]$, $h$ keeps increasing and reaches its maxima at $t_1$,  $z$ decreases from its positive maxima to 0; in the interval $[t_1, t_2]$,
both $h$ and $z$ decrease, and $z$ reaches its negative minima at $t_2$; in the interval $[t_2, \infty)$, $h$ decreases and $z$ increases, and $(h, z)\to (0, 0)$ as $t\to \infty$. This can be summarized as follows:
$$z(t) \uparrow_0 , \  \ h(t) \uparrow_0  ,  \ \ t\in [0,t_0]$$
$$ z(t) \downarrow_0,  \ \  h(t)  \uparrow,  \ \  t\in[t_0,t_1] $$
$$ z(t) \downarrow^0,  \ \  h(t)  \downarrow,  \ \  t\in[t_1,t_2] $$
$$ z(t)  \uparrow^0,  \ \  h(t)  \downarrow_0,  \ \  t\in[t_2,\iy). $$
It gives from the above analysis that
 there exists a finite constant $C=C(\ga, M)$ such that
\be\label{s5.2}
0\le h(t) \le C  \ \ {\rm for} \  \  t\ge 0.
\ee
In view of \ef{bareta} and \ef{teeta}, we then see that
$$
\lt(1 +   t \rt)^{ {1}/({3\ga-1})}  \le \te\eta_r(t) \le K \lt(1 +   t \rt)^{{1}/({3\ga-1})}.
$$
On the other hand, equation \ef{pomt} can be rewritten as
\be\label{s5.3}\begin{split}
&  \te\eta_{rtt} + \te\eta_{rt} -  \te\eta_{r}^{2-3\ga} /(3\ga-1) =0, \ \ t >0,  \\
&  \te\eta_r(t=0)= 1, \ \  \te\eta_{rt}(t=0)= 1/(3\ga-1) .
 \end{split}\ee
Then, we have by solving \ef{s5.3} that
\be\label{s5.4}\begin{split}
\te\eta_{rt}(t)= \frac{1}{3\ga - 1}   e^{-t} +  \frac{1}{3\ga - 1} \int_0^t e^{-(t-s)} \te\eta_{r}^{2-3\ga}(s) ds  \ge 0. \ \  \Box
\end{split}\ee

\vskip 0.5cm
\noindent{\bf Proof of $\ef{decay}_2$}.  We use the mathematical induction to prove $\ef{decay}_2$. First, it follows from \ef{s5.4} that
\be\label{s5.5}\begin{split}
({3\ga-1})  \te\eta_{rt}(t)  = &e^{-t}   + \int_0^{t/2} e^{-(t-s)}\te\eta_{r}^{2-3\ga}(s) ds +   \int_{t/2}^t e^{-(t-s)}\te\eta_{r}^{2-3\ga}(s)ds \\
\le & e^{-t} +  e^{-t/2}\int_0^{t/2}(1+  s)^{\f{2-3\ga}{3\ga-1}}ds
+\lt(1+ \frac{t}{2} \rt)^{\f{2-3\ga}{3\ga-1}} \int_{t/2}^t e^{-(t-s)}  ds \\
\le & e^{-t} + C e^{-t/2} \lt(1+  \frac{t}{2}\rt)^{\f{1}{3\ga-1}}
+\lt(1+ \frac{t}{2} \rt)^{\f{2-3\ga}{3\ga-1}} \\
\le  &  C \lt(1+ {t} \rt)^{\f{2-3\ga}{3\ga-1}} ,  \ \  t\ge 0,
\end{split}\ee
for some constant $C$ independent of $t$. This proves $\ef{decay}_2$ when $k=1$.   Suppose that $\ef{decay}_2$
holds for all $k=1,2,\cdots,m-1$, that is,
\be\label{s5.6}\begin{split}
 \lt|\f{d^k\tilde \eta_{r}(t)}{dt^k}\rt| \le C(m)\lt(1 +   t \rt)^{\frac{1}{3\ga-1}-k},   \ \ k=1, 2,  \cdots, m-1.
 \end{split}\ee
It suffices to prove $\ef{decay}_2$ holds for $k=m$. We derive from \ef{s5.3} that for $m=1, \cdots,k$,
$$
\f{d^{m+1} }{dt^{m+1}}\tilde\eta_{r}(t)+\f{d^{m}}{dt^{m}}\tilde \eta_{r}(t)- \frac{1}{3\ga-1}\frac{d^{m-1}}{dt^{m-1}}\tilde \eta_r^{2-3\ga}(t)=0 , \ \  t\ge 0, $$
so that
\be\label{s5.7}
   \f{d^{m}  }{dt^{m}}  \tilde\eta_{r}(t) = e^{-t} \f{d^{m}}{dt^{m}}\tilde\eta_{r}(0)  + \frac{1}{3\ga-1} \int_0^{t} e^{-(t-s)}\frac{d^{m-1}\tilde \eta_r^{2-3\ga}}{ds^{m-1}}(s)ds,  \  \ t \ge 0,\ee
where $ ({d^{m} / dt^{m})\tilde \eta_{r}}(0) $ is finite, which can be determined by the equation inductively. In view of $\ef{decay}_1$ and \ef{s5.6}, we see that
\bee\label{}
\lt|\frac{d}{dt} \tilde \eta_r^{2-3\ga}(t)\rt|\lesssim \lt|\tilde\eta_r^{1-3\ga}(t) \frac{d}{dt} \tilde \eta_r (t)\rt|\lesssim (1+t)^{\frac{1}{3\ga-1}-2},
  \eee
\bee\label{}
\lt|\frac{d^2}{dt^2} \tilde \eta_r^{2-3\ga}(t)\rt|\lesssim \lt|
 \tilde\eta_r^{ -3\ga}(t)  \lt(\frac{d}{dt} \tilde \eta_r (t)\rt)^2\rt|
+\lt|\tilde\eta_r^{1-3\ga}(t)
   \frac{d^2 }{dt^2} \tilde \eta_r (t) \rt|
\lesssim (1+t)^{\frac{1}{3\ga-1}-3}
  \eee
$$\cdots$$
 \be\label{s5.8}
\lt|\frac{d^{m-1}}{dt^{m-1}} \tilde \eta_r^{2-3\ga}(t)\rt|
\le  C(\ga, m)(1+t)^{\frac{1}{3\ga-1}-m}.
  \ee
Similar to deriving \ef{s5.5}, we can obtain, noting  \ef{s5.7} and \ef{s5.8}, that
\bee\label{}\begin{split}
 \lt|\f{d^m\tilde \eta_{r}(t)}{dt^m}\rt| \le C(\ga, m)\lt(1 +   t \rt)^{\frac{1}{3\ga-1}-m}.
 \end{split}\eee
This finishes the proof of $\ef{decay}_2$. $\Box$

\vskip 0.5cm
\noindent{\bf Proof of \ef{decayforh}}.  We may write the equation for $h$, \ef{pomt}, as
\be\label{s5.9}
h_t + \frac{1}{3\ga-1}(1+t)^{-\frac{3\ga-2}{3\ga-1}} \lt[1-\lt(1+(1+t)^{-\f{1}{3\ga-1}} h\rt)^{2-3\ga}\rt]   =- \tilde\eta_{rtt}  , \ \ t >0\ee
Notice that
$$\lt(1+(1+t)^{-\f{1}{3\ga-1}} h\rt)^{2-3\ga} \le 1+(2-3\ga) (1+t)^{-\f{1}{3\ga-1}} h +\frac{(2-3\ga)(1-3\ga)}{2} (1+t)^{-\f{2}{3\ga-1}} h^2,$$
due to $h\ge 0$. We then obtain, in view of $\ef{decay}_{2}$, that
$$
h_t+\frac{3\ga-2}{3\ga-1}(1+t)^{-1}h\le \frac{  3\ga-2}{2} (1+t)^{-\f{3\ga}{3\ga-1}} h^2 +C (1+t)^{\frac{1}{3\ga-1}-2}.$$
Thus,
\be\label{s5.10}
h(t)\le C(1+t)^{-\frac{3\ga-2}{3\ga-1}} \int_0^t \lt((1+s)^{-\frac{2}{3\ga-1}}h^2(s)+(1+s)^{-1}\rt)ds .\ee
We use an iteration to prove \ef{decayforh}.
First, since $h$ is bounded due to \ef{s5.2}, we have
\be\label{s5.11}h(t)\le C(1+t)^{-\frac{3\ga-2}{3\ga-1}} \int_0^t (1+s)^{-\frac{2}{3\ga-1}}ds \le C(1+t)^{-\frac{1}{3\ga-1}}. \ee
Substituting this into \ef{s5.10}, we obtain
$$h(t) \le C(1+t)^{-\frac{3\ga-2}{3\ga-1}} \int_0^t \lt((1+s)^{-\frac{4}{3\ga-1}}+(1+s)^{-1}\rt)ds ;$$
which implies
\bee\label{}h(t) \le \lt\{ \begin{split} & C (1+t)^{-\frac{3\ga-2}{3\ga-1}}\ln(1+t) & {\rm if} \ \  \ga \le 5/ 3,\\
& C (1+t)^{-\frac{3}{3\ga-1}}  & {\rm if} \ \  \ga > 5/ 3. \end{split}\rt. \eee
If $\ga\le  5/ 3$, then the first part of \ef{decayforh} has been proved. If $\ga>5/3$, we repeat this procedure and obtain
$$h(t) \le C(1+t)^{-\frac{3\ga-2}{3\ga-1}} \int_0^t \lt((1+s)^{-\frac{8}{3\ga-1}}+(1+s)^{-1}\rt)ds ;$$
which implies
\bee\label{}h(t) \le \lt\{ \begin{split} & C (1+t)^{-\frac{3\ga-2}{3\ga-1}}\ln(1+t) & {\rm if} \ \  \ga \le 3,\\
& C (1+t)^{-\frac{7}{3\ga-1}}  & {\rm if} \ \  \ga > 3. \end{split}\rt. \eee
For general $\gamma$, we repeat this procedure $k$ times   to  obtain
$$h(t) \le C (1+t)^{-\frac{3\ga-2}{3\ga-1}}\ln(1+t). $$
This, together with \ef{s5.2}, proves the first part of \ef{decayforh}, which in turn  implies the second part of \ef{decayforh}, by virtue of  \ef{s5.9} and $\ef{decay}_2$. $\Box$

\vskip 0.5cm
\noindent{\bf Proof of \ef{zuihou}}. Recall that $ j\ge 0$, $i\ge 1$ and  $i+j\le l$. Let
$n \in [0, j]$, $m\in [0,i-1]$ and $q\in [0, n]$ be  integers.  Denote
$$\mathcal{H}:= \lt\|  r  \sigma ^{\frac{\aa+i-1}{2}} \lt(\lt|r \pl_t^{n-q} \pl_r^{m+1}\zeta\rt| + \lt|  \pl_t^{n-q} \pl_r^{m }\zeta\rt| \rt) \lt(\lt| \sigma  r \pl_t^{j-n}\pl_r^{i-m +1}  \zeta  \rt|  +  \sum_{\iota=0}^{i-m }  \lt|\pl_t^{j-n}\pl_r^{\iota}  \zeta  \rt|\rt)\rt\|^2.$$

{\bf Case 1}. Assume $2n+4m\ge 2i+j+q$. We first note that
\be\label{s5.12'}
\aa + ( 2m + n ) -(i+j)+2 \ge \aa  -\frac{j}{2} + \frac{q}{2} +2 \ge \aa -\frac{l-1}{2} +2 \ge 0,
\ee
\be\label{s5.12}\begin{split}
    i+j-(n+m)\le l-2.
\end{split}\ee
(Indeed, if $i+j-(n+m)=l$, then $i+j=l$ and $n+m=0$, so that it is a contradiction due to
$$
0=4(n+m) \ge  2n+4m \ge 2i+j+q \ge i+j =l;
$$
if $i+j-(n+m)=l-1$, then $i+j=l-1$ and $n+m=0$ or
$i+j=l$ and $n+m=1$, so that it is also a contradiction because of
$$
0=4(n+m) \ge  2n+4m \ge 2i+j+q \ge i+j =l-1>0
$$
 or
$$4 = 4(n+m) \ge 2n+4m\ge 2i+j+q \ge i+ i+j\ge 1 + l = 5 +[\aa] \ge 5 . $$
So, \ef{s5.12} holds.)

When $ 2i+j\le 2m + n + 3 $, it follows from \ef{apb} and \ef{s5.12} that
\be\label{s5.13}\begin{split}
    \mathcal{H}\lesssim  & \ea_0^2 (1+t)^{2n-2j} \lt\|  r  \sigma ^{\frac{\aa+i-1}{2}} \lt(\lt|r \pl_t^{n-q} \pl_r^{m+1}\zeta\rt| + \lt|  \pl_t^{n-q} \pl_r^{m }\zeta\rt| \rt)  \rt\|^2 \\
   \lesssim &  \ea_0^2 (1+t)^{2q-2j} \lt(\mathcal{E}_{n-q, m+1} +  \mathcal{E}_{n-q, m }\rt).
\end{split}\ee
When $ 2i+j \ge  2m + n + 4 $, it follows from \ef{apb} and \ef{s5.12} that
\bee\label{}\begin{split}
     \mathcal{H}\lesssim  & \ea_0^2 (1+t)^{2n-2j}\lt\|  r  \sigma ^{\frac{\aa+i-1}{2}-\frac{j+2i-(n+2m)-3}{2}} \lt(\lt|r \pl_t^{n-q} \pl_r^{m+1}\zeta\rt| + \lt|  \pl_t^{n-q} \pl_r^{m }\zeta\rt| \rt)  \rt\|^2 \\
   =&\ea_0^2 (1+t)^{2n-2j} \lt\|  r  \sigma ^{\frac{\aa + m +\lt(n+m-(i+j)+2\rt)}{2} } \lt(\lt|r \pl_t^{n-q} \pl_r^{m+1}\zeta\rt| + \lt|  \pl_t^{n-q} \pl_r^{m }\zeta\rt| \rt)  \rt\|^2;
\end{split}\eee
which implies for $n+m-(i+j)+2\ge 0$,
\be\label{s5.15}\begin{split}
    \mathcal{H}\lesssim &\ea_0^2 (1+t)^{2n-2j} \lt\|  r  \sigma ^{\frac{\aa + m  }{2} } \lt(\lt|r \pl_t^{n-q} \pl_r^{m+1}\zeta\rt| + \lt|  \pl_t^{n-q} \pl_r^{m }\zeta\rt| \rt)  \rt\|^2 \\
    \lesssim  & \ea_0^2 (1+t)^{2q-2j} \lt(\mathcal{E}_{n-q, m+1} +  \mathcal{E}_{n-q, m }\rt)
\end{split}\ee
and for $n+m-(i+j)+2 \le -1 $,
\be\label{s5.16}\begin{split}
       \mathcal{H}\lesssim  &  \ea_0^2 (1+t)^{2n-2j} \lt(\lt\|     r^2  \pl_t^{n-q} \pl_r^{m+1}\zeta  \rt\|_{L^2(\mathcal{I}_o)}^2
       + \lt\|    r \pl_t^{n-q} \pl_r^{m }\zeta    \rt\|_{L^2(\mathcal{I}_o)}^2  \rt. \\
       & \lt. + \lt\|     \sigma ^{\frac{\aa + m +\lt(n+m-(i+j)+2\rt)}{2} } \lt(\lt|  \pl_t^{n-q} \pl_r^{m+1}\zeta\rt| + \lt|  \pl_t^{n-q} \pl_r^{m }\zeta\rt| \rt)  \rt\|^2_{L^2(\mathcal{I}_b)}\rt) \\
       \lesssim & \ea_0^2 (1+t)^{2n-2j} \lt(\lt\|     r^2 \sa^{\f{\aa+m+1}{2}}  \pl_t^{n-q} \pl_r^{m+1}\zeta  \rt\|_{L^2(\mathcal{I}_o)}^2
       + \lt\|    r \sa^{\f{\aa+m-1}{2}} \pl_t^{n-q} \pl_r^{m }\zeta    \rt\|_{L^2(\mathcal{I}_o)}^2  \rt)  \\
       +& \ea_0^2 (1+t)^{2q-2j}\sum_{h=m}^{i+j-n+q-1} \mathcal{E}_{n-q, h}
       \lesssim \ea_0^2 (1+t)^{2q-2j}\sum_{h=m}^{i+j-n+q-1} \mathcal{E}_{n-q, h} .
\end{split}\ee
Here we have used \ef{s5.12'} and the Hardy inequality   \ef{hardy2} to derive
\bee\label{}\begin{split}
   &  \lt\|     \sigma ^{\frac{\aa + m +\lt(n+m-(i+j)+2\rt)}{2} } \lt(\lt|  \pl_t^{n-q} \pl_r^{m+1}\zeta\rt| + \lt|  \pl_t^{n-q} \pl_r^{m }\zeta\rt| \rt)  \rt\|^2_{L^2(\mathcal{I}_b)} \\
     \lesssim  & \sum_{h=0}^{1}  \lt\|     \sigma ^{\frac{\aa + m +\lt(n+m-(i+j)+2\rt)}{2} } \pl_t^{n-q} \pl_r^{m+h}\zeta   \rt\|^2_{L^2(\mathcal{I}_b)}
  \lesssim  \sum_{h=0}^{2}  \lt\|     \sigma ^{\frac{\aa + m +\lt(n+m-(i+j)+2\rt) +2 }{2} } \pl_t^{n-q} \pl_r^{m+h}\zeta   \rt\|^2_{L^2(\mathcal{I}_b)}  \\
       \lesssim  & \cdots \lesssim  \sum_{h=0}^{i+j-(n+m)+q }  \lt\|     \sigma ^{\frac{\aa + m +\lt(n+m-(i+j)+2\rt) +2(i+j-(n+m)+q-1) }{2} } \pl_t^{n-q} \pl_r^{m+h}\zeta   \rt\|^2_{L^2(\mathcal{I}_b)} \\
       = & \sum_{h=0}^{i+j-(n+m)+q }  \lt\|     \sigma ^{\frac{\aa + i+j -n +2q  }{2} } \pl_t^{n-q} \pl_r^{m+h}\zeta   \rt\|^2_{L^2(\mathcal{I}_b)} \\
       \lesssim & \sum_{h=0}^{i+j-(n+m)+q }  \lt\|  r^2   \sigma ^{\frac{\aa + i+j -n + q  }{2} } \pl_t^{n-q} \pl_r^{m+h}\zeta   \rt\|^2_{L^2(\mathcal{I}_b)}
       \lesssim (1+t)^{2q-2n}\sum_{h=m}^{i+j-n+q-1} \mathcal{E}_{n-q, h},
      \end{split}\eee
which implies \ef{s5.16}. Therefore, it gives from \ef{s5.13}, \ef{s5.15} and \ef{s5.16} that
\be\label{s5.17}\begin{split}
       \mathcal{H}
       \lesssim \ea_0^2 (1+t)^{2q-2j} \lt(  \mathcal{E}_{n-q, m }+ \mathcal{E}_{n-q, m+1} + \sum_{h=m}^{i+j-n+q-1} \mathcal{E}_{n-q, h}  \rt).
\end{split}\ee

{\bf Case 2}. Assume $2n+4m < 2i+j+q$. In this case, we can use the similar way to dealing with Case 1 to obtain
\be\label{s5.18}\begin{split}
       \mathcal{H}
       \lesssim \ea_0^2 (1+t)^{2q-2j}  \ \sum_{ h=0 }^{i+n-j} \mathcal{E}_{j-n, h } (t).
\end{split}\ee

In view of \ef{s5.17} and \ef{s5.18}, we prove \ef{zuihou}. $\Box$

\vskip 0.5cm

\noindent{\bf Acknowledgements}  The author is grateful to Professor Tao Luo for his many helpful discussions. This work was supported in part by NSFC under grant  11301293.

\bibliographystyle{plain}

\noindent Huihui Zeng\\
Mathematical Sciences Center\\
Tsinghua University\\
Beijing, 100084, China;\\
Center of Mathematical Sciences and Applications \\
Harvard University \\
Cambridge, MA 02318, USA.\\
E-mail: hhzeng@mail.tsinghua.edu.cn

\end{document}